\font\smallrm=cmr8

\font\script=rsfs10

\def\pitt#1{\mathop{\vtop{\ialign{##\crcr$\hfil\displaystyle
{#1}\hfil$\crcr\noalign{\kern1pt\nointerlineskip}\crcr\noalign
{\kern-2pt}}}}\limits}

\long\def\symbolfootnote[#1]#2{\begingroup%
\def\thefootnote{\fnsymbol{footnote}}\footnote[#1]{#2}\endgroup} 

\documentclass[11pt]{article}
\usepackage[english]{babel}
\usepackage{epsfig}
\usepackage{amsmath, amsthm, amsfonts}
\numberwithin{equation}{section}
\usepackage{graphicx}
\usepackage{rawfonts}
\usepackage{pslatex}
\usepackage{float}
\usepackage[textwidth=14cm,textheight=22cm]{geometry}

\begin{document}

\vglue 1cm

\centerline{\LARGE\bf On ringing effects near jump discontinuities}
\vskip 1mm
\centerline{\LARGE\bf for periodic solutions to dispersive partial}
\vskip 1mm
\centerline{\LARGE\bf differential equations} 

\vskip 4mm 

\centerline {
B\scalebox{0.7}{Y} K\scalebox{0.7}{ENNETH} D. T.-R. M\scalebox{0.7}{C}L\scalebox{0.7}{AUGHLIN}${}^{1,2}$
\scalebox{0.7}{AND} N\scalebox{0.7}{IGEL} J. E. P\scalebox{0.7}{ITT}${}^{2,3,}$
\!\!\symbolfootnote[1]{\small Author for correspondence (pitt@mat.unb.br)}}

\vskip 2mm

\centerline {\it ${}^1$Department of Mathematics, University of Arizona, Tucson, AZ 85721, USA}

\centerline {\it ${}^2$Departamento de Matem\'atica, Universidade de Bras\'{\i}lia, DF 70910-900, Brazil}

\centerline {\it ${}^3$School of Mathematics, Institute for Advanced Study, Princeton, NJ 08540, USA}

\vskip 4mm

\vbox{
\baselineskip=11pt
\noindent {\small 
We consider weak solutions to dispersive partial differential equations with 
periodic boundary conditions and initial data with jump discontinuities.  
These are already known to be continuous at irrational times and piecewise 
constant at rational times; we show that as time approaches a rational value 
the solution exhibits a ringing effect, with the characteristic overshoot of 
fixed amplitude near the discontinuities.  Furthermore this effect is the 
same whether the sequence of times follows rational or irrational values.}}

\vskip 1mm

\centerline{\small\bf Keywords: ringing effect; dispersive equations; 
weak solutions, periodic boundary conditions;}
\centerline{\small\bf Diophantine approximations; Weyl shift method}

\vskip 2mm

\hrule 

\vskip 6mm

\setcounter{section}{1}

\centerline{\bf 1. Introduction}

\vskip 1mm

Consider the equation 
\begin{equation}
U^{(0,n)}(t,x)=(2\pi i)^{n-1}U^{(1,0)}(t,x)
\end{equation}
for $n\ge 2$.  Any classical solution to this equation satisfies the 
condition that
\begin{equation}
\pitt{\int{\hskip -2mm}\int}_{{\mathbb R}^2} U(t,x)\left\lbrace
F^{(0,n)}(t,x)-(-2\pi i)^{n-1}F^{(1,0)}(t,x)\right\rbrace dx\,dt
= 0
\end{equation}
for any ${\hbox{\script C}}^{\,\,\infty}$ function $F$ of compact 
support in ${\mathbb R}^2$, as can be seen using integration by parts, 
and we call a function $U(t,x)$ satisfying $(1.2)$ a weak
solution to $(1.1)$.  We are interested here in weak solutions with periodic boundary conditions
$U(t,x+1)=U(t,x)$ for all $x$.  It is natural to study such 
solutions through their representations by Fourier series, which 
are of course ubiquitous in mathematics, and we will use techniques from 
analytic number theory to analyse these series and show some asymptotic 
properties of the weak solution $U(t,x)$.

To motivate this, for the present consider the simpler case of $(1.1)$ 
for $n=2$ with vanishing boundary conditions $U(t,x)\rightarrow 0$ as 
$x\rightarrow\pm\infty$, and initial data with jump discontinuities; 
for instance $U(0,x)=\chi(x)$, where $\chi$ is the 
characteristic function of the interval $[-\gamma,\gamma]$ for 
$0<\gamma < 1/2$.   By considering the Fourier transform in $x$ of 
$U(t,x)$ one can show that
\begin{equation}
U(t,x)=\int\limits_{-\infty}^{\infty}
{\hat U}(0,k)e(tk^2+xk)dk
=\int\limits_{-\infty}^{\infty}\frac{\sin 2\pi \gamma k}{\pi k}
e(tk^2+xk)dk
\end{equation}
where we have used the notational 
convention $e(\xi)=e^{2\pi i\xi}$, as throughout this paper. With $t=0$, 
the integral represents the initial data $\chi(x)$ in $L^{2}$.  
Moreover it is not hard to prove that the truncated integral
$$
U(0,x;N):= \int_{-N}^{N} \frac{ \sin{ 2 \pi \gamma k}}{\pi k} e\left( 
x k \right) dk
$$
converges pointwise to $\chi(x)$.  Of course, the convergence cannot be 
uniform and indeed one has
$$
\lim_{N\to \infty}U\left(0, \gamma + s/N;N \right) = \frac{1}{2} - \frac{1}{\pi} \int_{0}^{s} \frac{\sin{2 \pi w}}{w} dw.
$$
This is the well-known Gibbs phenomenon, which is an oscillatory ringing 
produced by truncation of the Fourier representation of a 
function with a discontinuity.  

There is, however, an entirely different phenomenon, 
which is an oscillatory ringing produced not by truncation but 
rather by a different mechanism 
which, as we will explain, should be thought of as a dispersive 
regularization of a discontinuity.  Indeed, by contour deformation or 
otherwise, the integral solution (1.3) can be seen to be 
analytic in $t$ or $x$ if $t>0$.  Away from $x=\pm \gamma$ 
the integral converges nicely to $U(0,x)$ as $t\rightarrow 0^+$, however 
in the rescaled variable $s$ where $x=\pm\gamma +s t^{1/2}$ one can 
obtain the asymptotic expression (see DiFranco \& McLaughlin 2005) 
\begin{equation}
U\left(t,\gamma +st^{1/2}\right)
=\frac{1}{2}-\frac{1}{2} \,{\hbox{\rm Erf}}\,\left({\sqrt{\frac{\pi}{2}}} e^{-i\pi/4}s\right) +O\left(t^{1/2}\right)
\end{equation}
where ${\hbox{\rm Erf}}$ denotes the usual error function (see Abramovitz 
\& Stegun 1972).  Similar phenomena can be shown to occur for $n>2$, 
although with different ringing functions.  

Our goal here is to study 
similar effects in the more complicated situation where $U(t,x)$ obeys 
periodic boundary conditions, which have not been observed previously, 
and to give asymptotic expressions for these effects.  Here Fourier 
series expansions replace the Fourier transform (or if one prefers, 
series replace integrals), and the resulting structure is more 
complicated, presenting ringing effects at all rational times, not 
just at zero as in the ``whole line'' case.  The general setup is as follows.  Suppose now that $U(t,x)$ is 
defined on $[0,\delta]\times {\mathbb R}$ for some $\delta > 0$ and is 
periodic in $x$, that is, $U(t,x)=U(t,x+1)$ for all $t$ and $x$.  
Such a function has a Fourier series 
\begin{equation}
U(t,x) \sim \sum_{k=-\infty}^{\infty} c_k(t)e(kx)\ ,\ 
c_k(t)= \int\limits_0^1 U(t,\xi)e(-k\xi)d\xi.
\end{equation}
If $U(t,x)$ is a square-integrable weak solution to $(1.1)$ then 
standard arguments show that $c_k(t)=c_ke(tk^n)$ where the constants 
$c_k$ are the Fourier coefficients 
of $U(0,x)$, and where the convergence of the series to $U(t,x)$ 
is understood in the $L^2$ sense.  
Note that periodicity in $x$ has forced periodicity in $t$ also, which is
the root cause of the special behaviour at rational times.
We will say a function $f(x)$ is in class ${\hbox{\script D}}$ if it is
integrable, periodic of period 1, piecewise continuously differentiable, and 
$f(x)=\left\lbrace f\left(x^+\right)+f\left(x^-\right)\right\rbrace /2$
for all $x$.  It is a well-known theorem of harmonic analysis (see 
Katznelson 2004, for instance) that 
functions $f$ of class ${\hbox{\script D}}$ have Fourier series that 
converge pointwise to $f$, in the sense that
$$f(x)=\lim_{K\rightarrow\infty} \sum_{|k|\le K}c_k e(kx)\ ,$$
that is, the series converges pointwise once written in terms of 
sines and cosines rather than exponentials.  (The series of exponentials 
may diverge, in the literal sense).  We will consider periodic 
discontinuous initial conditions
\begin{equation}
U(0,x)=\sum_{l=-\infty}^{\infty}\chi_{[-\gamma,\gamma]}(x+l)
\end{equation}
where $0<\gamma<1/2$, which is in class ${\hbox{\script D}}$, 
with $c_k=(\pi k)^{-1}\sin 2\pi \gamma k$, so the Fourier series $(1.5)$ 
for $U(t,x)$ is
\begin{equation}
U(t,x)\sim 
\sum_k \frac{\sin 2\pi \gamma k}{\pi k} e(tk^n+xk)\ .
\end{equation}
The series clearly converges in the $L^2$ sense, and represents the unique 
weak solution in that space.  We will show some structural properties of 
$U(t,x)$ by identifiying it with this series representation and
considering the properties of the series.

To understand why the situation should be so much more complicated
requires reviewing some known results.  In a variety of applied sciences, 
the term ``dispersive'' describes a localized quantity which spreads 
out as time passes.  In the context of wave phenomena, an equation is 
called dispersive if its simple oscillatory solutions propagate at 
velocities that depend on their spatial frequencies.
This is the case here; rewriting the coefficients as 
$e(k ( x + k^{n-1} t))$ they can be seen as a family 
of traveling waves with velocity $k^{n-1}$, which is clearly dependent on
the spatial frequency $k$.  Thus if one starts with localised 
initial data the various terms of the Fourier series of the solution 
spread out in space, since they move at different velocities.  As time 
passes the solution will appear more extended, with the 
highest spatial frequency components most apparent in the furthest reaches 
of the solution.  Properly speaking, dispersion requires sufficient 
extent to permit the observation of spreading.  This may, or may 
not, be the case for a partial differential equation with periodic 
boundary conditions, but nonetheless such an equation is still called 
dispersive if its family of solutions have frequency dependent velocity.
In this context it is striking that periodic boundary conditions should 
cause a dispersive equation like (1.1) to show recurrence, as well as 
radically different behaviour at rational and irrational times.
This can be seen experimentally in the case $n=2$, as reported by
Talbot 1836, who shone light on a periodic grating and looked at the images
produced by it.  He reports {\it ``...a regular alternation of numerous
lines or bands of red and green colour, having their direction parallel 
to the lines of the grating.  On removing the lens a little further from 
the grating, the bands gradually changed their colours, and became 
alternately blue and yellow.''}  (Talbot 1836).   This is now known as the 
Talbot effect, and has been extensively studied; we refer the reader to 
Berry \& Klein 1996 for details of the effect, but also
for two contributions which are pertinent to our discussion here.
Supposing that the period of the grating is $a$ and that the wavelength 
of the light is $\lambda$, the images seen by Talbot form at regular 
multiples of $a^2/\lambda$.  (In the context of $(1.1)$,
we can view Talbot's grating as initial data which is periodic, being
$1$ on the slits of the grating and zero elsewhere).  
Berry and Klein observe that at 
rational multiples $(p/q)a^2/\lambda$ of this distance {\it ``...these 
fractional Talbot images consist of $q$ equally spaced copies of the transmission function of the 
grating, which superpose coherently when they add up.''}.
Furthermore, they show that these translates have phases given 
by Gauss sums, so the solutions can be written down explicitly 
at rational times, and are piecewise constant functions of $x$.  
(This was also discussed in Olver 2010).  They also observe that 
this is in sharp contrast with irrational times; considering 
the solution at a sequence of rational times tending to an irrational $t$, 
they show that 
{\it ``the graph of a function with power spectrum 
$|g_{n}|^{2}$ proportional to $n^{-\beta}$ is a fractal curve with fractal 
dimension $D=(5- \beta)/2$.  A smooth curve has $D = 1$, and a 
curve with $D = 2$ is so jagged that it is almost area filling; 
curves with $1 < D < 2$ are continuous but non-differentiable.
Thus, since the fourier coefficients of the initial data have 
power spectrum decaying like $n^{-2}$, the fractal dimension must be $3/2$.''}

A number of researchers from a variety of different areas of analysis, 
including Arkhipov \& Oskolkov 1989, Stein \& Wainger 1990, Oskolkov 1992, 
Kapitanski \& Rodnianski 1999, Rodnianski 1999 and Rodnianski 2000, put 
these ideas on a rigorous footing by establishing 
convergence properties of the associated Fourier series in the space 
of continuous functions.  One result which is particularly pertinent to
the discussion here is due to Rodnianski 2000: this discusses 
periodic solutions to $(1.1)$ with $n=2$ with initial data of bounded 
variation but not in the space $\cup_{\epsilon>0} H^{1/2+\varepsilon}$.
(Here the Sobolev space $H^{s}$ is the space of functions whose 
Fourier coefficients $\hat{f}(k)$ satisfy 
$ \sum_{k=1}^{\infty} | \hat{f}(k)|^{2} ( 1 + |k|^{2})^{s/2} 
\ < \infty$.)  The result states that for almost all irrational $t$,
including all algebraic $t$, the fractal dimension of the graphs of 
the real and imaginary parts of $U(t,x)$ 
(which are continuous functions of $x$) is 3/2.

In fact there is a dichotomy between rational and irrational times 
for all $n\ge 2$, which is reflected in our calculations.  
In contrast to the considerations in Berry \& Klein 1996 or 
Rodnianski 2000, which consider sequences of rational times tending
to an irrational limit, to understand the nature of the function
at irrational times, we consider sequences of times, rational or irrational, 
tending to a rational limit, to understand the ringing phenomena exhibited.

As observed above $(1.7)$ converges in $L^2$ to $U(t,x)$, however 
convergence in any stronger sense is more complicated, 
depending greatly for the reasons mentioned above on whether $t$ is rational or 
irrational.  For irrational times $t$ we will use rational approximations;
we speak of approximants as being reduced fractions $u/q$ such that 
\begin{equation}
\left| t-\frac{u}{q}\right| <  \frac{1}{q^2}\ .
\end{equation}
There are infinitely many approximants to any irrational $t$, 
some of them given by the convergents $u_j/q_j\in {\mathbb Q}$ from the 
continued fraction expansion to $t$, which obey the inequality
\begin{equation}
\left| t-\frac{u_j}{q_j}\right| <  \frac{1}{q_jq_{j+1}}\ .
\end{equation}
(See Hua 1982 or Hardy \& Wright 1979).  Note also that a celebrated
theorem in Roth 1955 states that if $t$ is an algebraic irrational
then for any $\eta > 0$ 
there are only finitely many solutions to 
\begin{equation}
\left|t-\frac{u}{q}\right| < \frac{1}{q^{2+\eta}}.
\end{equation}

We need to introduce some sets of irrationals described by their rational
approximations.  To this end $n\ge 2$ will always represent 
the order of the equation $(1.1)$, and $\Delta$ will be fixed once and for 
all in $(0,1)$, and in $(0,1/2)$ in the case $n=2$.  Note that several 
of the sets and implied constants below will depend on $n$ and $\Delta$ 
without this being explicitly mentioned or shown in the notation.

\vskip 2mm

{\bf Definition 1.1.} 
{\it Define 
${\hbox{\script A}}$ to be the set of $t\in [0,1]$ such that 
for all sufficiently large $M$ there is an approximant $u/q$ to $t$ with 
\begin{equation}
M^{\Delta} < q\le M^{n-\Delta}
\end{equation}
and define ${\hbox{\script A}}_{m}$ to be the set of 
$t\in [0,1]$ such that for all $M\ge m$ there is an approximant 
$u/q$ to $t$ satisfying $(1.11)$.  Further, for 
$\alpha>0$ define ${\hbox{\script B}}_{m,\alpha}$ to be the set of 
$t\in [0,1]$ such that for all $M\ge m t^{-\alpha}$ there is an
approximant $u/q$ to $t$ satisfying $(1.11)$. }

\vskip 2mm

Although some parts are not strictly necessary below, we prove the 
following lemma.

\vskip 2mm

{\bf Lemma 1.2.} {\it Let $\mu$ denote Lebesgue measure, and 
suppose $(n-\Delta)^{-1}<\alpha< \Delta^{-1}$.

\noindent {\rm 1.} If $m_1\le m_2$ then ${\hbox{\script A}}_{m_1} 
\subseteq {\hbox{\script A}}_{m_2}$.  

\noindent {\rm 2.} If $t\in {\hbox{\script A}}$ then there 
exists $m$ such that $t\in {\hbox{\script A}}_{m}$, so 
${\hbox{\script A}} = \bigcup_{m \ge 2} {\hbox{\script A}}_{m}$.

\noindent {\rm 3.} There is a constant $c_{n,\Delta}$ such that 
$\mu\left({\hbox{\script A}}_{m}\right)
\ge 1-c_{n,\Delta}m^{1+2\Delta-n}$, 
which tends to $1$ as $m\rightarrow \infty$, so 
$\mu\left({\hbox{\script A}}\right)=1$.   

\noindent {\rm 4.} ${\hbox{\script A}}$ contains all 
algebraic irrationals.  

\noindent {\rm 5.} For all $m$, ${\hbox{\script A}}_{m}\subseteq 
{\hbox{\script B}}_{m,\alpha}$, so  
$\mu\left({\hbox{\script B}}_{m,\alpha}\right)\rightarrow 1$ 
as $m\rightarrow\infty$.  

\noindent {\rm 6.} There is a constant $c_{n,\Delta,\alpha}$ such that 
$$\mu\left( {\hbox{\script B}}_{m,\alpha}
\cap [0,t_0]\right)
\ge t_0 \left(1-c_{n,\Delta,\alpha}t_0^{(n-1-2\Delta)\alpha}
m^{-n+\Delta +1}\right).$$}

\vskip 2mm

The convergence of the series $(1.7)$ is described in the following result.

\vskip 2mm

{\bf Theorem 1.3.} {\it Let
\begin{equation}
S_K=\sum_{0<|k|\le K} \frac{e(tk^n+xk)}{k}\ .
\end{equation}

\noindent {\rm 1}. If $t$ is rational then $S_K$ converges pointwise 
in $x$ to a piecewise constant function.  

\noindent {\rm 2}. If $t\in {\hbox{\script A}}$
then the sequence $S_K$ converges uniformly in $x$. }
 
\vskip 2mm

Note that Arkhipov and Oskolkov 1989 have shown that a more
general class of series converges pointwise, and hence the partial
sums $(1.12)$ converge pointwise in $x$ for all times.  Their proof
uses Vinogradov's method for exponential sums in place of the Weyl 
shift method described below and used here; we state and prove Theorem 1.3
since it is necessary to Theorem 1.5 below, for which we are as yet 
unable to apply the Vinogradov method.

\begin{figure}[h]\label{q=7,real}
\begin{center}
\scalebox{0.2}{\includegraphics{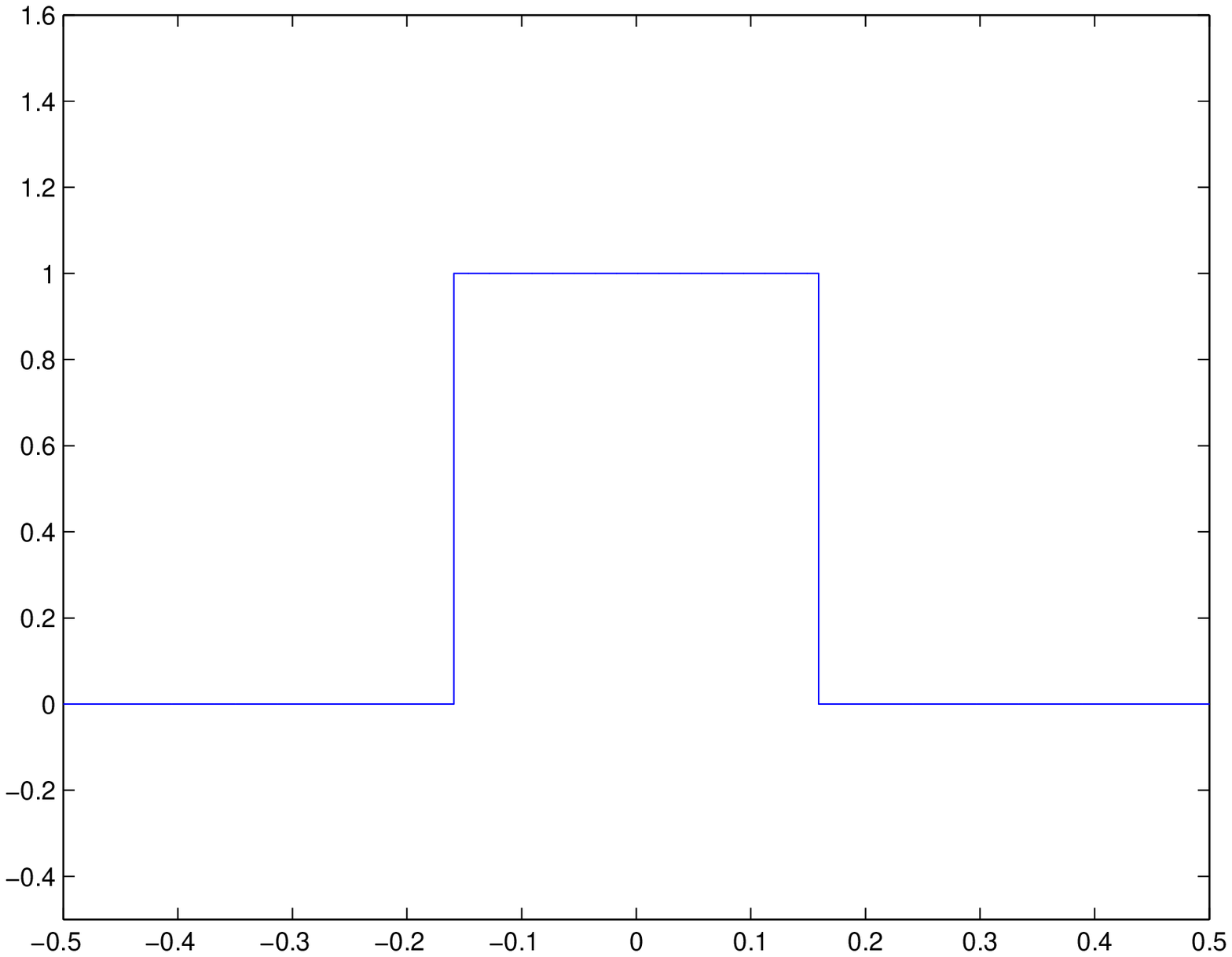}
\includegraphics{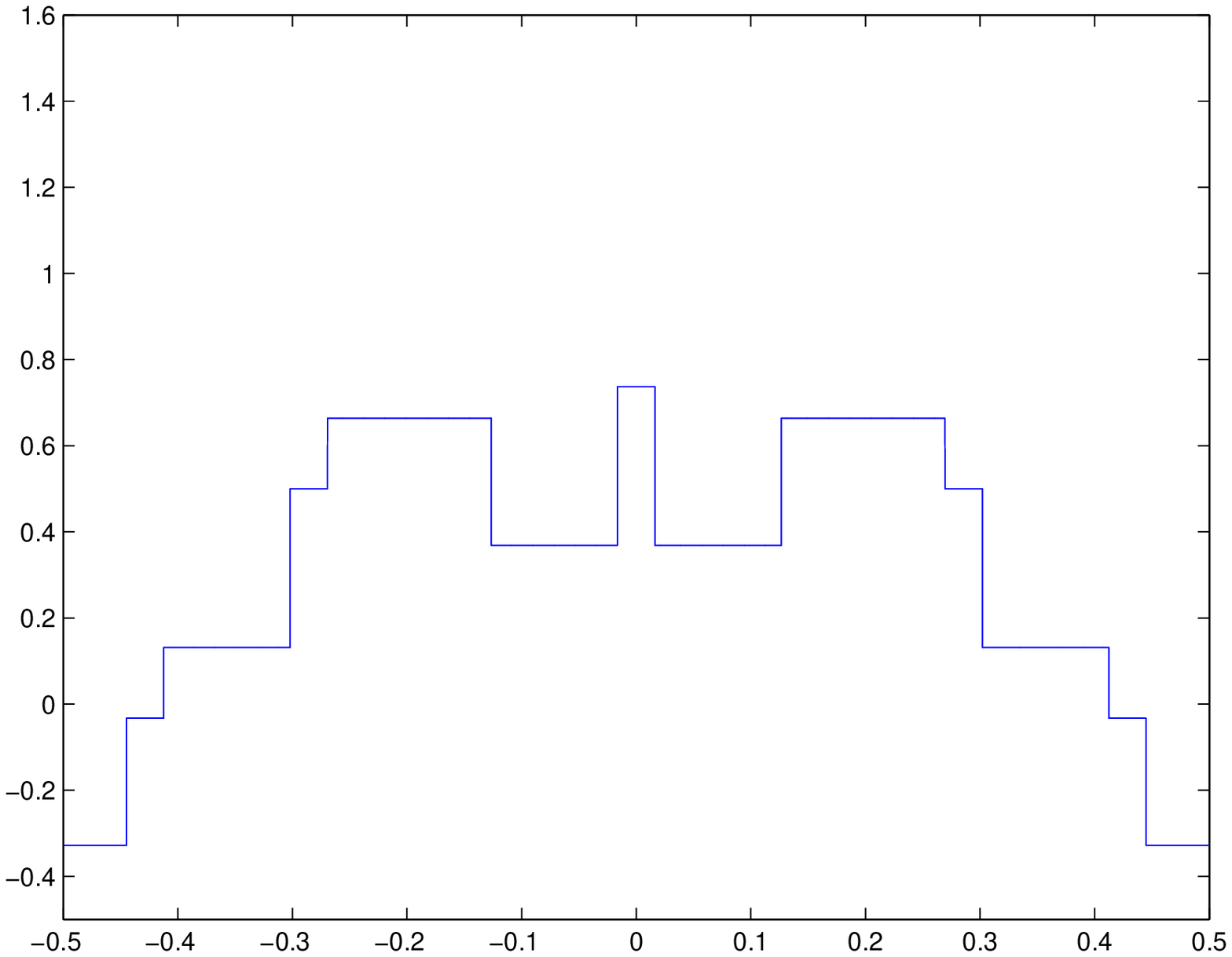}
\includegraphics{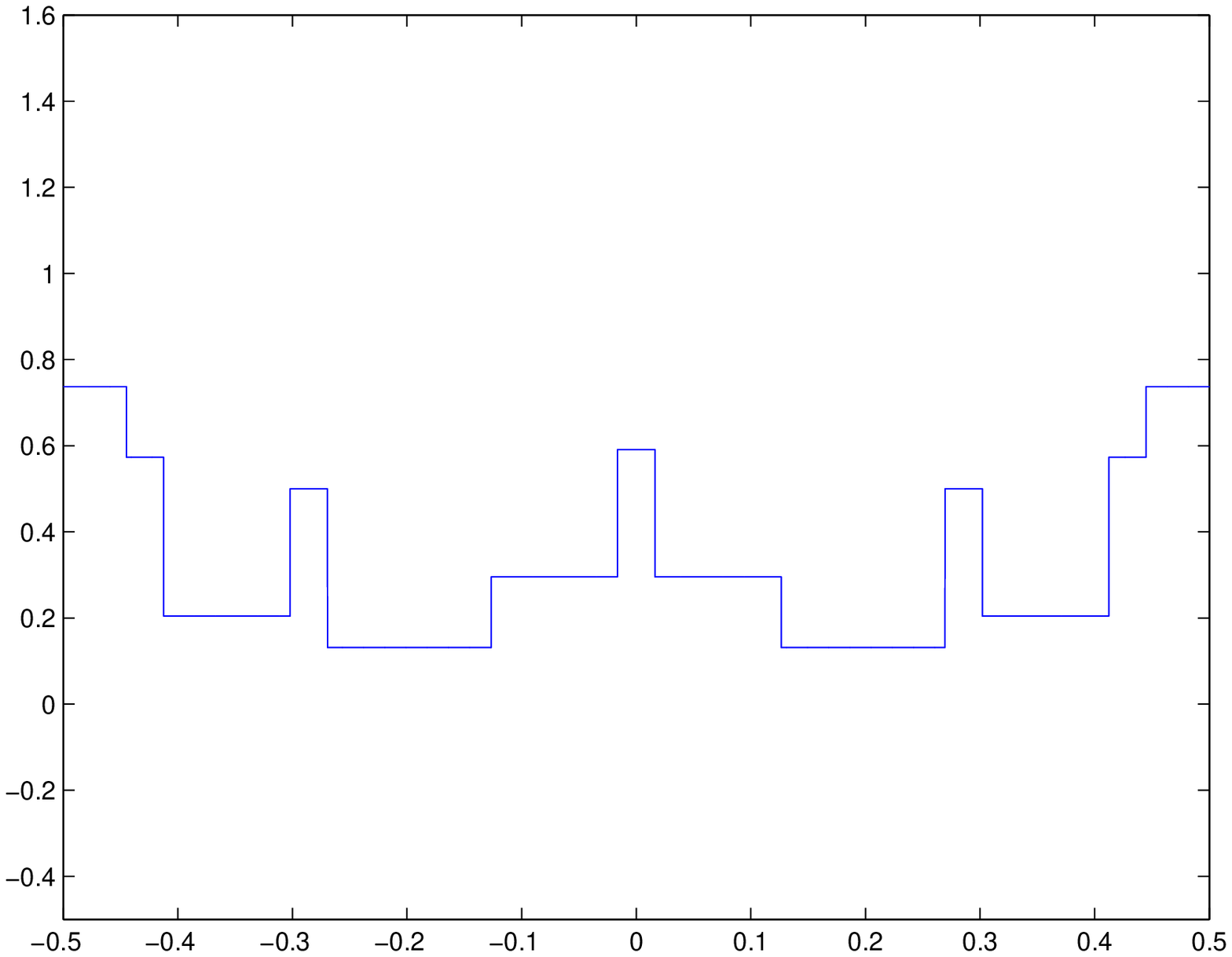}
\includegraphics{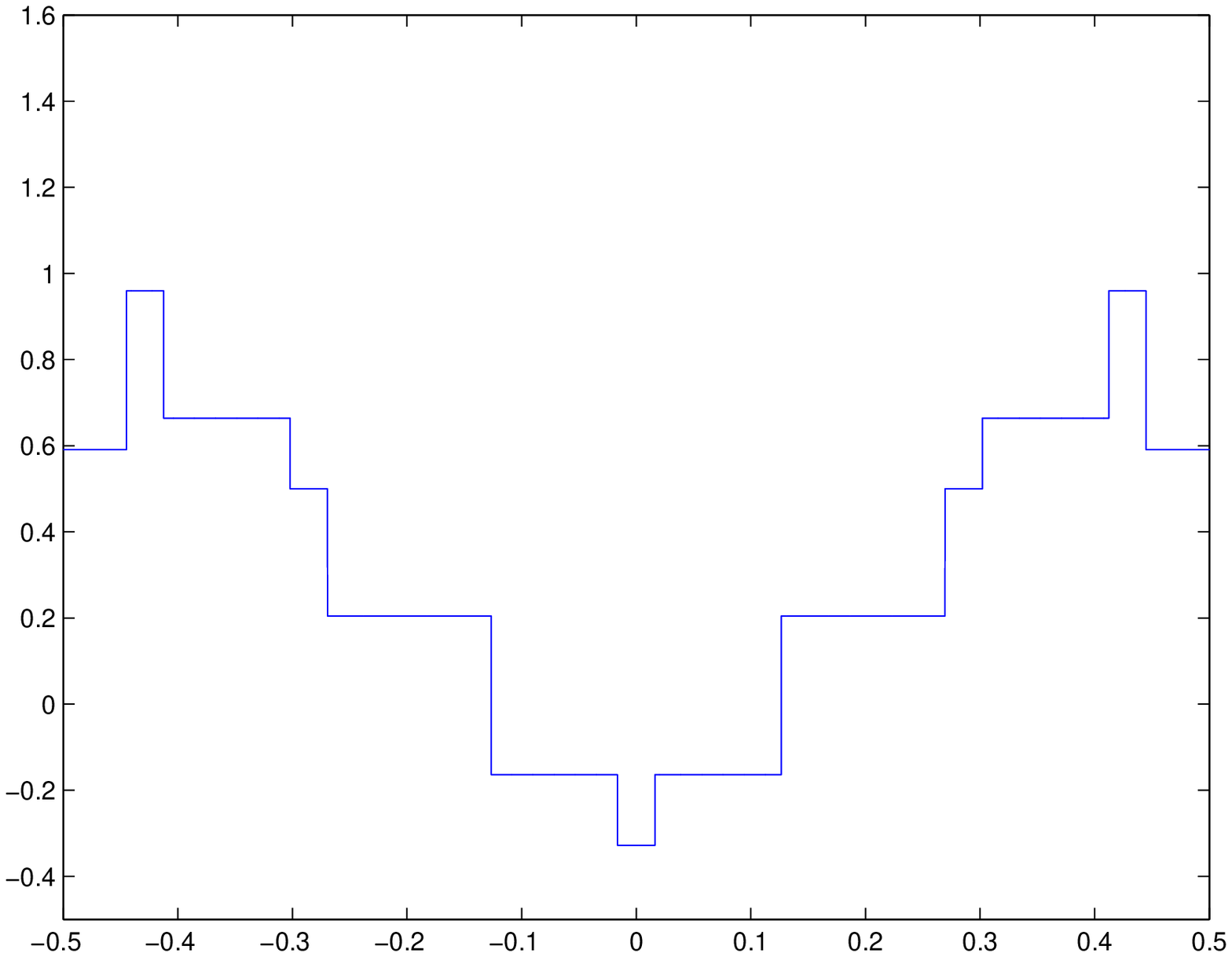}}
\end{center}
\begin{center}
\scalebox{0.2}{\includegraphics{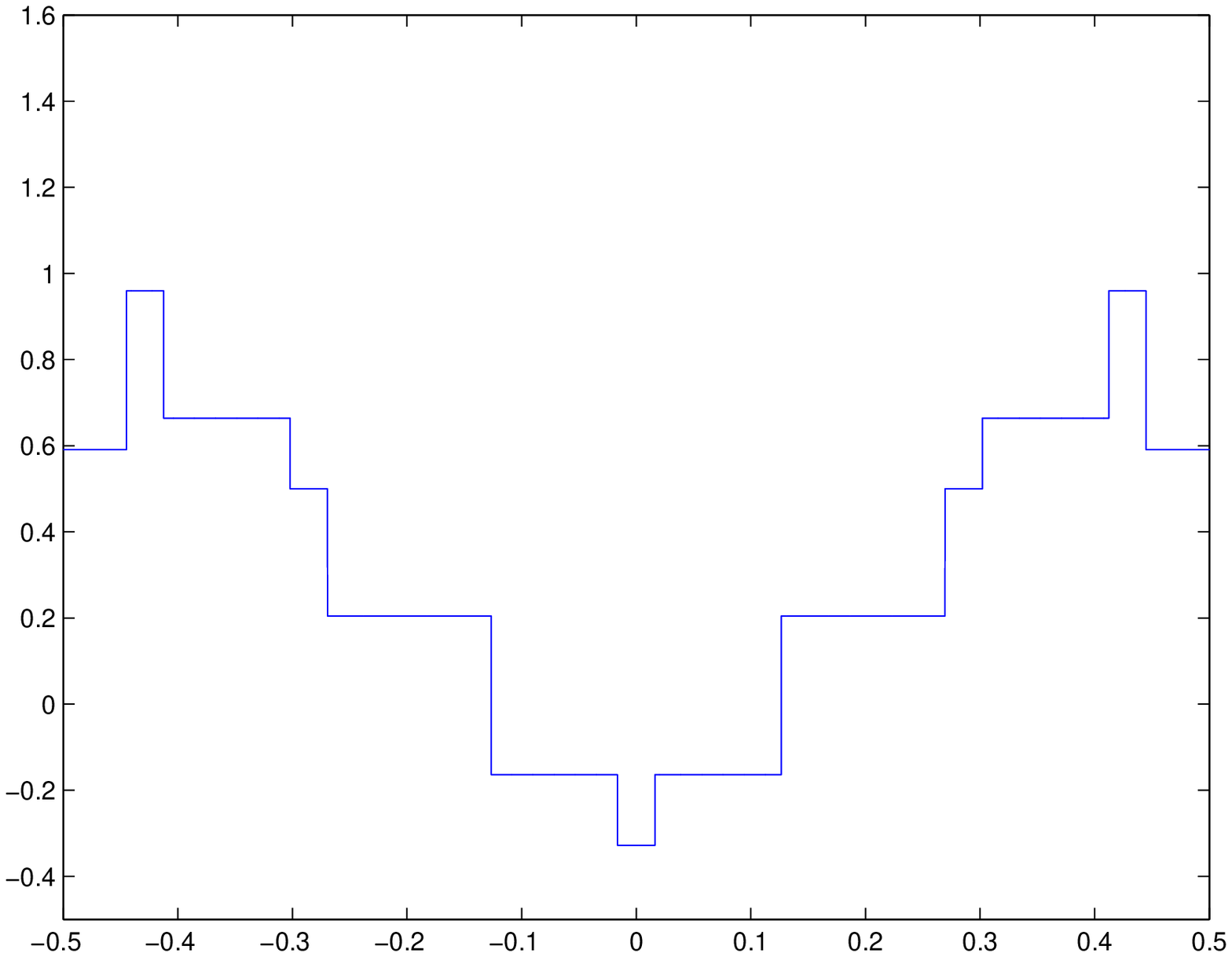}
\includegraphics{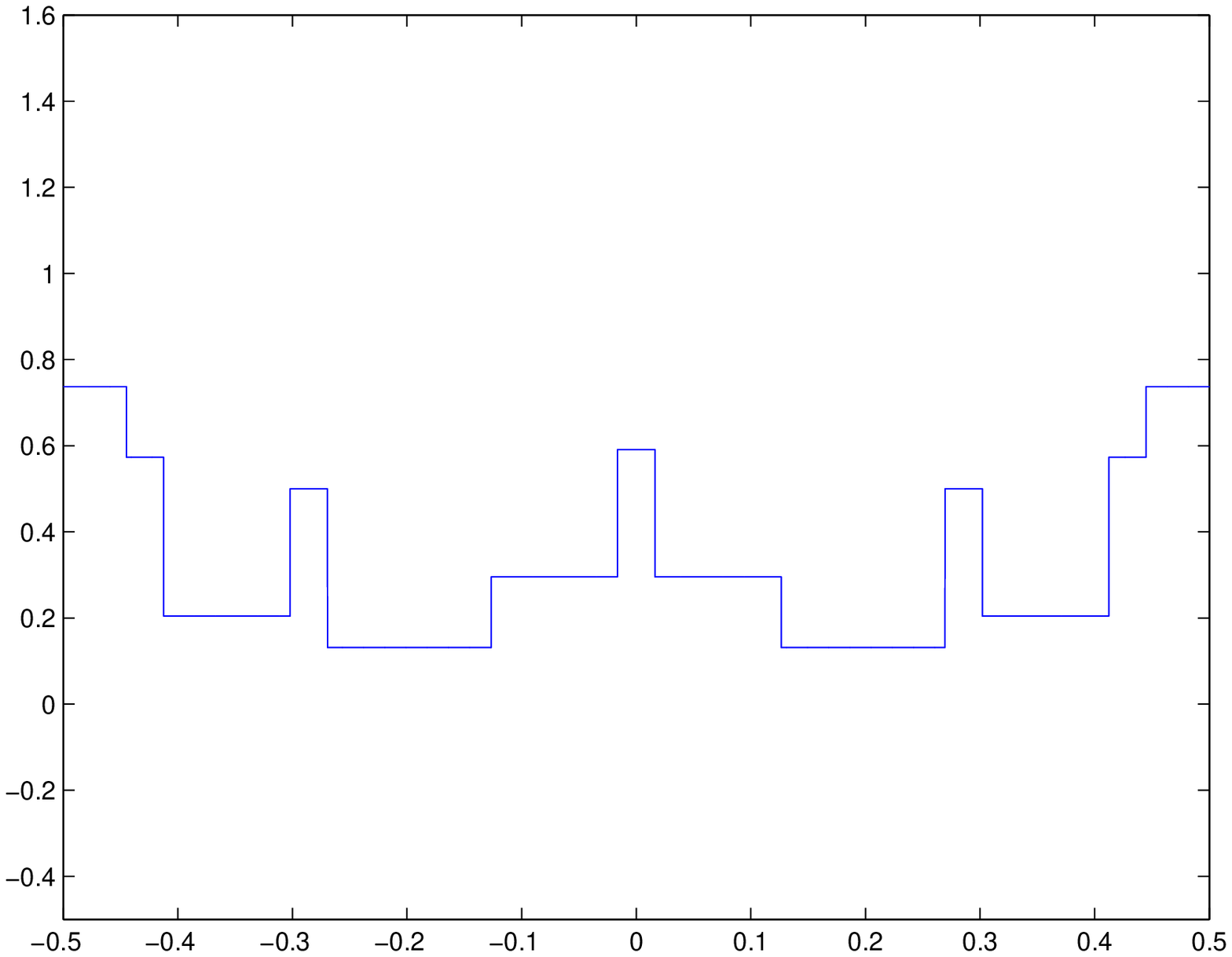}
\includegraphics{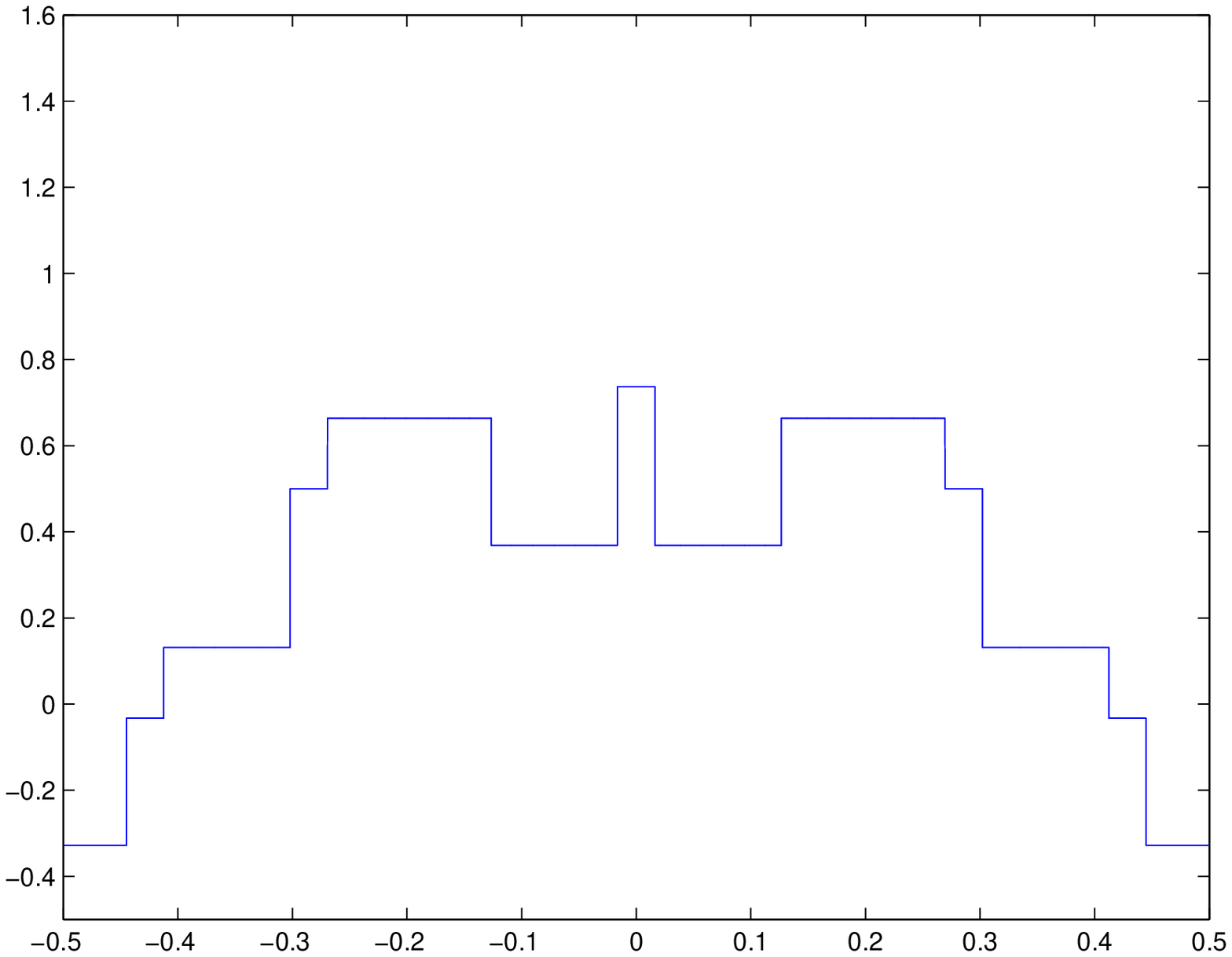}
\includegraphics{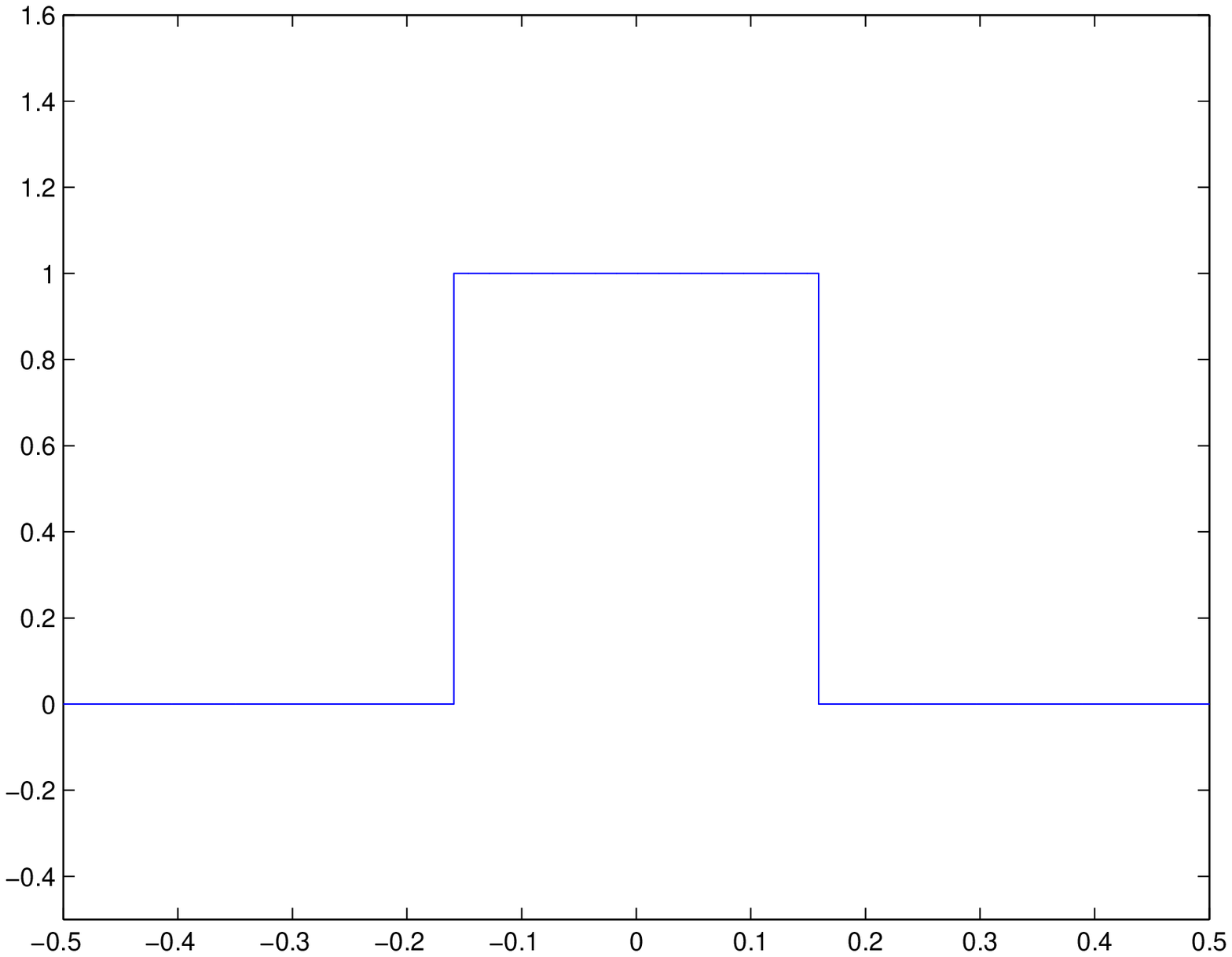}}
\end{center}
\vskip -5mm
\caption{The real part of $U(t,x)$ for $n=2$
at $u/7$ for $u=0,1,2,3$ and $u=4,5,6,7$}
\end{figure}

\begin{center}
\begin{figure}[h]\label{q=7,imag}
\vskip 3mm
\begin{center}
\scalebox{0.2}{\includegraphics{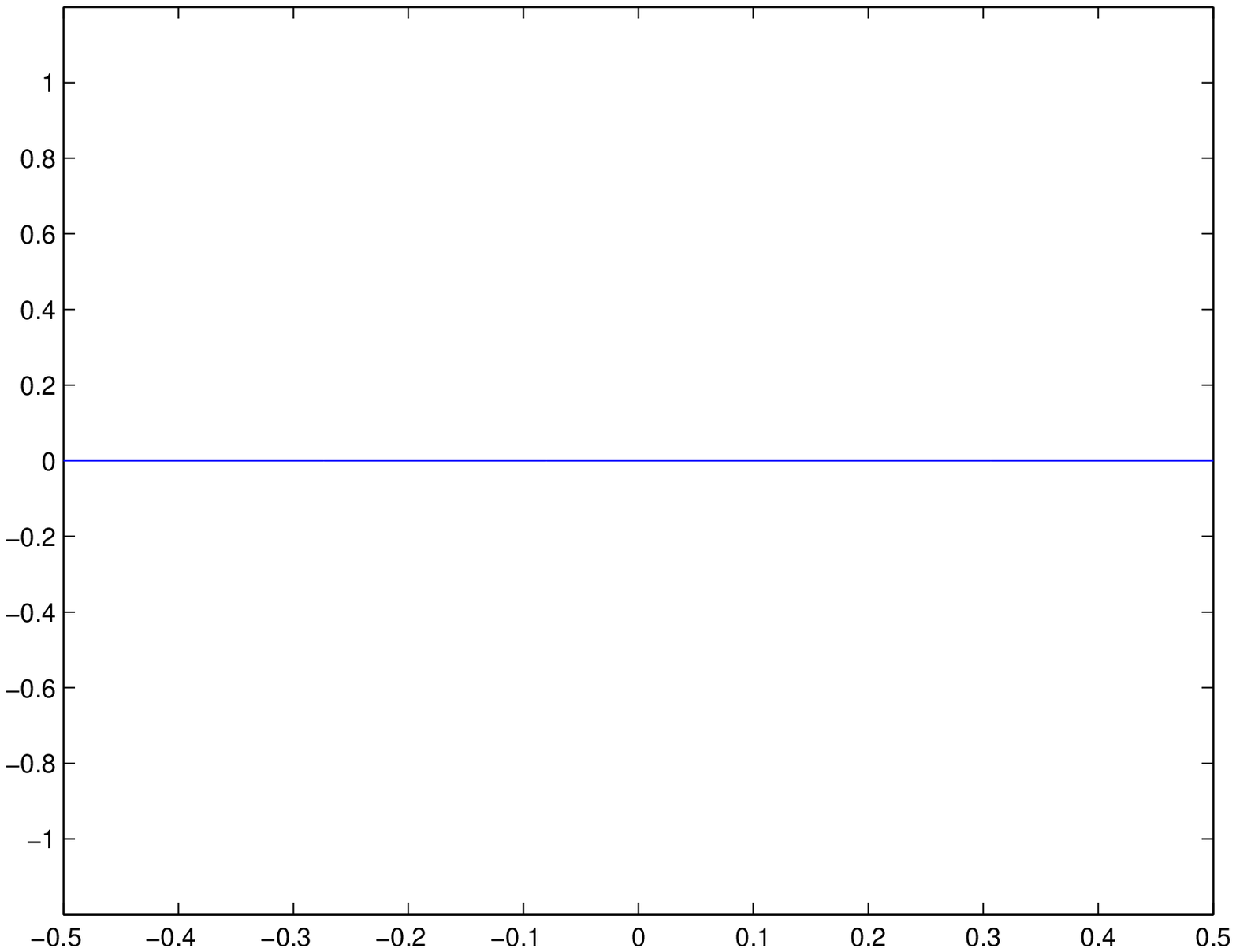}
\includegraphics{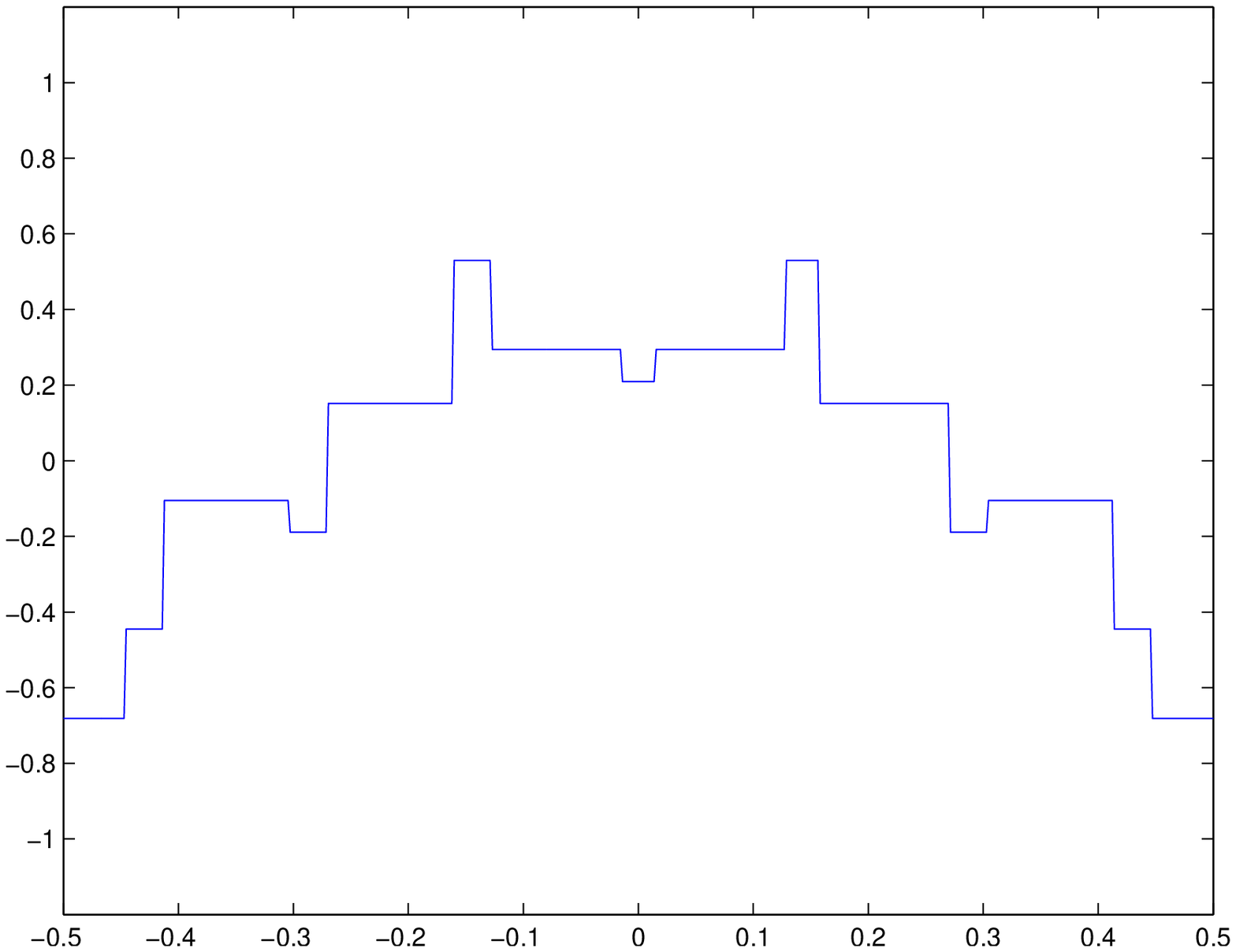}
\includegraphics{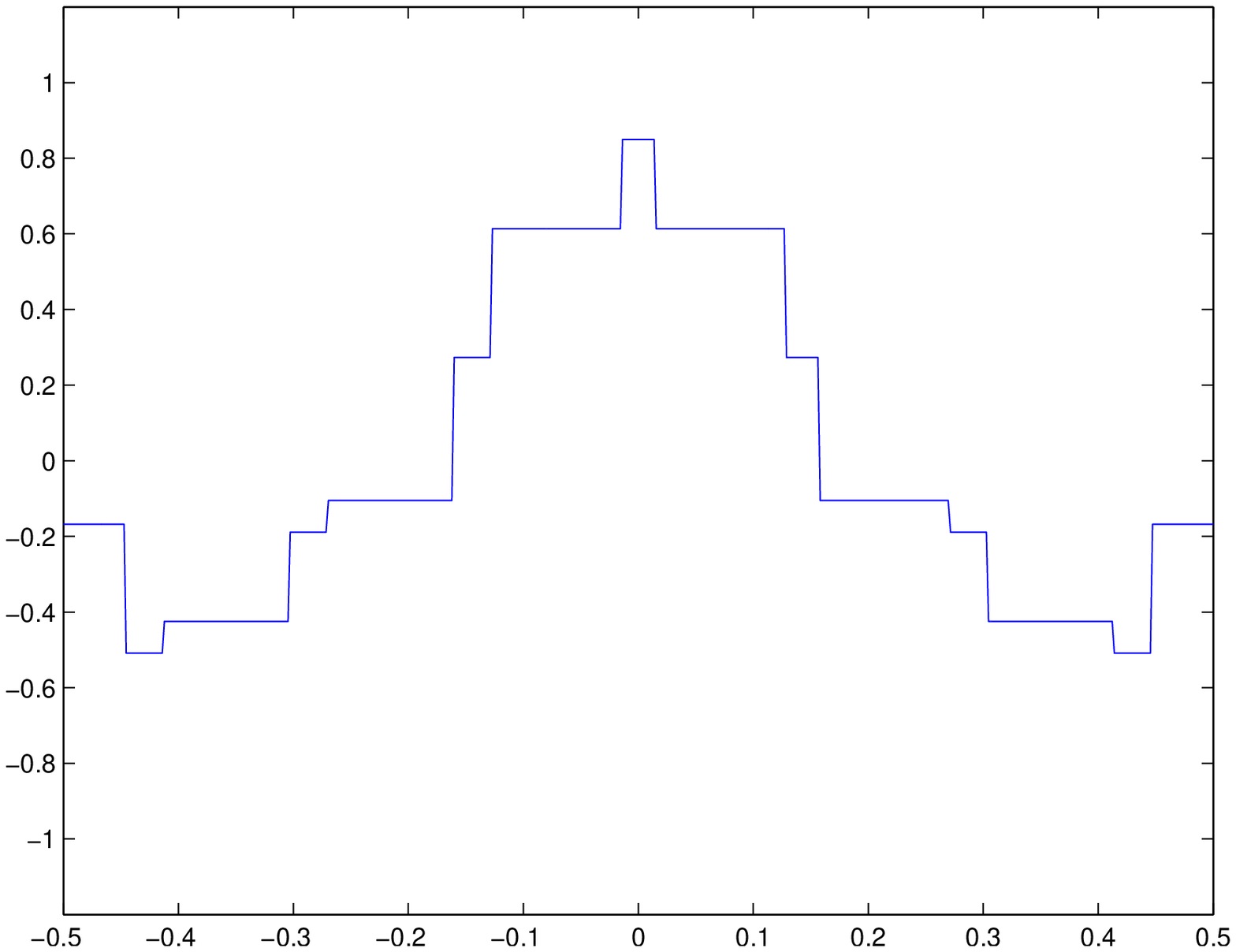}
\includegraphics{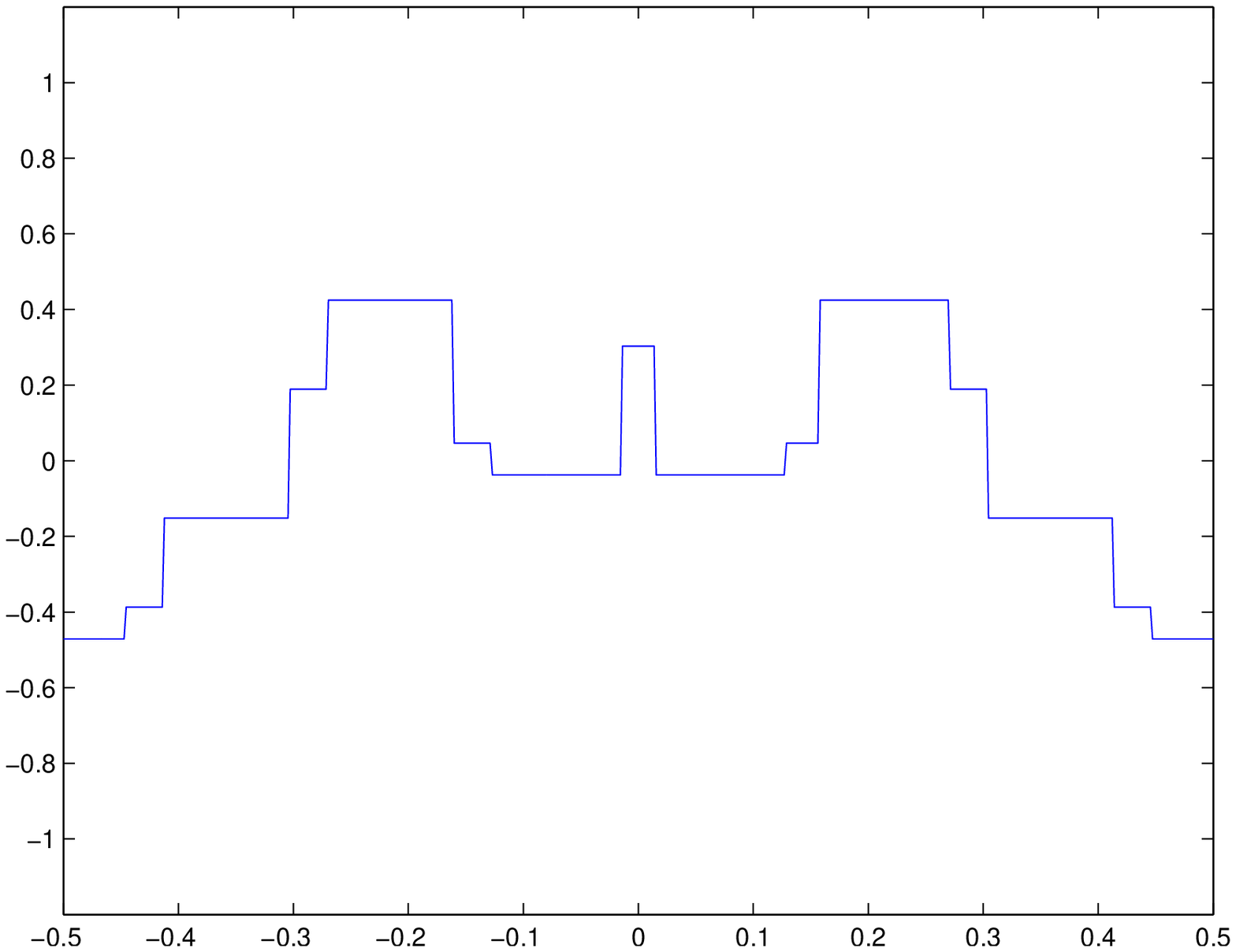}}
\end{center}
\begin{center}
\scalebox{0.2}{\includegraphics{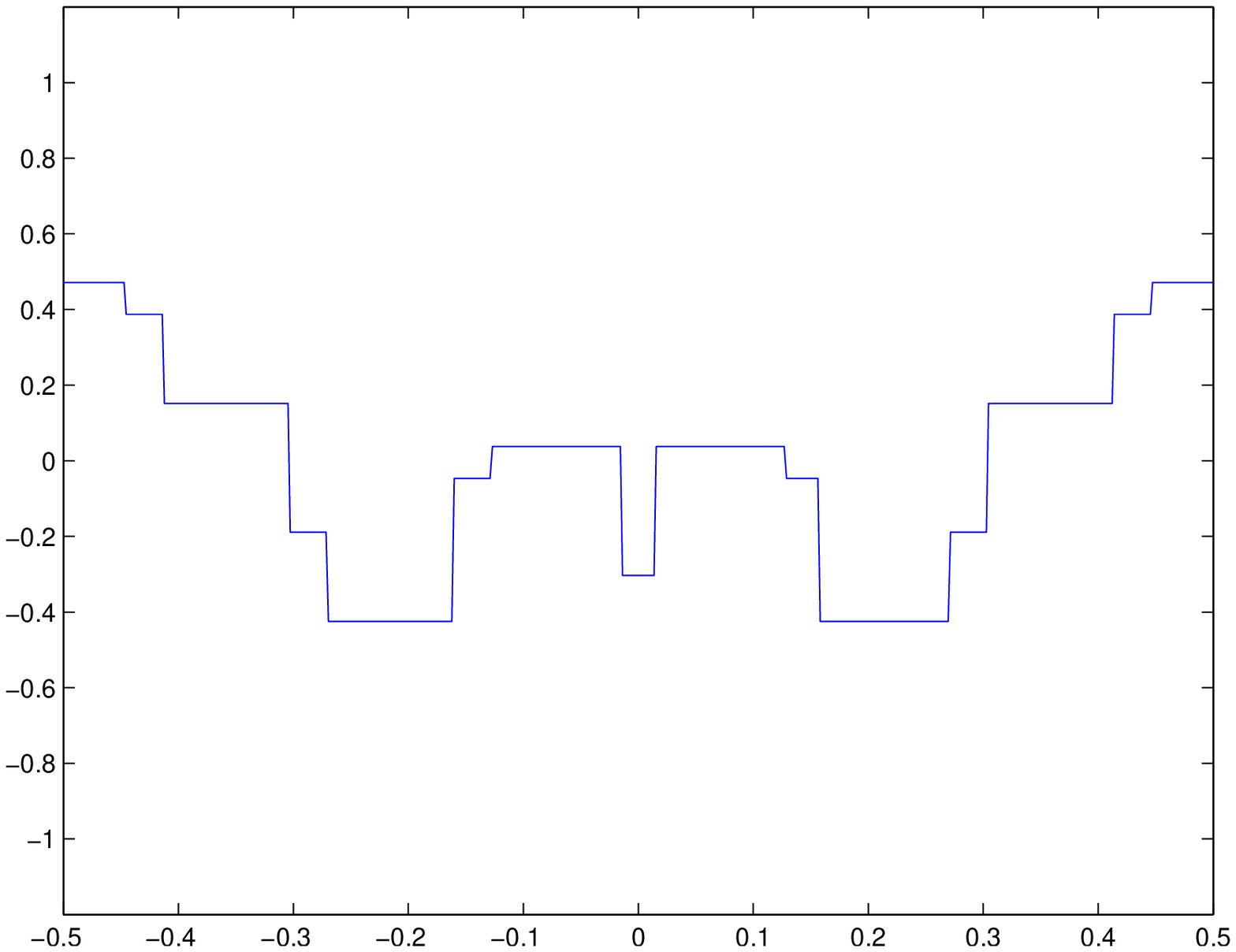}
\includegraphics{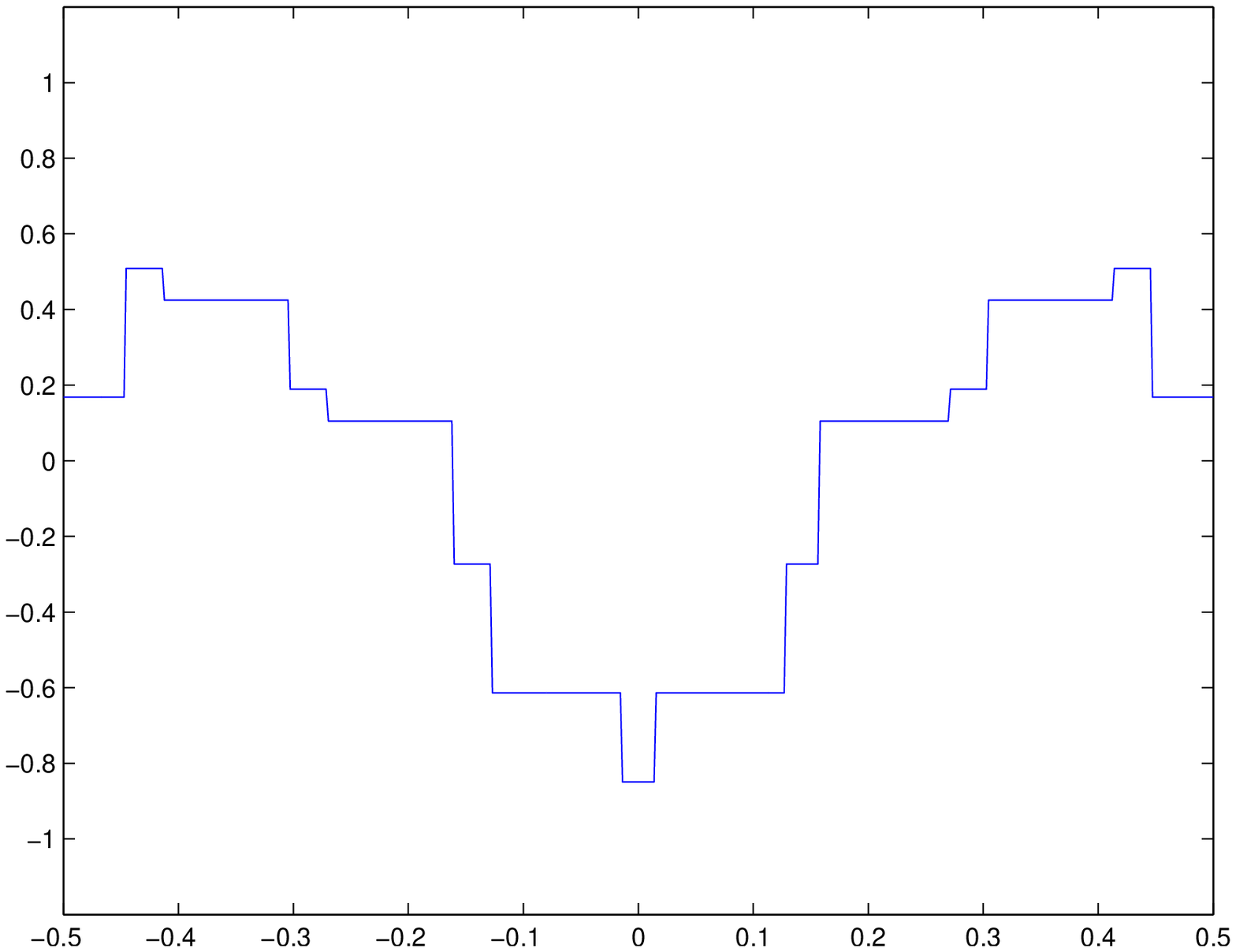}
\includegraphics{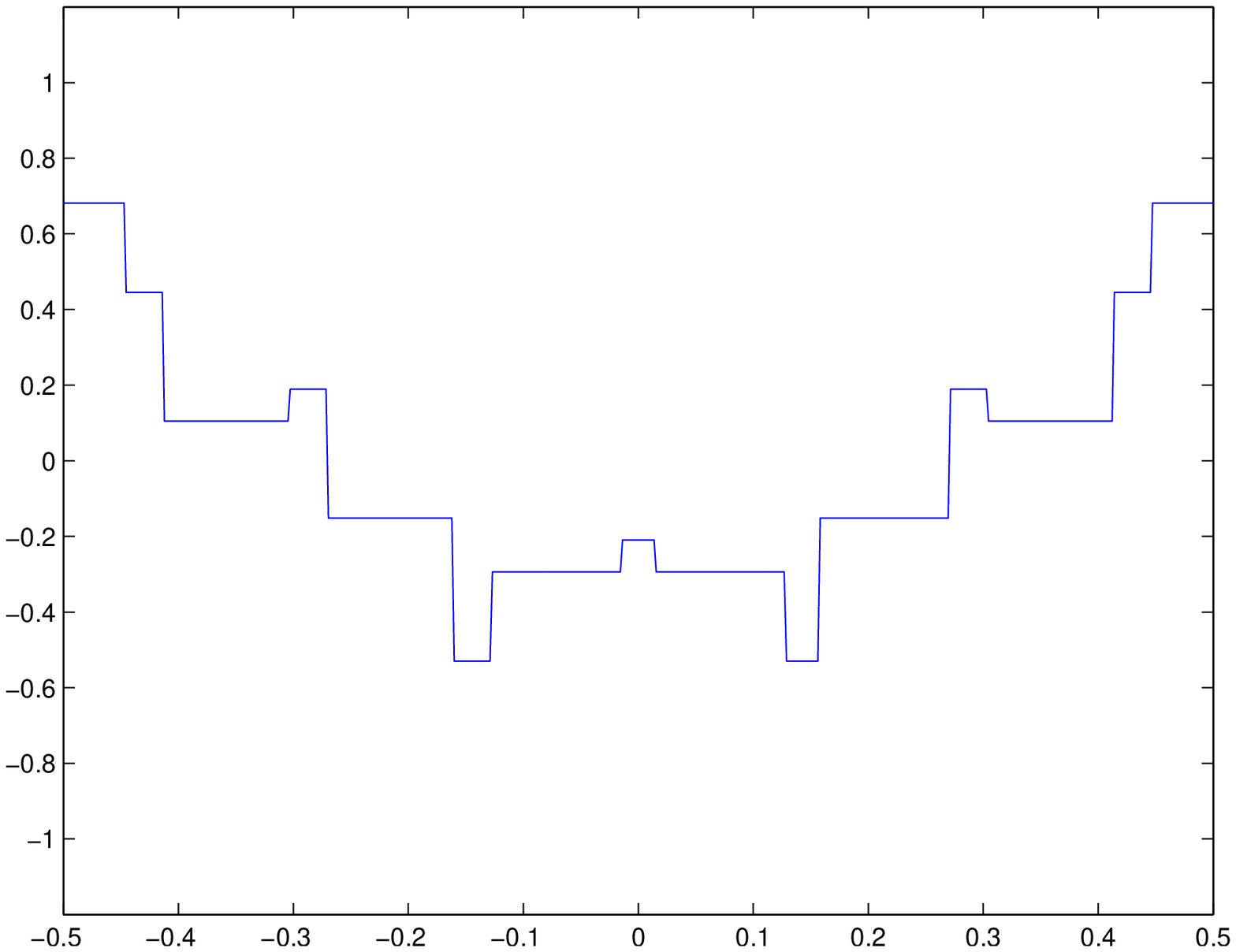}
\includegraphics{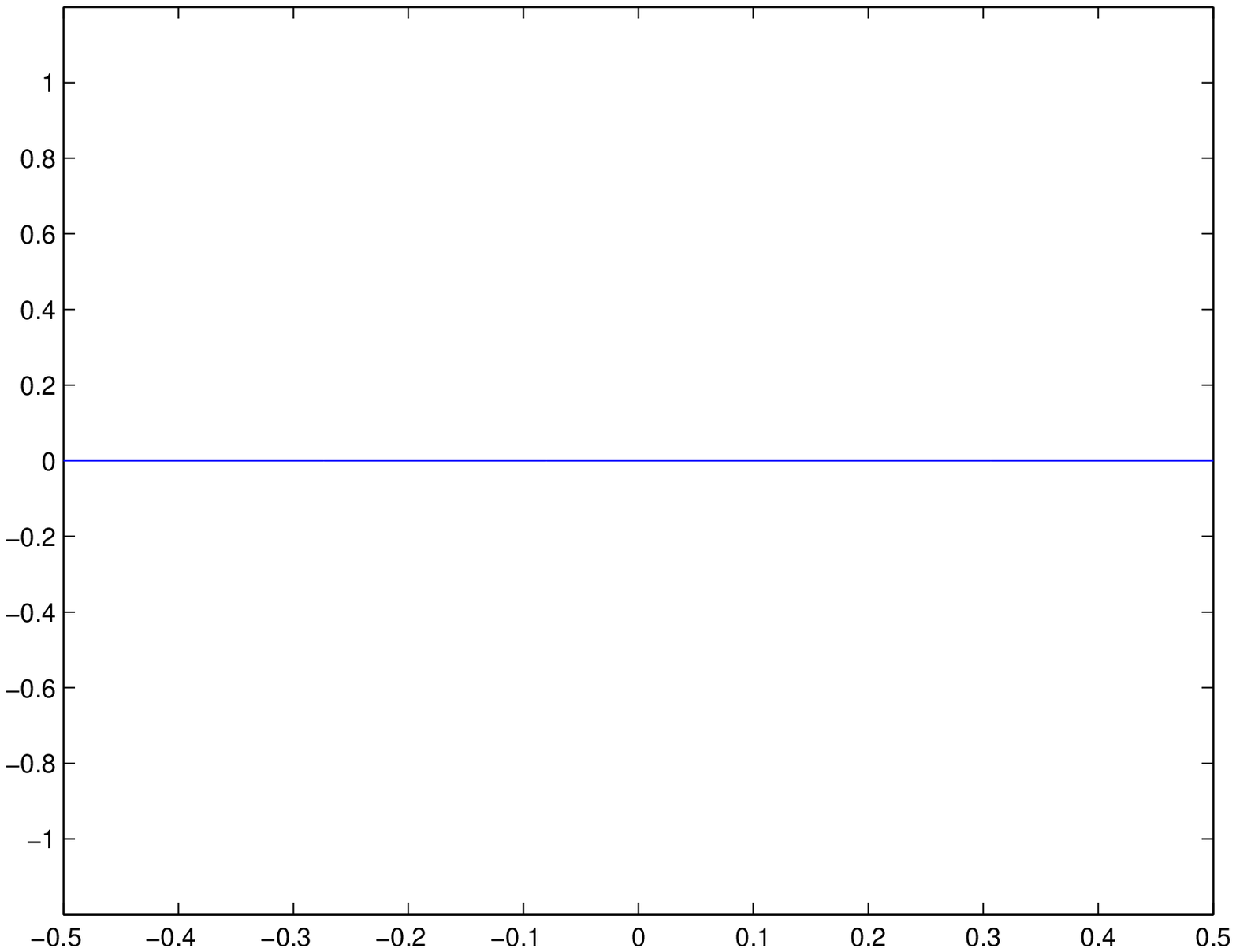}}
\end{center}
\vskip -5mm
\caption{The imaginary part of $U(t,x)$ for $n=2$
at $u/7$ for $u=0,1,2,3$ and $u=4,5,6,7$}
\end{figure}
\end{center}

As discussed above, Part 1 of this Theorem has been noted 
elsewhere.
The same result is true, however, in a much more general context 
than we can prove Part 2, so we state and prove it separately, 
and deduce 1 of Theorem 1 as an immediate corollary.

\vskip 2mm

{\bf Theorem 1.4.}  {\it Let $P$ be a polynomial with integer 
coefficients and $L$ be the differential operator 
$$L=2\pi i P\left(\frac{1}{2\pi i} \frac{\partial}{\partial x}\right),$$
and consider the initial value problem 
\begin{equation}
L U(t,x)=U_t(t,x)
\end{equation}
with $U(0,x)=f(x)$, where $f\in{\hbox{\script D}}$ and $f(x+1)=f(x)$.  Denoting
\begin{equation}
G(u,v;q)=\sum_{w({\hbox{\smallrm mod}}\ q)}e_q(uP(w)-vw),
\end{equation}
at the rational time $t=u/q$ we have
\begin{equation}
U(t,x)
=\frac{1}{q}\sum_{v({\hbox{\smallrm mod}}\,q)}
G(u,v;q)
f\left(x+\frac{v}{q}\right) \ .
\end{equation}
Thus if $U(0,x)$ is piecewise constant, $U(t,x)$ is piecewise constant
at all rational times.}

\vskip 2mm

We can rephrase this theorem as saying that at rational times $t$ 
the solution to the initial value problem $(1.13)$ 
is a linear combination of translates of the initial data.  Note also 
that very generally (see Schmidt 2004) we have $G\ll q^{1/2+\varepsilon}$, 
a bound 
which is not necessary to prove Theorem 2, but which we will use elsewhere; 
in particular note that
\begin{equation}
\sum_{w({\hbox{\smallrm mod}}\,q)}
e_q\left({{uw^n-vw}}\right) \ll (u,v,q)^{1/2}q^{1/2+\varepsilon}
\end{equation}
where $(u,v,q)$ denotes the greatest common divisor of $u,v$ and 
$q$, and where the implied constant depends on $\varepsilon$ and $n$.
As noted in Berry \& Klein 1996, in the special case of the Schr\"odinger
equation these sums are Gauss sums and can be evaluated for general $q$, 
but in general this is not realistic and we must be satisfied with estimates.
We omit the proof of $(1.16)$, merely noting that it follows from Theorem 2.5 
in Schmidt 2004, which contains an extensive discussion of these issues.  

\begin{center}
\begin{figure}[ht]
\begin{center}
\scalebox{0.27}{\includegraphics{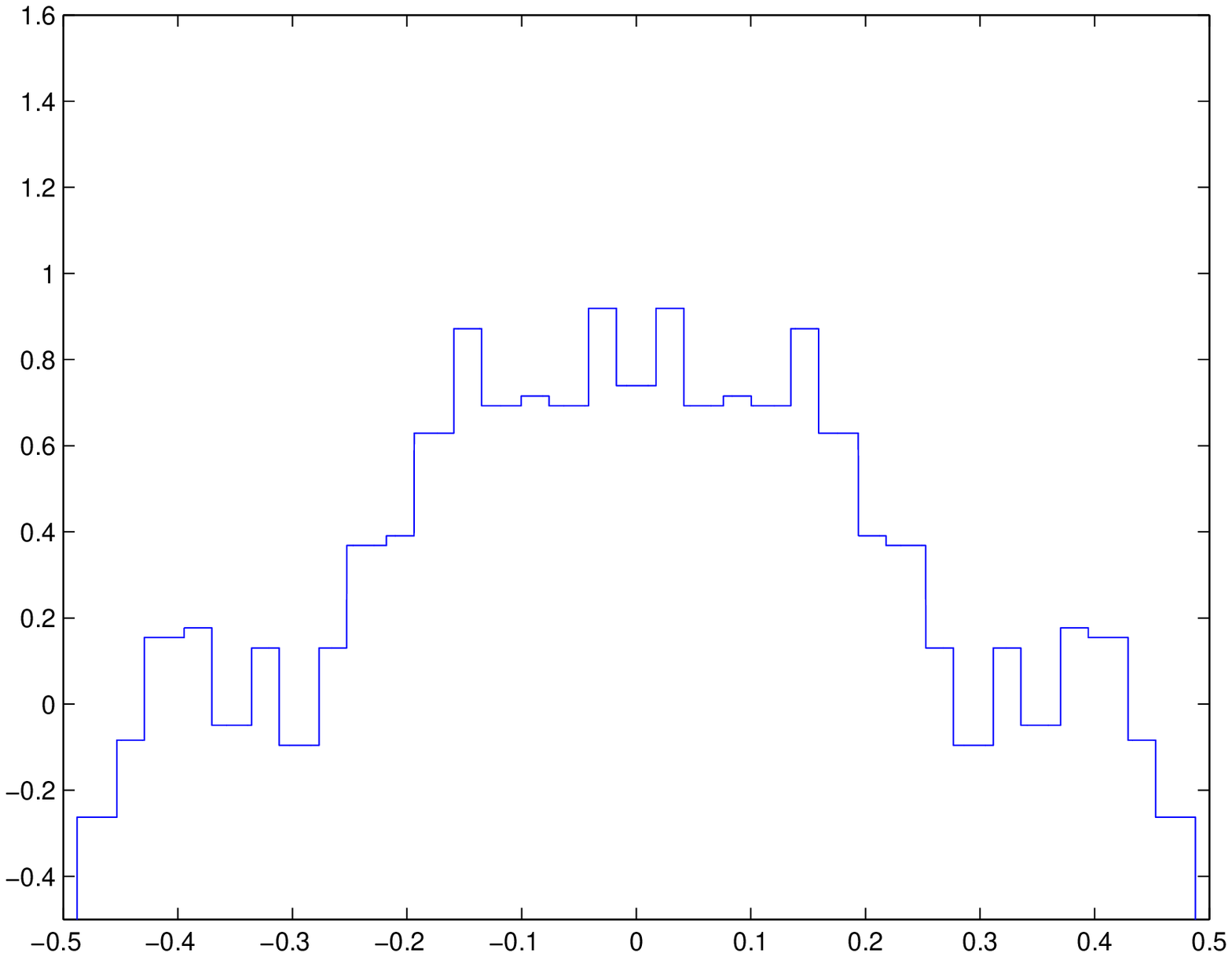}
\includegraphics{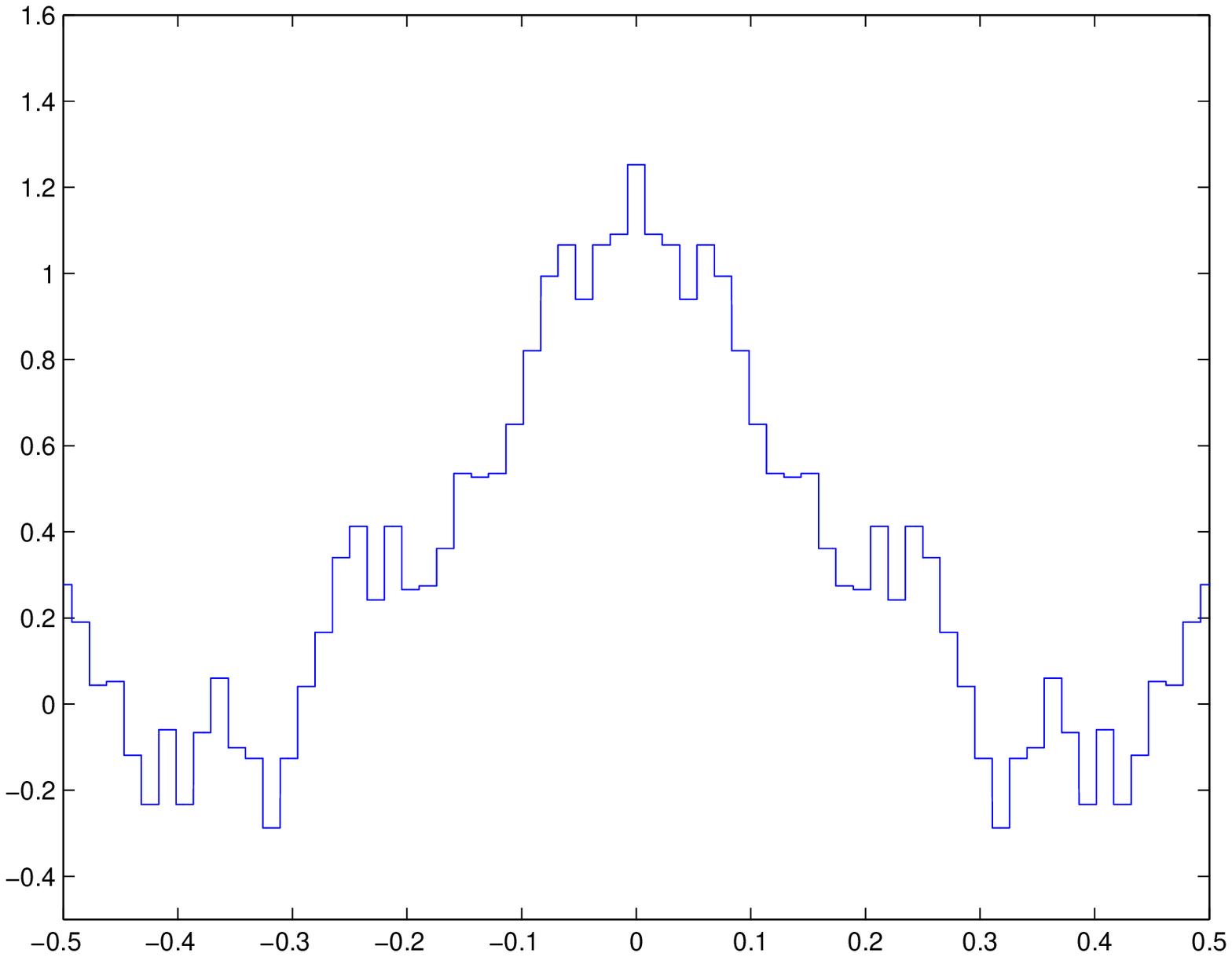}
\includegraphics{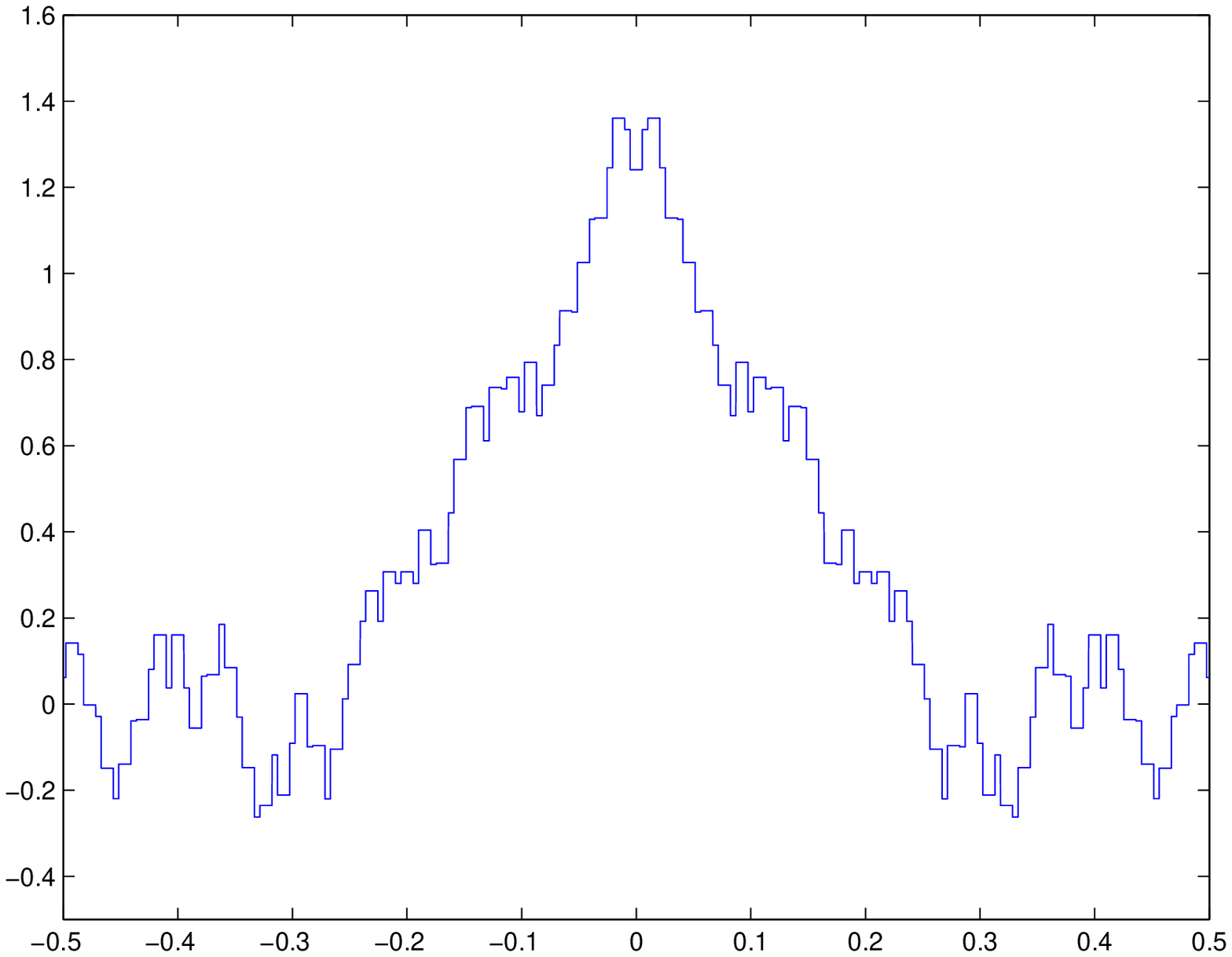}}
\end{center}
\begin{center}
\scalebox{0.27}{\includegraphics{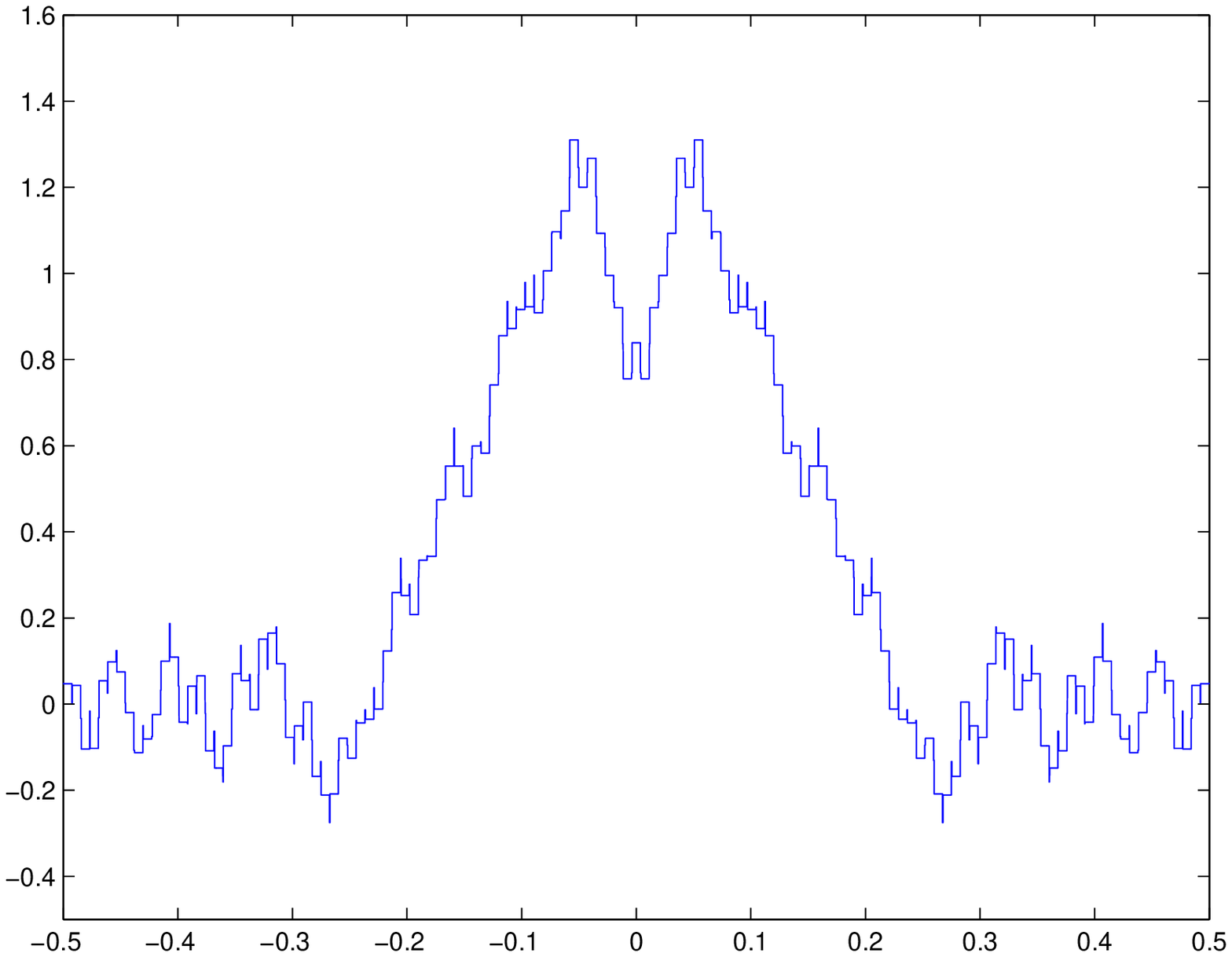}
\includegraphics{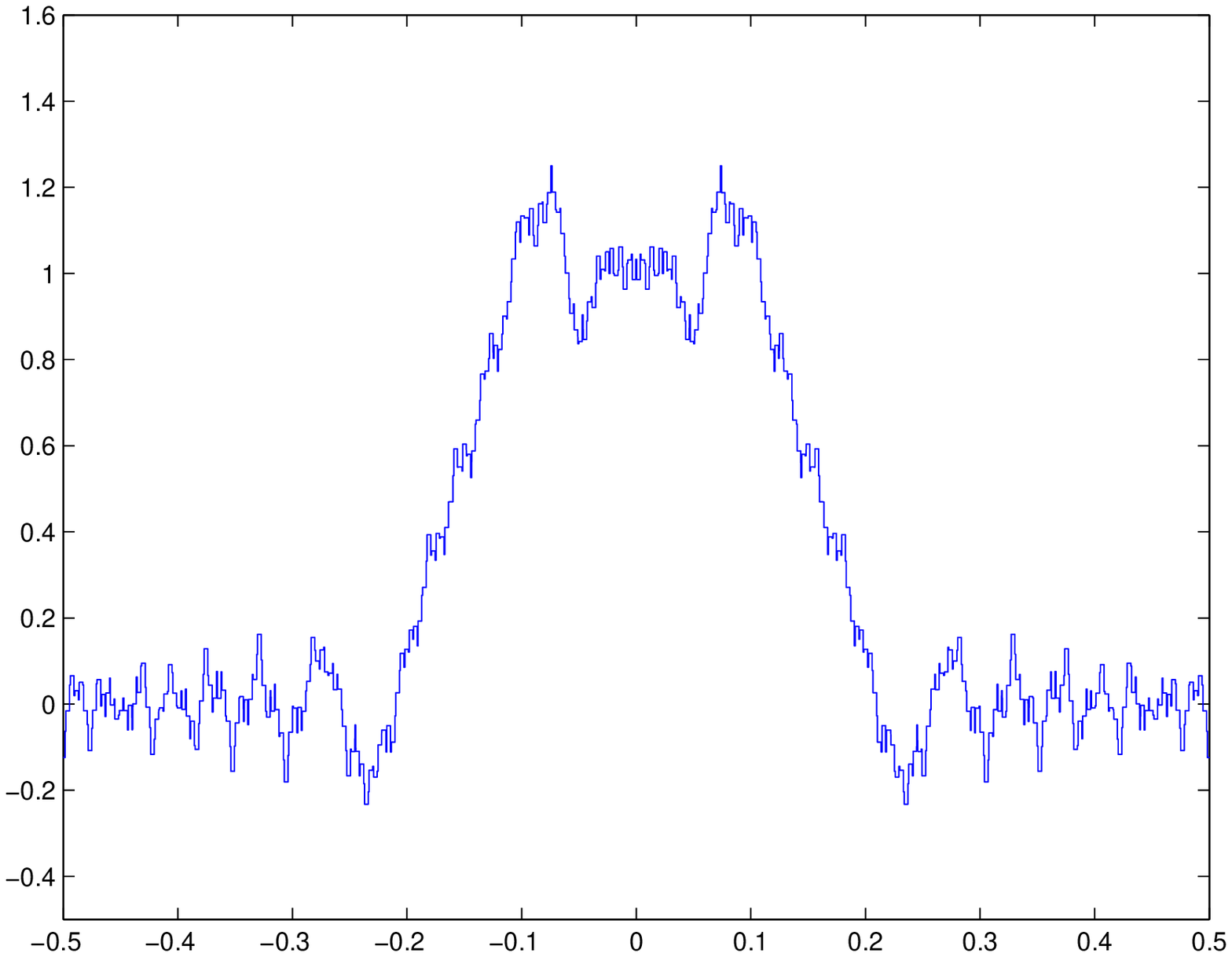}
\includegraphics{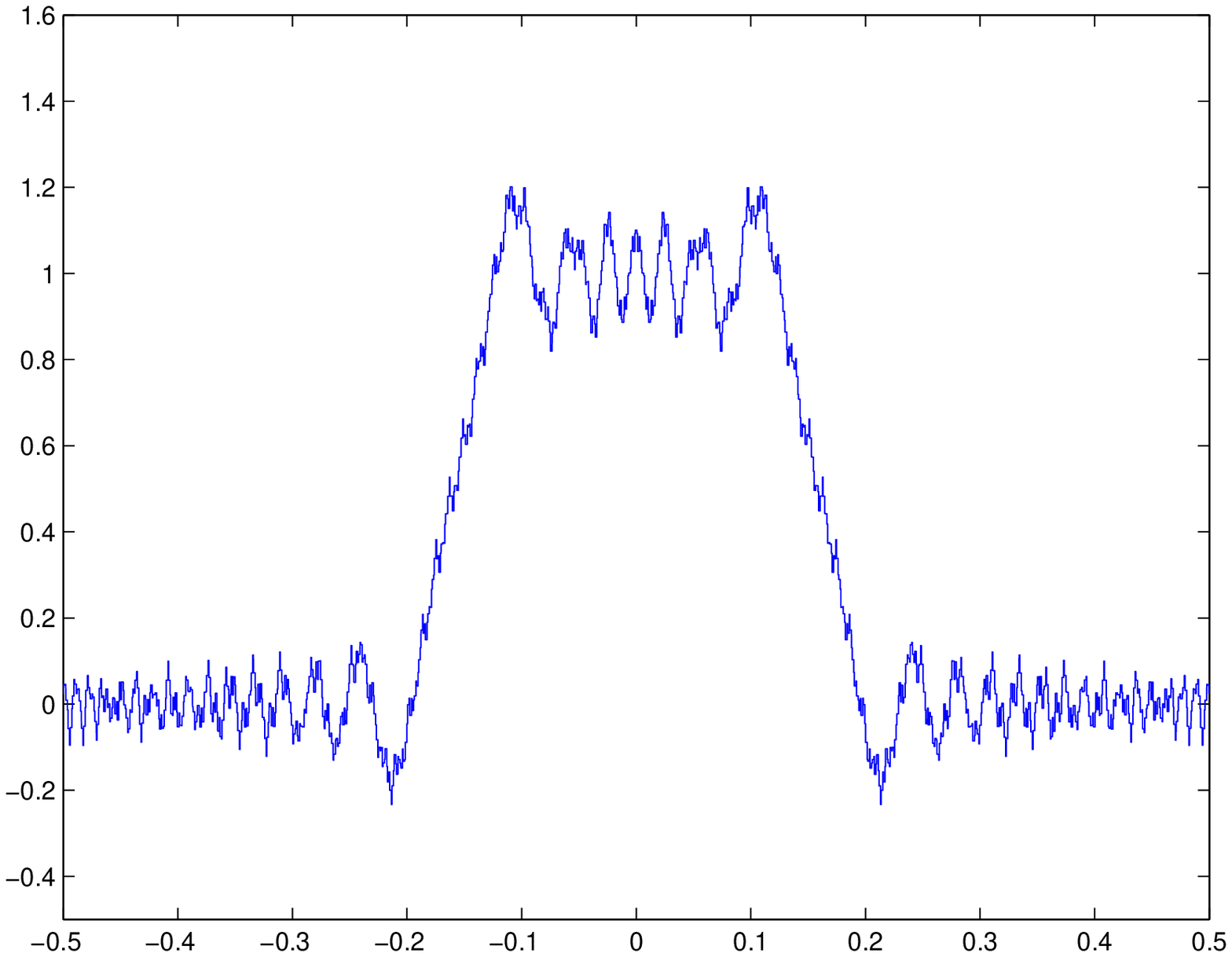}}
\end{center}
\begin{center}
\scalebox{0.27}{\includegraphics{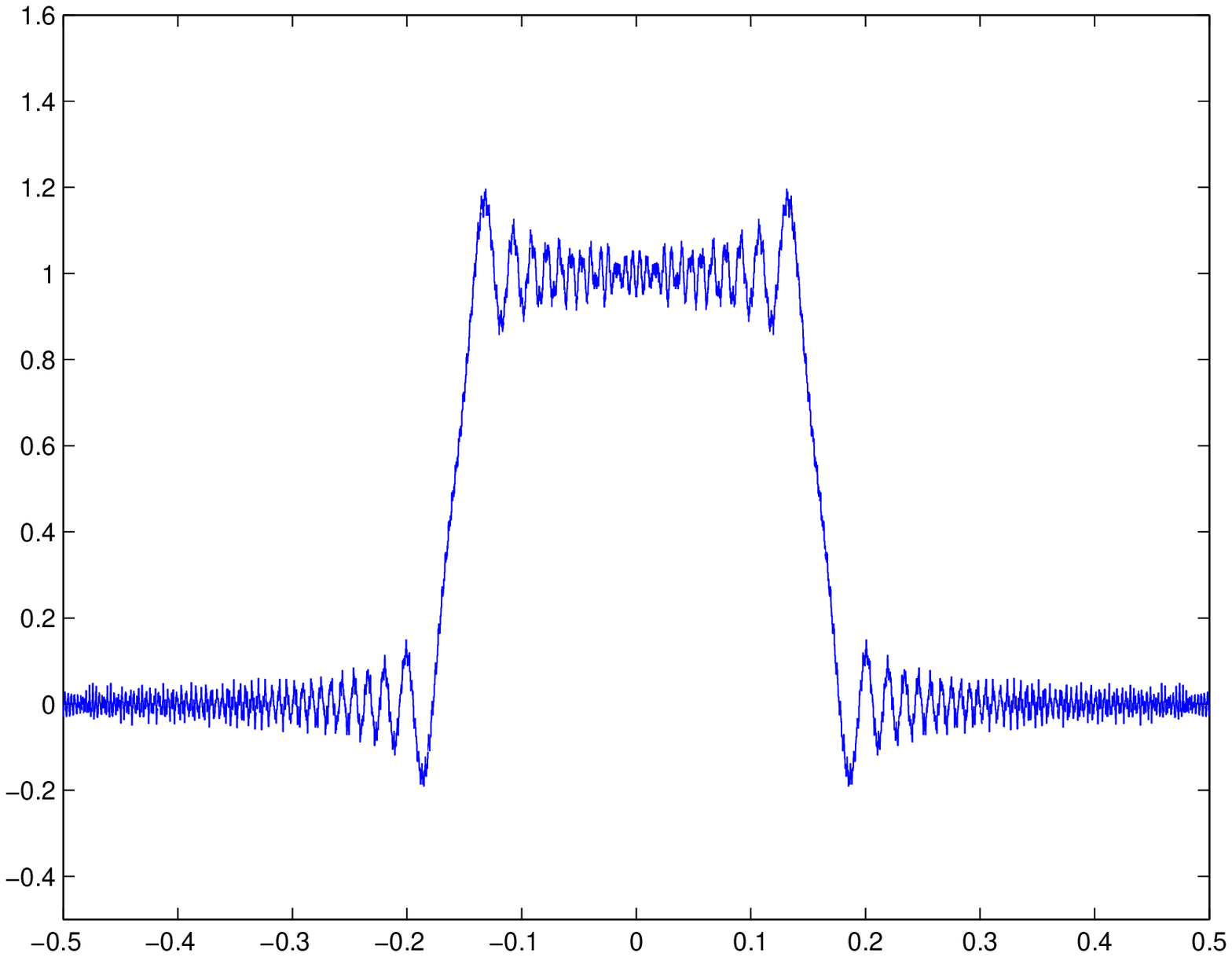}
\includegraphics{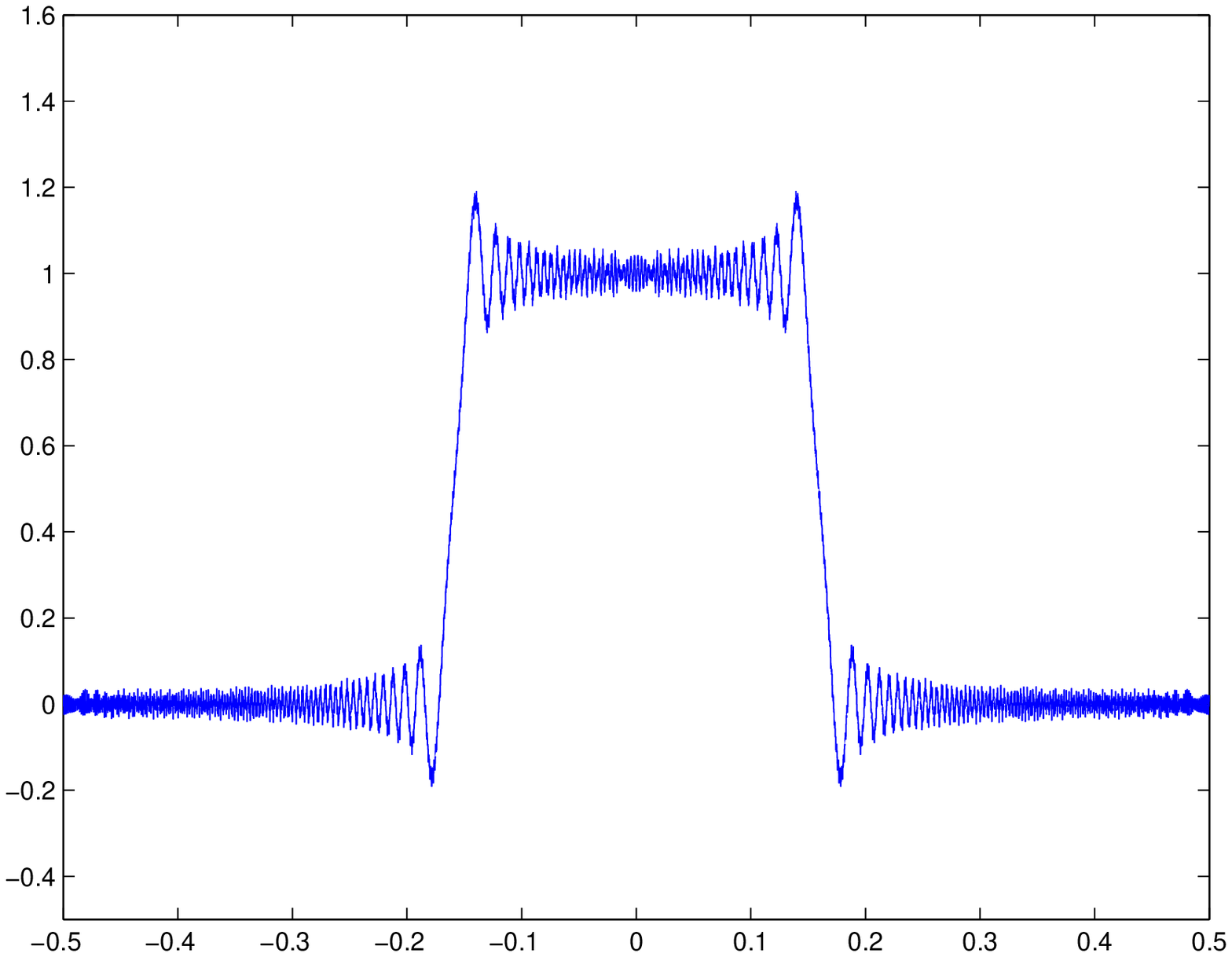}
\includegraphics{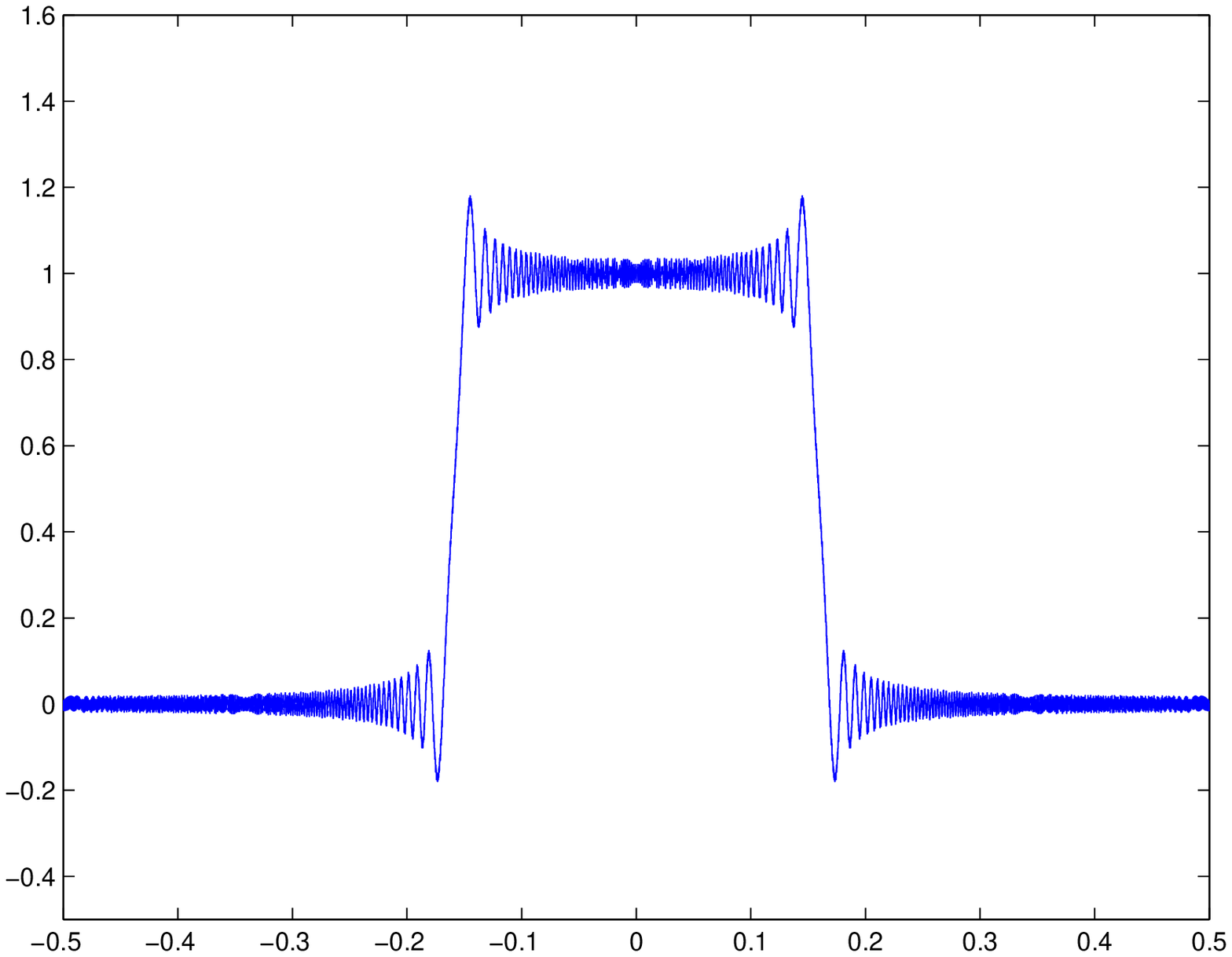}}
\end{center}
\vskip -5mm
\caption{$U(t,x)$ for $n=2$
and $u/q=1/17,1/33,1/65$, $1/129.1/257,1/513$, and $1/2049,1/4097,1/7374$}
\end{figure}
\end{center}

\begin{center}
\begin{figure}[ht]
\begin{center}
\scalebox{0.4}{\includegraphics{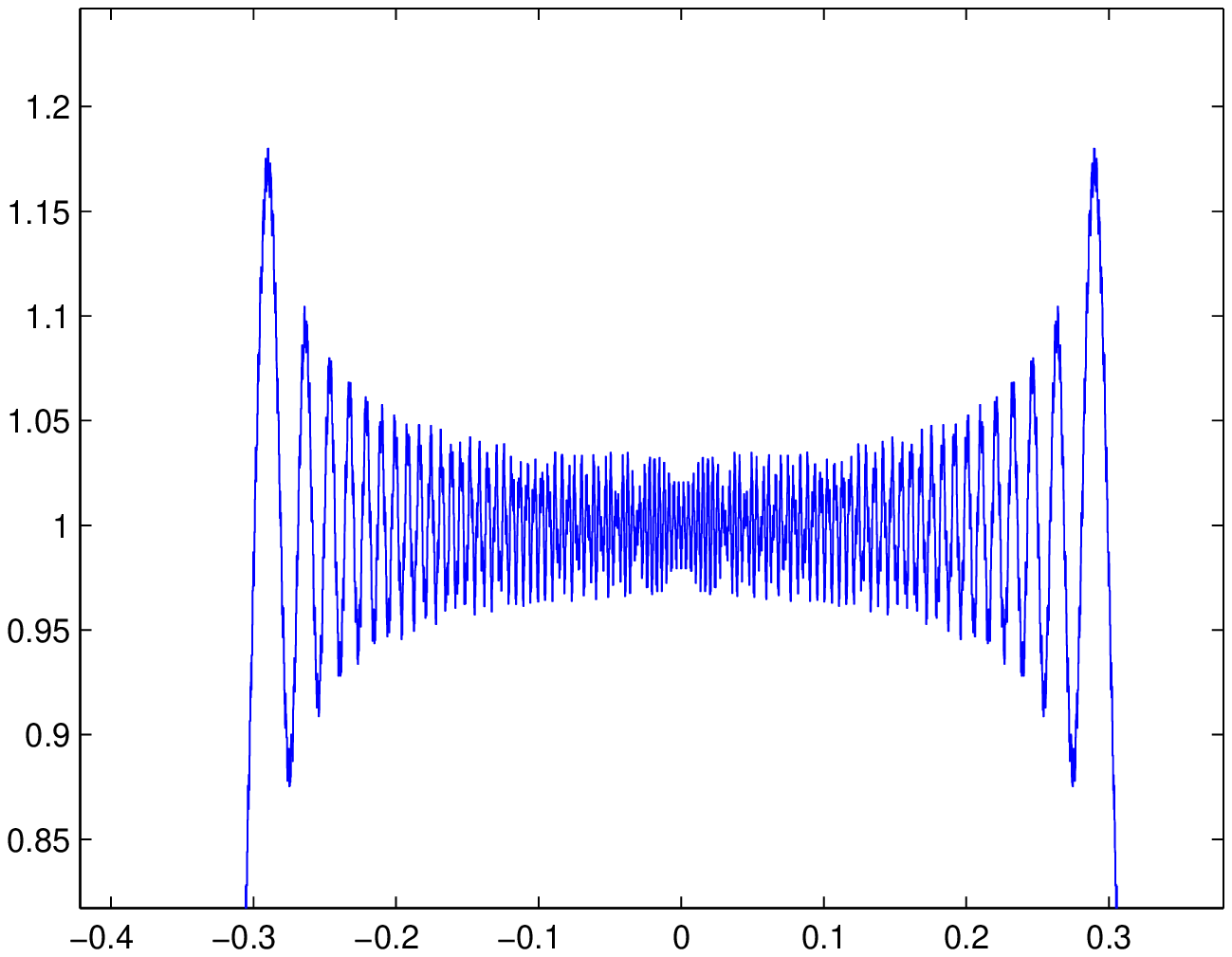}
\includegraphics{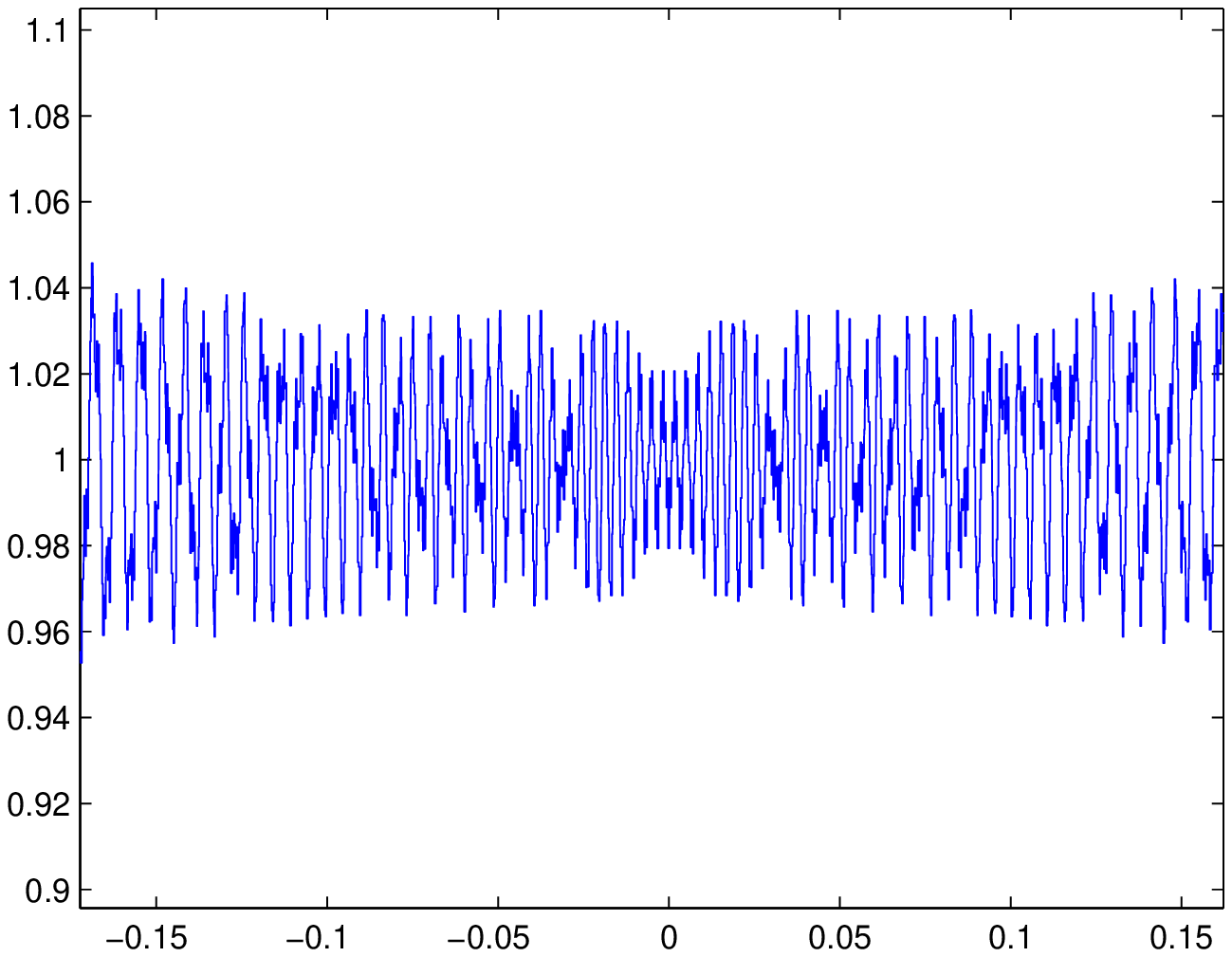}}
\end{center}
\begin{center}
\scalebox{0.4}{\includegraphics{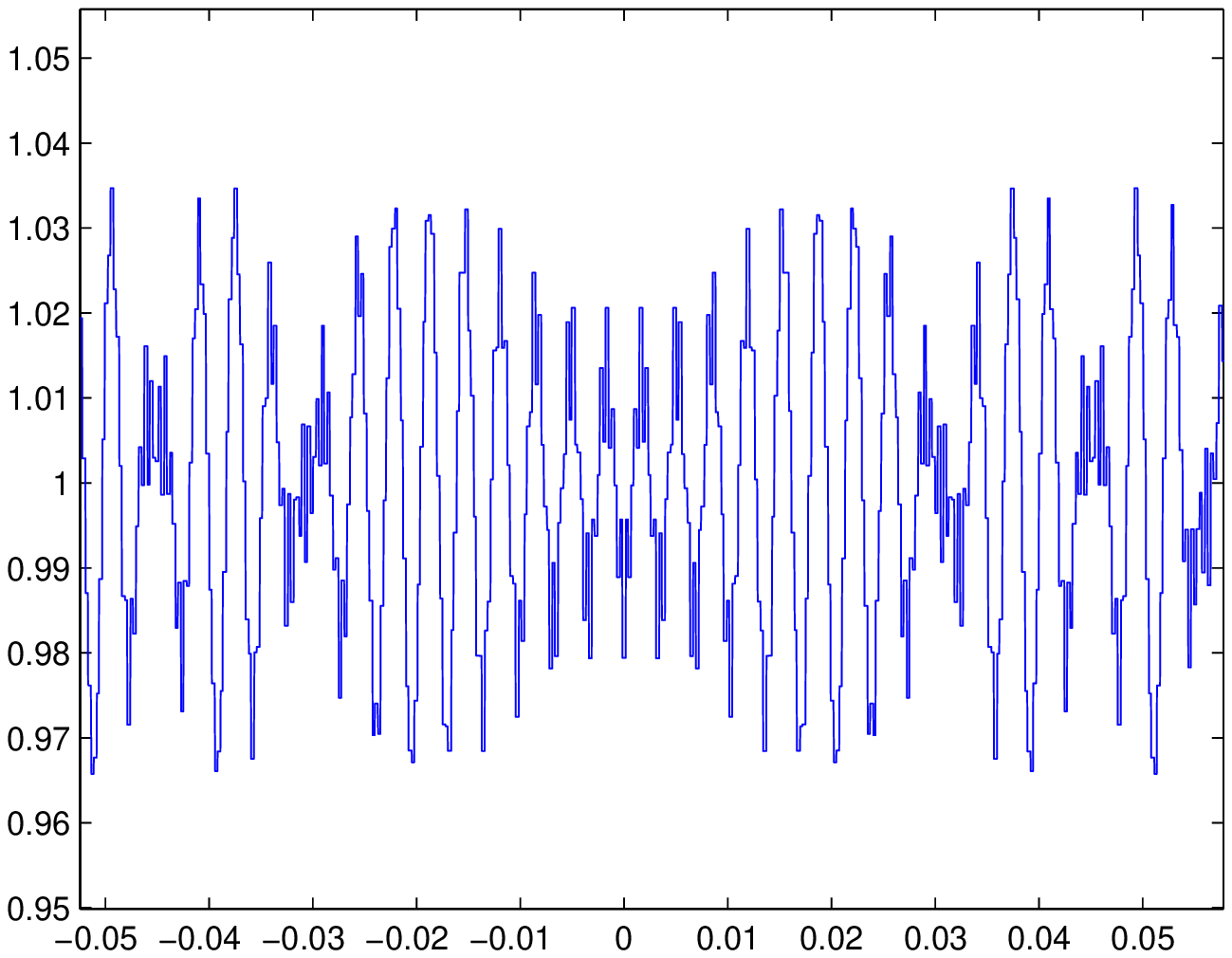}
\includegraphics{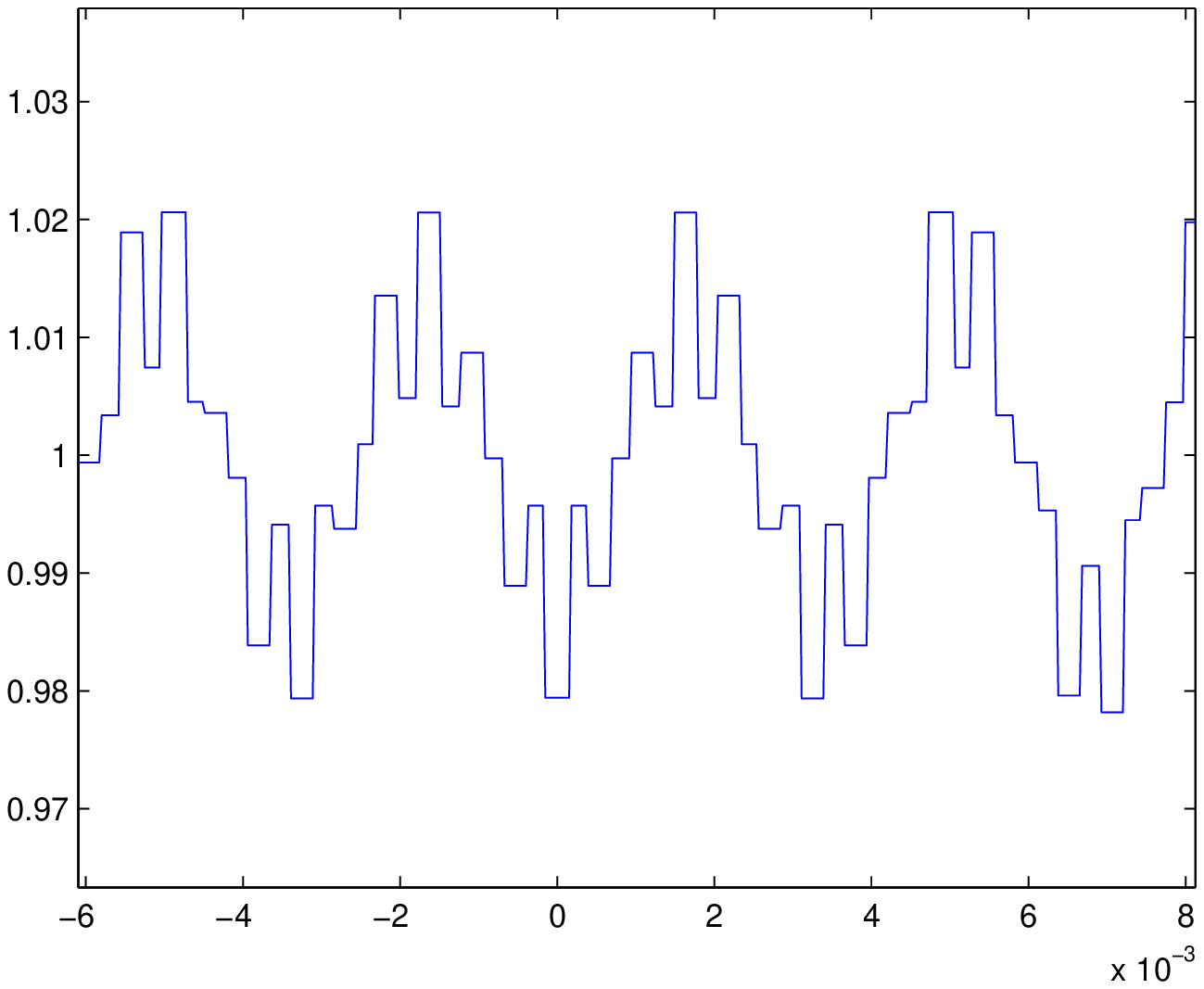}}
\end{center}
\vskip -5mm
\caption{Zoomed-in sections from $U(t,x)$ for $n=2$
and $u/q=1/7374$}
\end{figure}
\end{center}

An example is given by the solution to $(1.1)$ in 
the case $n=2$ and $\gamma=1/\pi$, 
at the rational times $t=u/7$.  Its real and imaginary
parts are plotted in Figure 1 and 2, and can be seen to start and finish 
equal to the initial data.  Note that for odd $n$, of course, the solutions 
are purely real. Note also that Figures 1 - 6 and 10 were produced with the 
aid of the explicit closed form expressions afforded by Theorem 1.4, in 
all cases for rational values of $t$.  Indeed, formula (1.15) represents 
the solution $U(x,t)$ as a finite sum.  The ringing effect which 
is the main topic of this paper occurs as $t$ tends to a rational value.  
To illustrate this effect in Figure 3 we plot the real part of 
$U(t,x)$ at rational times $1/q$ for various 
denominators $q$.  The first values $u/q=1/17,1/33,1/65$ are too large for
any asymptotic tendency to be evident, although they clearly show the 
piecewise constant behaviour described in Theorem 1.4, and the number 
of segments grows approximately linearly with $q$.  The initial data starts
to appear for $u/q=1/129,1/257,1/513$, and the characteristic overshoot of 
a ringing effect becomes evident for $u/q=1/2049,1/4097,1/7374$.
Of course, in these last figures the resolution is insufficient to show 
the piecewise constant behaviour, which is clear in zoomed-in plots 
from the case $q=7374$ shown in Figure 4.
The ringing is also clear in the case of odd $n$, which shows some features 
similar to that of even $n$ and others which are quite distinct.  In 
particular, the solutions are not even as functions of $x$, although they 
are purely real, and the overshoot is asymmetrical on the two sides  
of the discontinuities.  This can be seen in Figure 5, which 
graphs the solution for $n=3$, again with $\gamma=1/\pi$ at several times 
$1/q$.  

Figures 1 - 6 and 10 were produced using Matlab.  Fixing $n$ and $q$, 
we first produce a $q\times q$ matrix containing the values of 
$G(u,v;q)$.  Then a vector representing the initial 
data $f(x)$ is produced for a fixed fixed uniform grid of points $x$ 
between $-1/2$ and $1/2$, and the sum in (1.15) evaluated 
as a matrix-vector multiplication.  
Matlab calculates using IEEE 754 double-precision format 
(1 sign bit, 11 exponent bits, 52 mantissa bits), which translates 
into approximately 16 decimal digits of accuracy.  For our purposes, 
this is sufficiently accurate;  in fact the real limitation is not 
the precision but rather the time required to compute the matrix $G$, 
and we are limited by this restriction to roughly $q \approx 100000$.

The main result of this paper is the following theorem, which
confirms the two principal aspects of these plots: that the pointwise limit
is equal to the initial condition, and that there is a ringing effect,
given by an overshoot of fixed amplitude which migrates towards the 
discontinuity.

\vskip 2mm

{\bf Theorem 1.5.} {\it Let $U(t,x)$ be the solution to $(1,1)$
with periodic initial conditions $U(0,x)$ given by $(1.4)$.  Taking 
any $m\ge 2$ and any $\alpha$ between $1/(n-\Delta)$ and $1/(n-1)$, 
let $t\rightarrow 0^+$ following a sequence of points such that 
$t=u/q$ with $t\le mq^{-1/\alpha(n-\Delta)}$ if rational, and  
$t\in {\hbox{\script B}}_{m,\alpha}$ if irrational. 
Then for any $\varepsilon > 0$, 
$$U(t,x)= 
\begin{cases}
1 +O\left(t^{\alpha\Delta/2^{n-1}-\varepsilon}
+t^{\alpha} |x^2-\gamma^2|^{-1} \right)
&\text{if $|x|< \gamma$,} \\
O\left(t^{\alpha\Delta/2^{n-1}-\varepsilon}+
t^{\alpha} |x^2-\gamma^2|^{-1} \right)
 &\text{if $|x|> \gamma$} 
\end{cases}$$
and if $S>0$ is fixed then for $s\in [-S,S]$, 
$$U\left(t,\pm \gamma+st^{1/n}\right)
=\frac{1}{2}\mp \frac{1}{2\pi}\int\limits_{-\infty}^{\infty}
e(y^n) \sin 2\pi sy \frac{dy}{y} +O\left(
t^{\alpha\Delta/2^{n-1}-\varepsilon}
+t^{\alpha n -1}\right)$$
for even $n$ and
$$U\left(t,\pm \gamma+st^{1/n}\right)
=\frac{1}{2} \mp \frac{1}{2\pi}\int\limits_{-\infty}^{\infty}
e(y^n) \sin 2\pi(y^n+ sy) \frac{dy}{y} +O\left(
t^{\alpha\Delta/2^{n-1}-\varepsilon}+t^{\alpha n -1}\right)$$
for odd $n$, where the implied constants depends on $S$, $\gamma$, $m$ and $\varepsilon$.}

\vskip 2mm

Note that it is the dependence of the implied constant on $m$ that forces 
us to consider $t$ tending to zero in the fixed set 
${\hbox{\script B}}_{m,\alpha}$ rather than in the union of all such, 
which is ${\hbox{\script A}}$.
Nonetheless, there are plenty of points in the set 
${\hbox{\script B}}_{m,\alpha}$, as shown in part 6 of Lemma 1.2.
The restriction on rational values of $q$ is essentially a requirement
that the numerator in $t=u/q$ not remain large; that is, that as $q$
grows, $t$ tends to zero reasonably quickly.  Note that it is 
similar to the definition of the set ${\hbox{\script B}}_{m,\alpha}$
in the irrational case, although less restrictive on $q$. In fact, 
in the rational case the bounds proved are a little stronger than 
stated in the Theorem, as can be seen from the proof.

The ringing effect is described in Theorem 1.5 in a  
renormalised variable $s$; it is interesting to 
graph this renormalised behaviour; in the case $n=2$ 
Figure 6 shows the real and imaginary parts (left and right, 
respectively) of $U(t,s+\gamma t^{1/2})$ for $t=1/202$, $1/1616$, 
$1/6464$, and $1/51712$.  By way of comparison, we can also plot the 
integral to which $U$ is asymptotic; for $n=2$ the real 
and imaginary parts are shown in Figure 7, which appear very close 
to the plots of $U$ for large denominator.  
Although it is computationally difficult to produce plots of the function
$U(u/q,x)$ with large $q$ for $n>3$ larger than $3$, 
we can plot the integral from Theorem 1.5, which is much easier to 
calculate.  (This is, of course, the whole point of proving asymptotic
expressions!)  For instance, Figure $8$ shows the real and imaginary parts 
of the integral for $n=6$.  In the case of odd $n$ we can produce 
similar plots, although as above they do not show the same symmetry 
as for even $n$; for instance we have the graphs for $n=3$ and $n=17$;
Figures 7-9 were produced by Mathematica.  In fact the integrals 
appearing in Theorem 1.5 are, for each integer $n \ge 2$, special 
functions representable using the hypergeometric function 
${}_{p}F_{q}$ (typically with $p=1$ or $p=2$).

As mentioned above, while Theorem 1.5 describes the ringing effect
in detail for times tending to zero, in fact the same phenomenon 
is repeated at each of the jump discontinuities for rational time 
$u/q$, as $t\rightarrow u/q$.  This follows immediately from 
Theorems 1.4 and 1.5 considered together.  An example is shown in Figure 10, 
which shows the solution for $n=2$ and $\gamma=1/\pi$, at
$t=468/3277$, which is a close approximation to $t=1/7$.  Note that the graph
exhibits ringing effects near each of the jump discontinuities
in the graph of $U$ for $t=1/7$ as in Figure 1.

\begin{figure}[ht]
\begin{center}
\scalebox{0.27}{\includegraphics{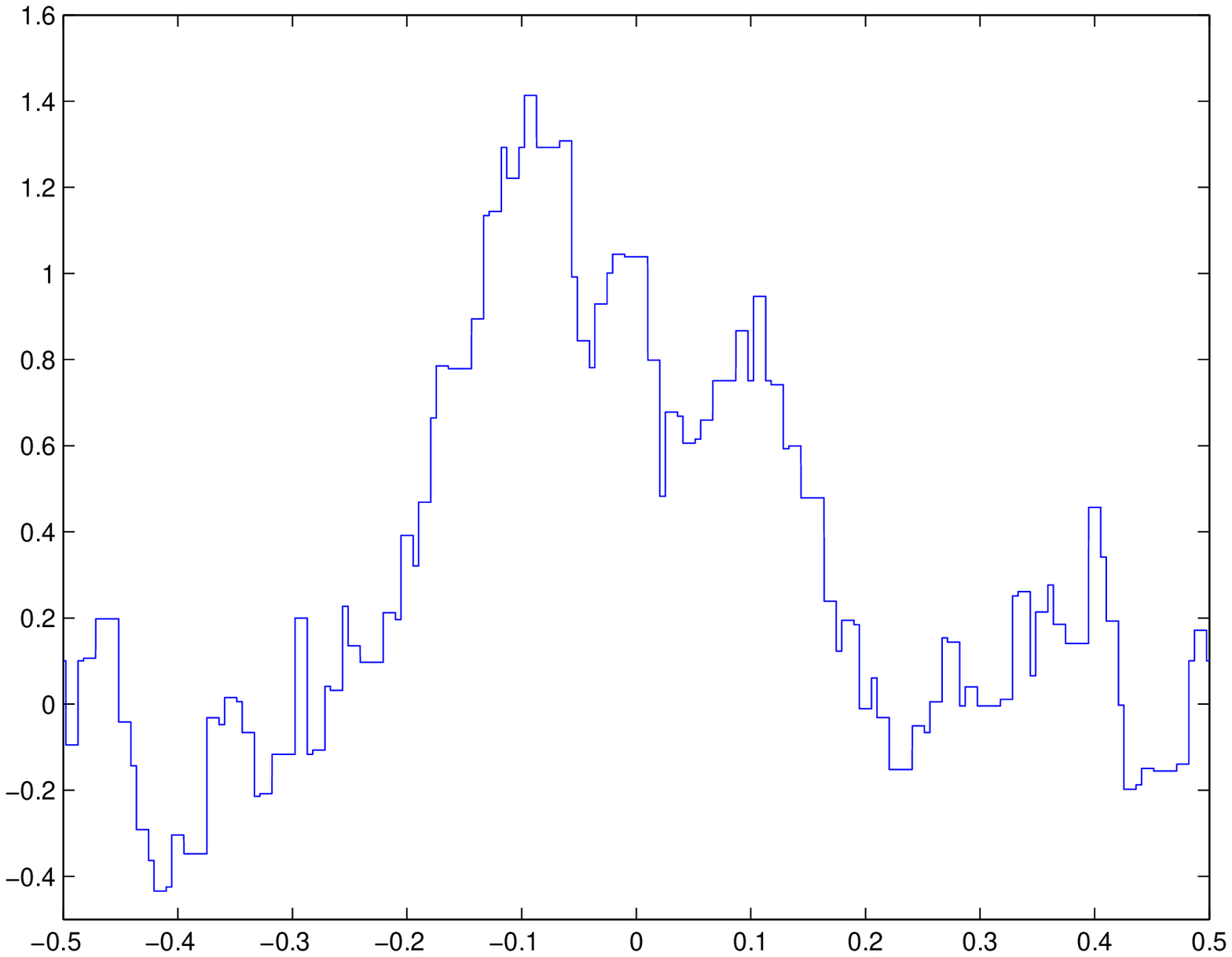}
\includegraphics{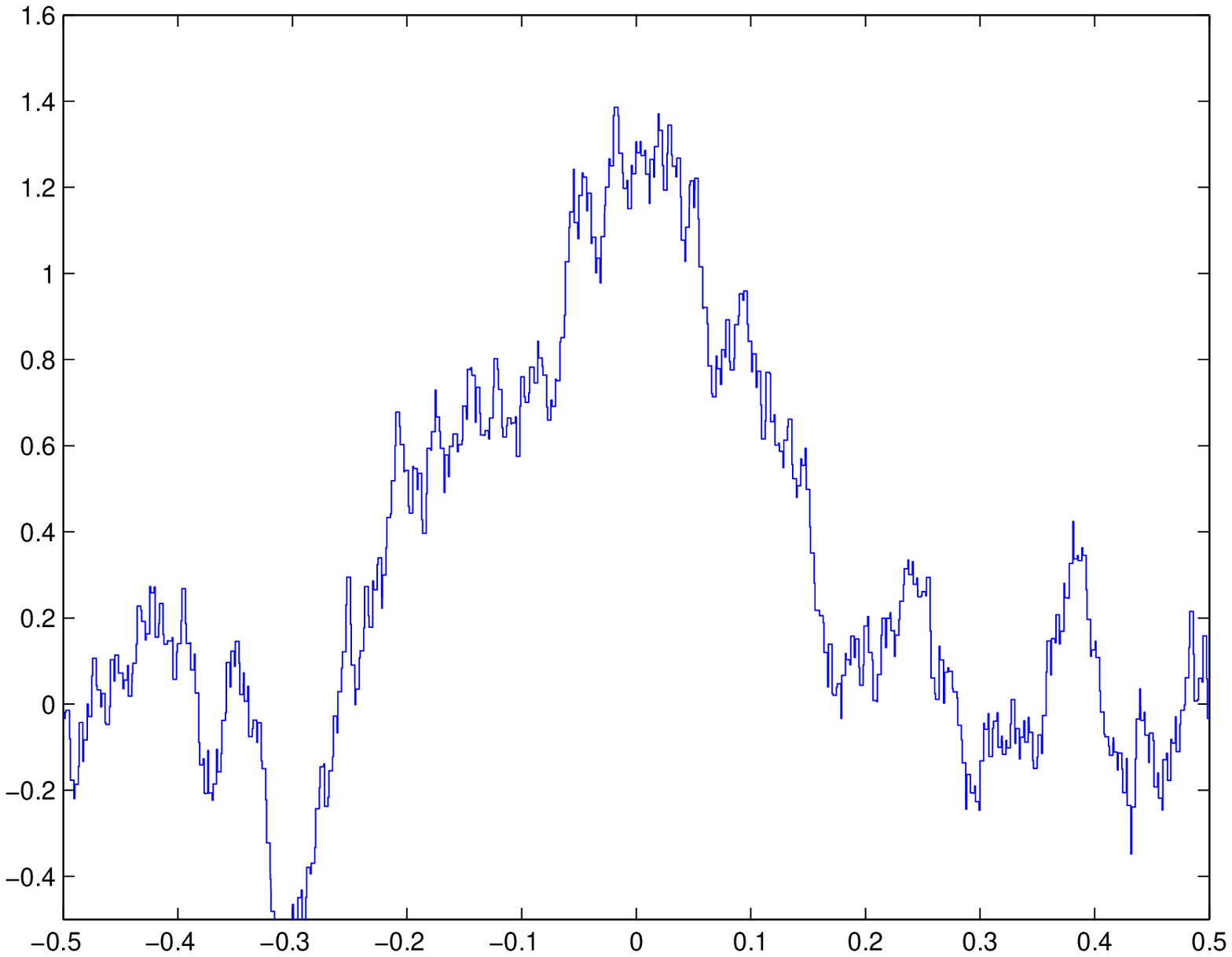}
\includegraphics{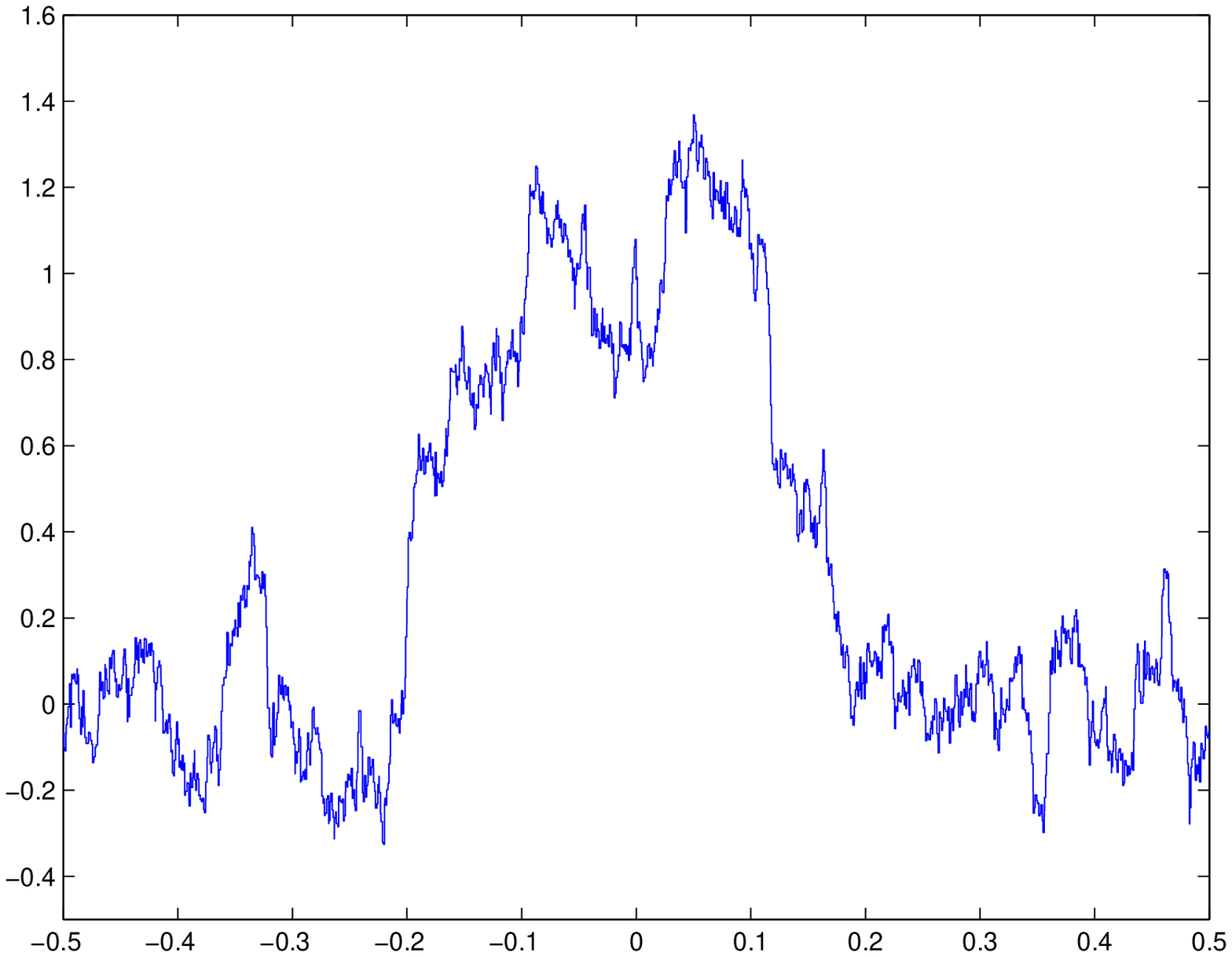}}
\end{center}
\begin{center}
\scalebox{0.27}{\includegraphics{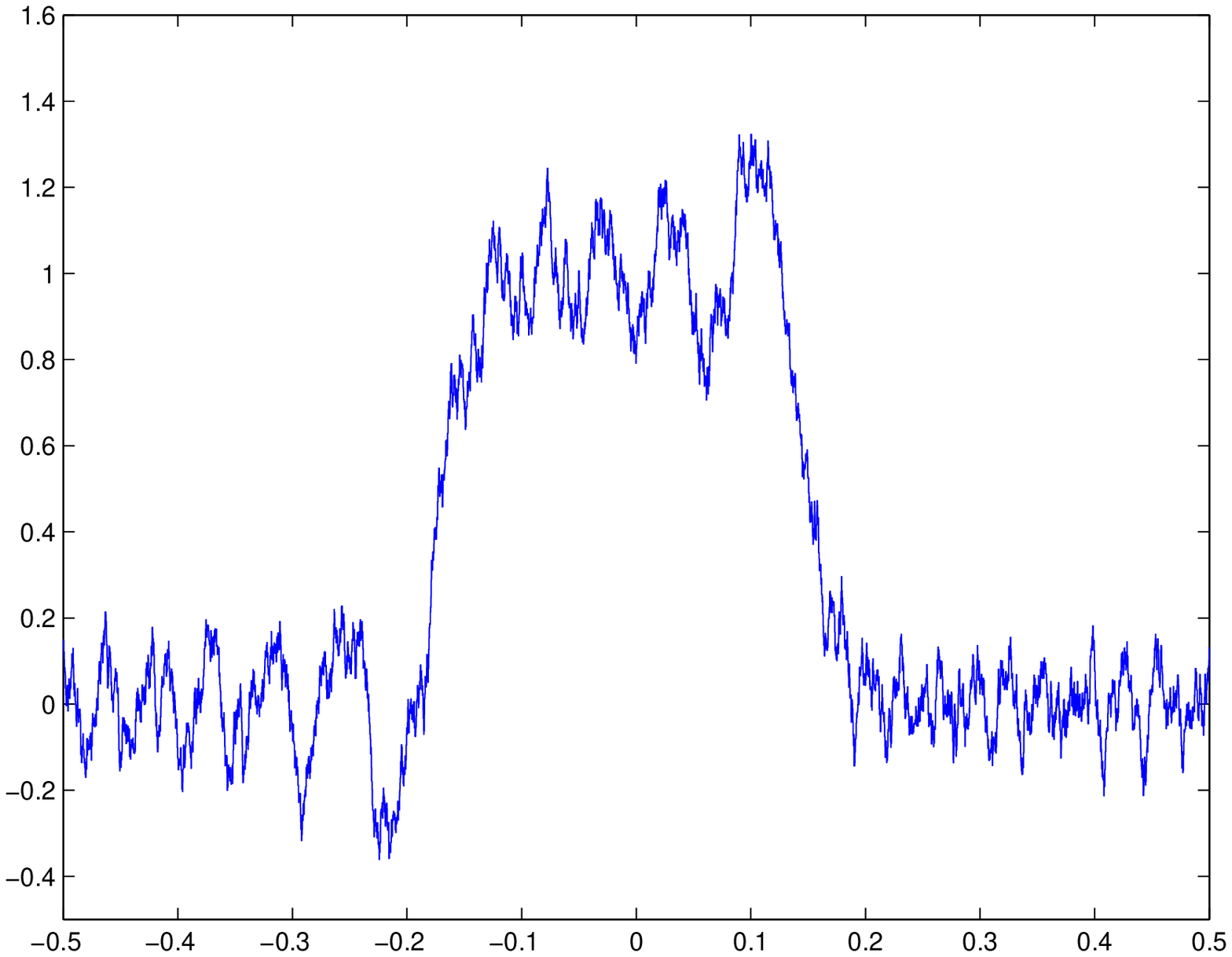}
\includegraphics{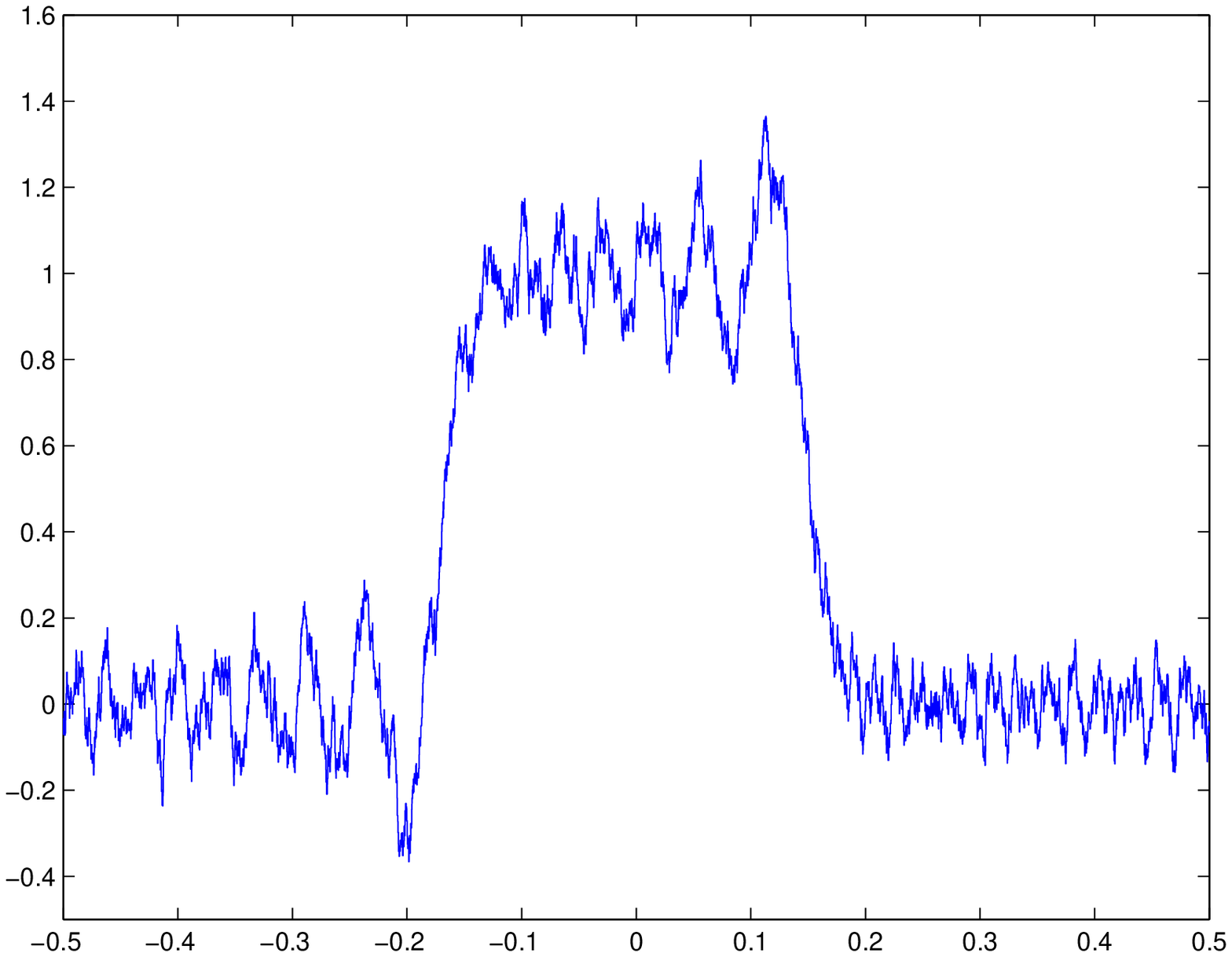}
\includegraphics{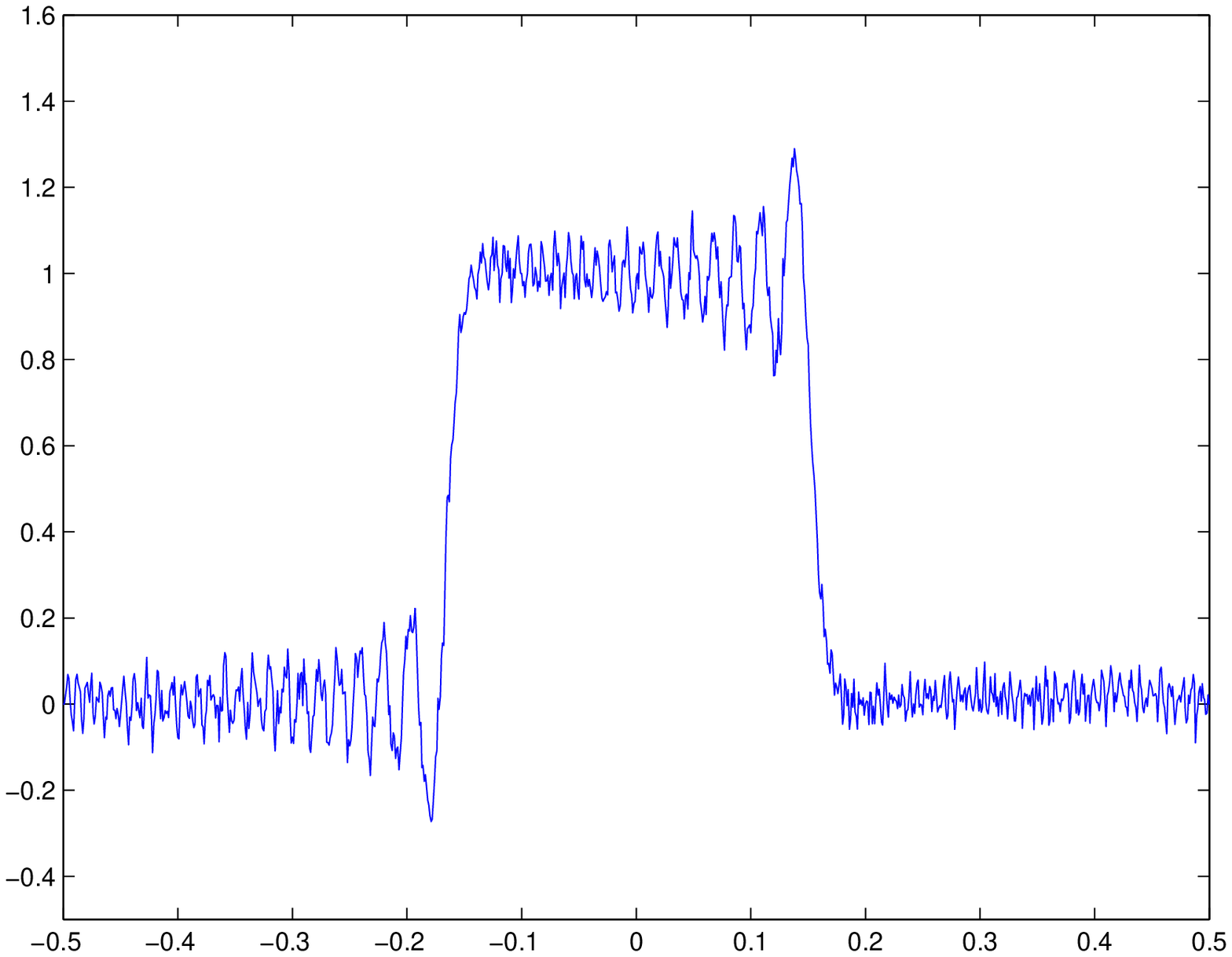}}
\end{center}
\vskip -5mm
\caption{$n=3$ and $u/q=1/65,1/257,1/1025$ and $1/65,1/257,1/1025$}
\end{figure}

\begin{figure}[p]
\begin{center}
\scalebox{0.35}{\includegraphics{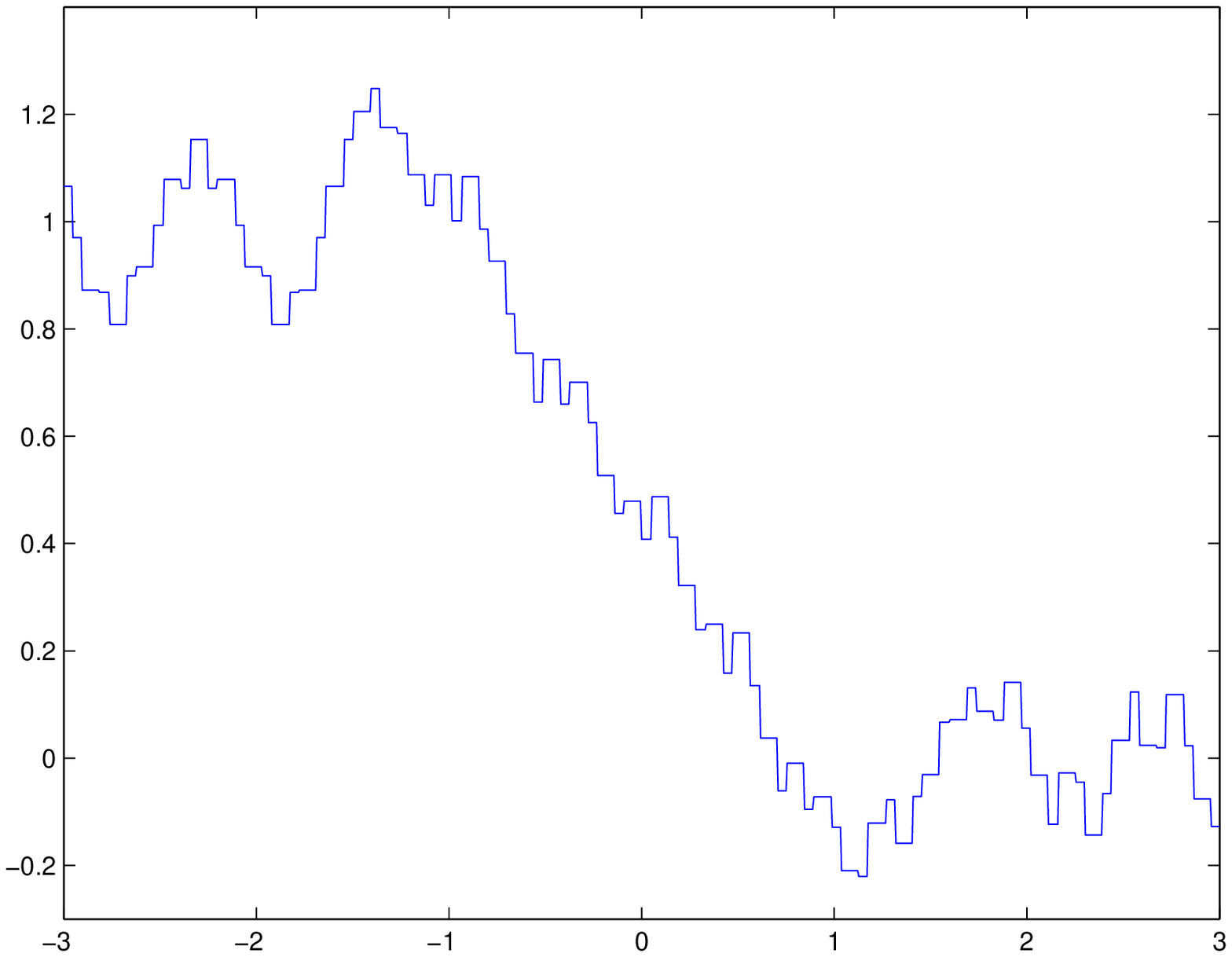}
\includegraphics{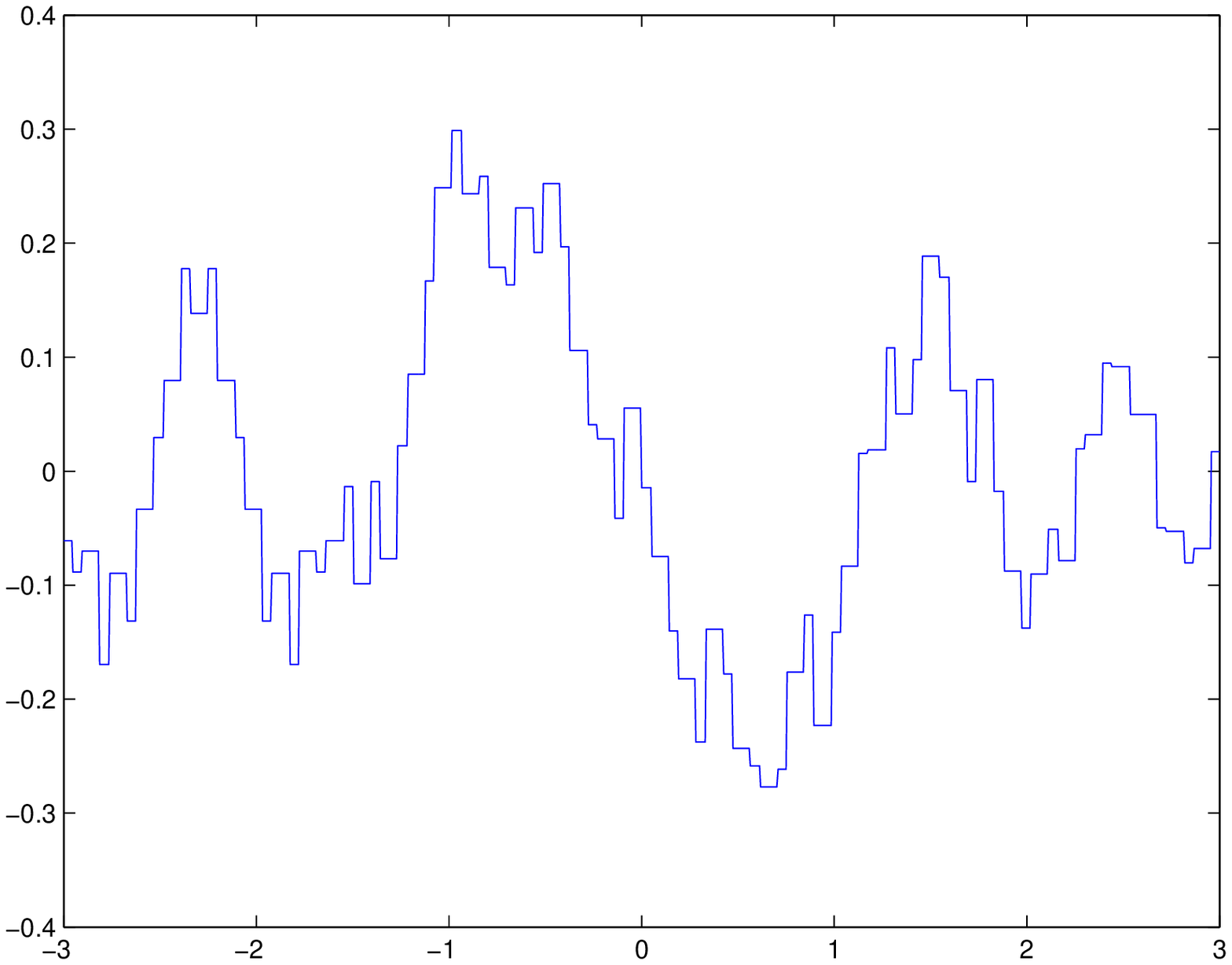}}
\end{center}
\begin{center}
\scalebox{0.35}{\includegraphics{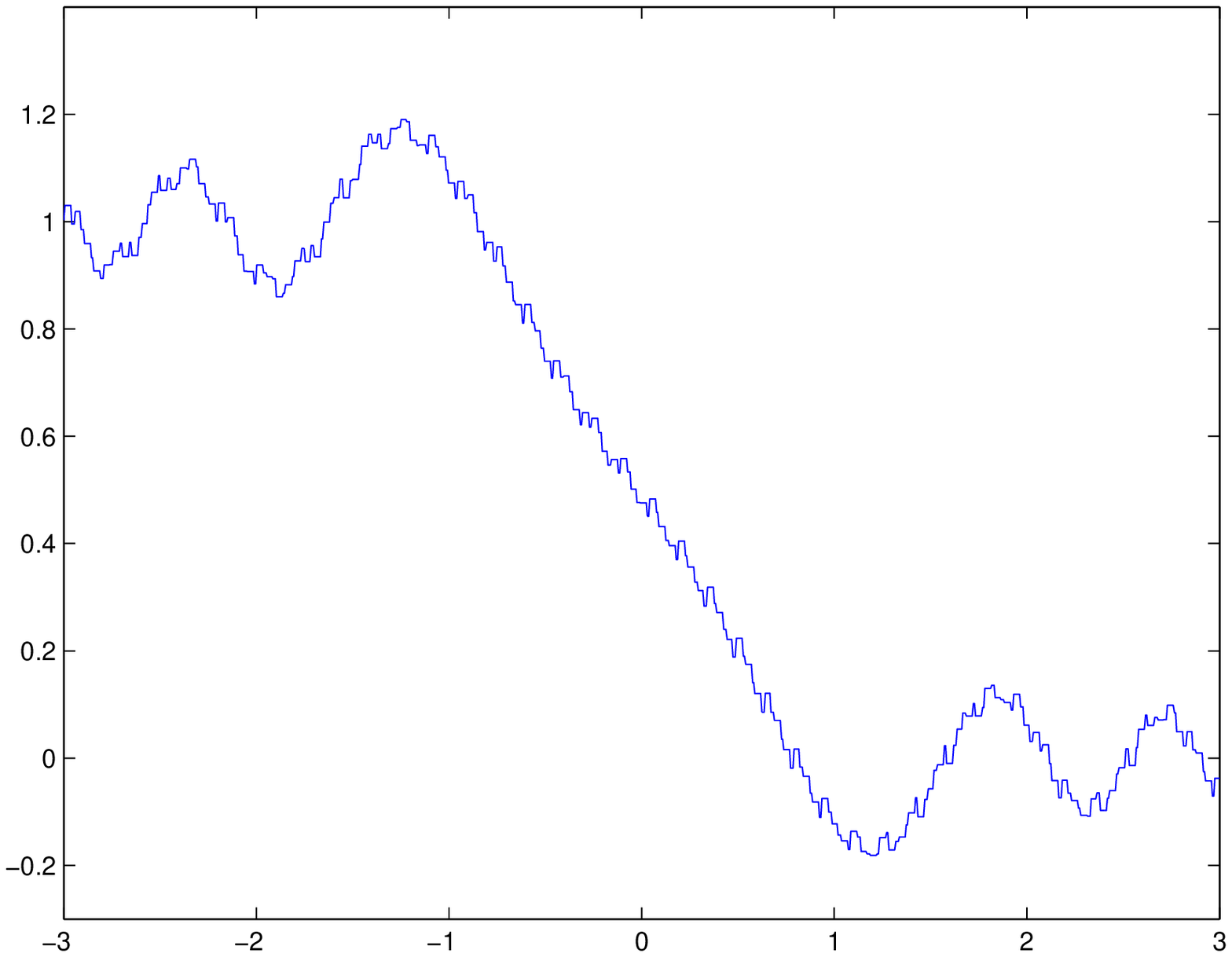}
\includegraphics{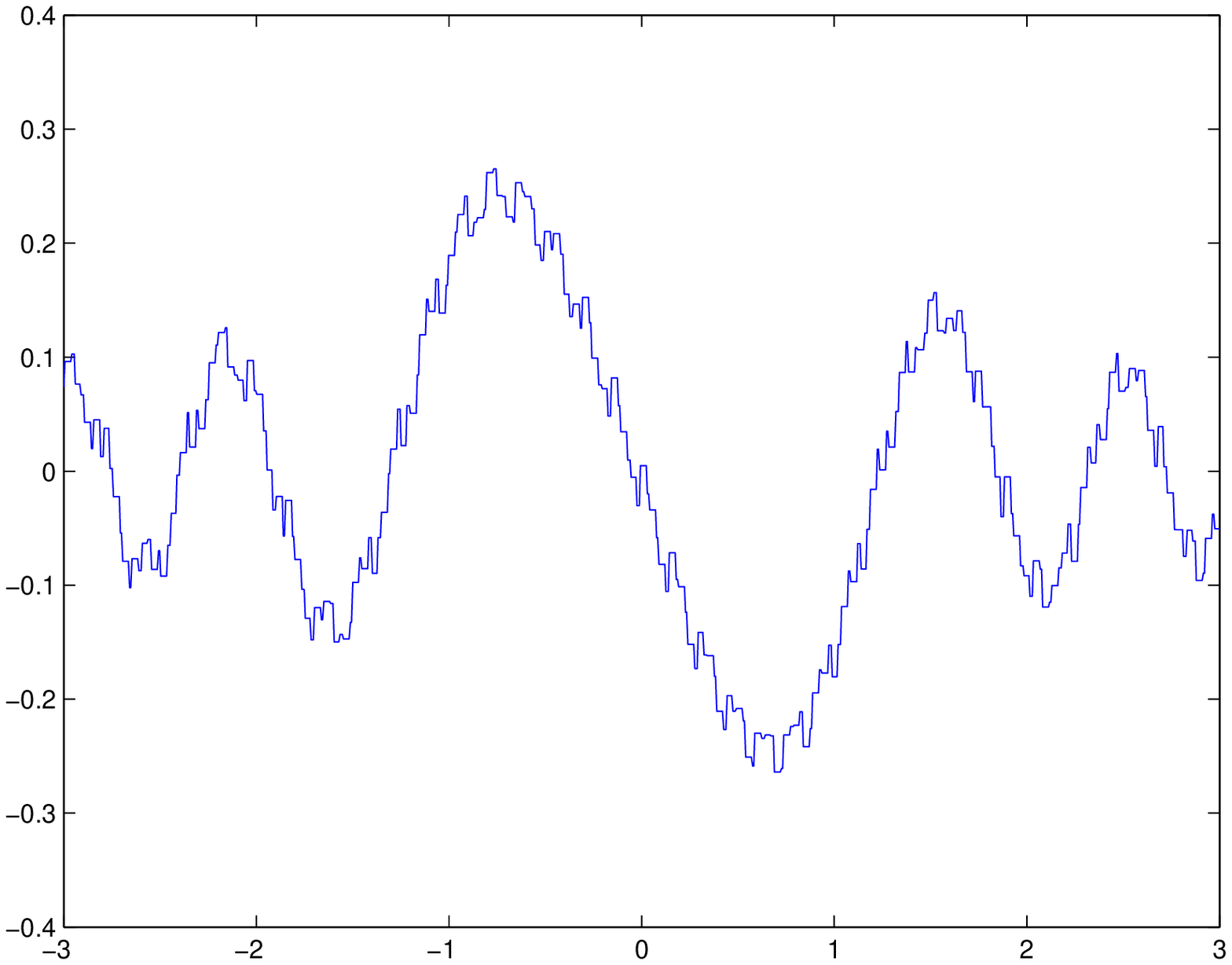}}
\end{center}
\begin{center}
\scalebox{0.35}{\includegraphics{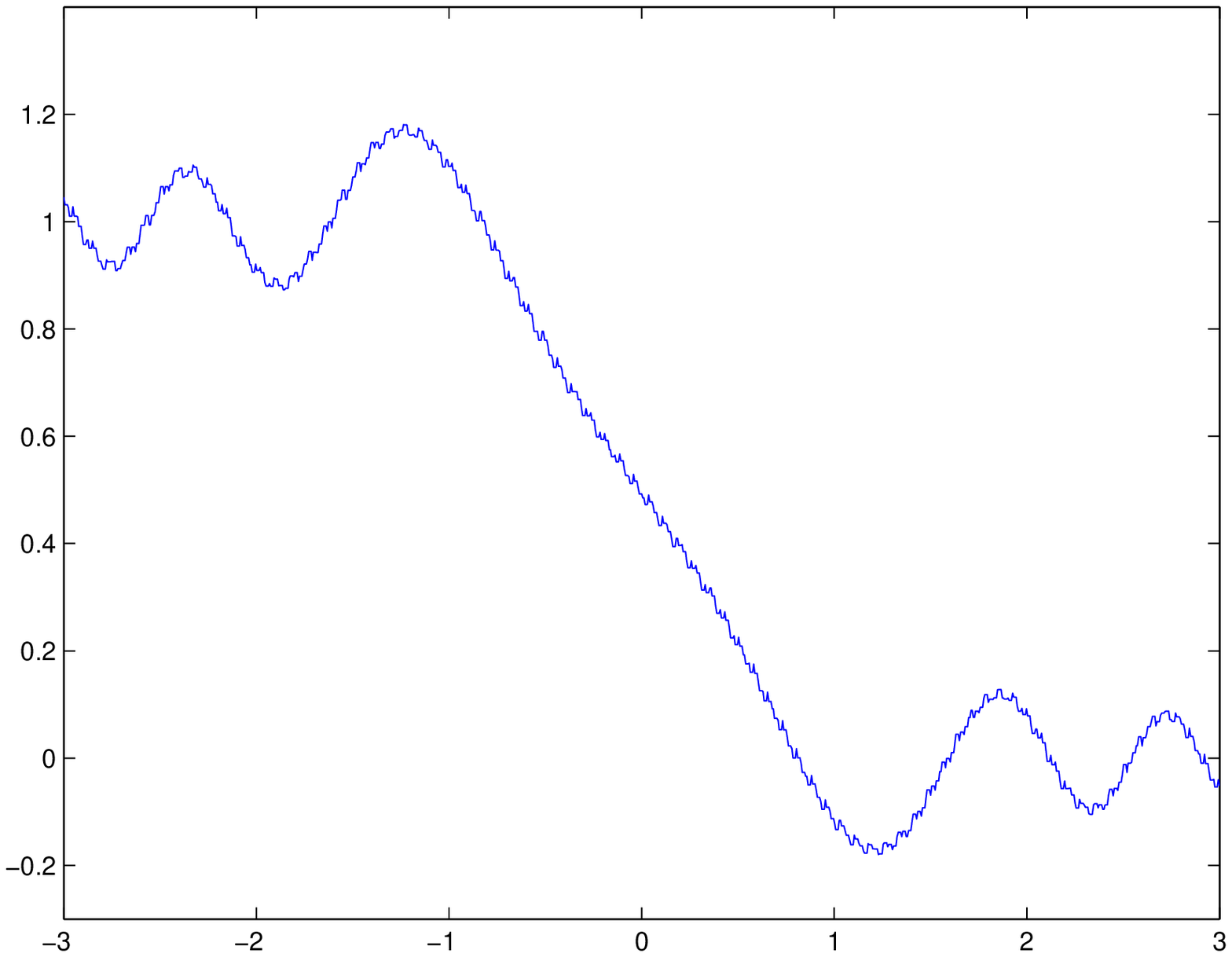}
\includegraphics{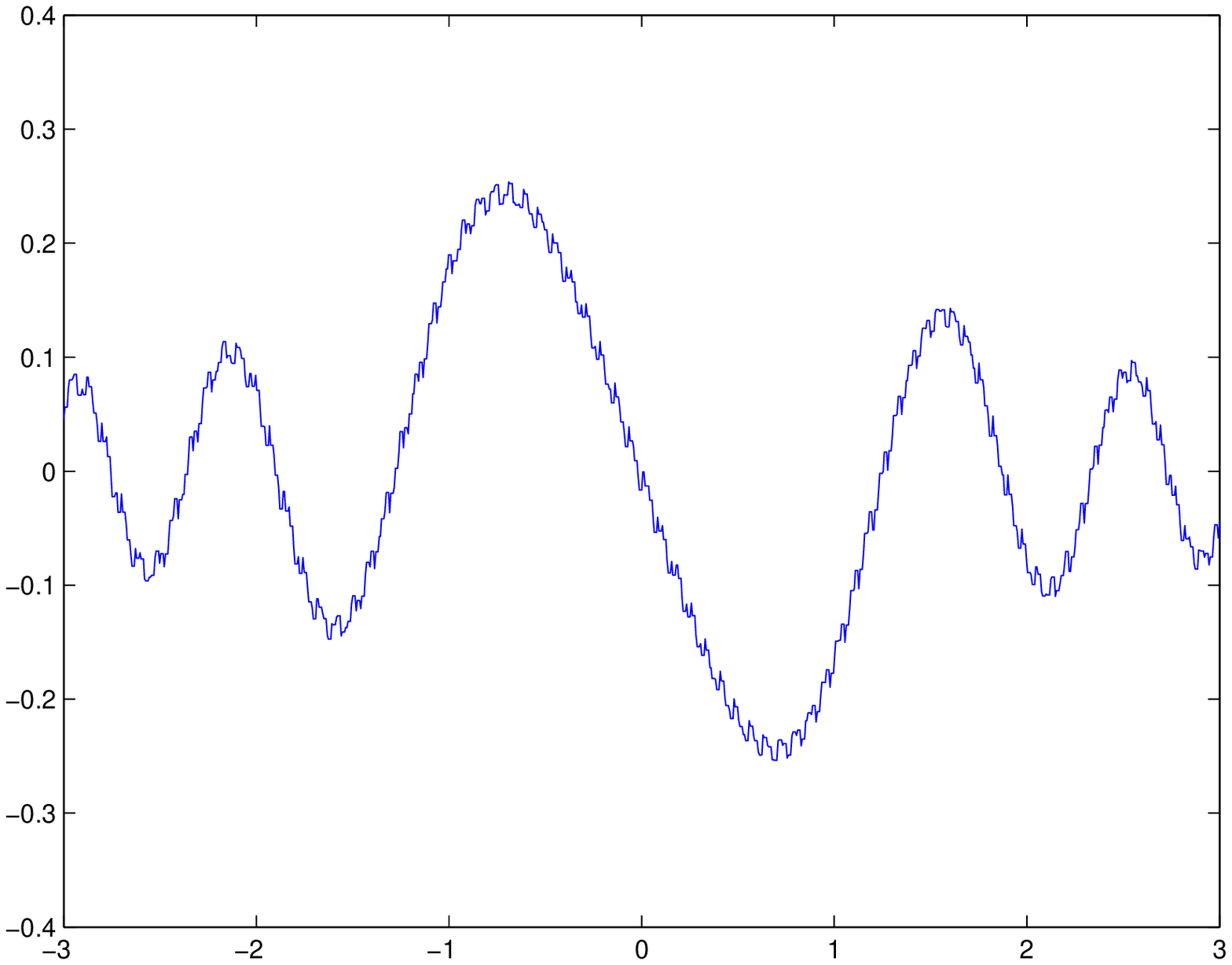}}
\end{center}
\begin{center}
\scalebox{0.35}{\includegraphics{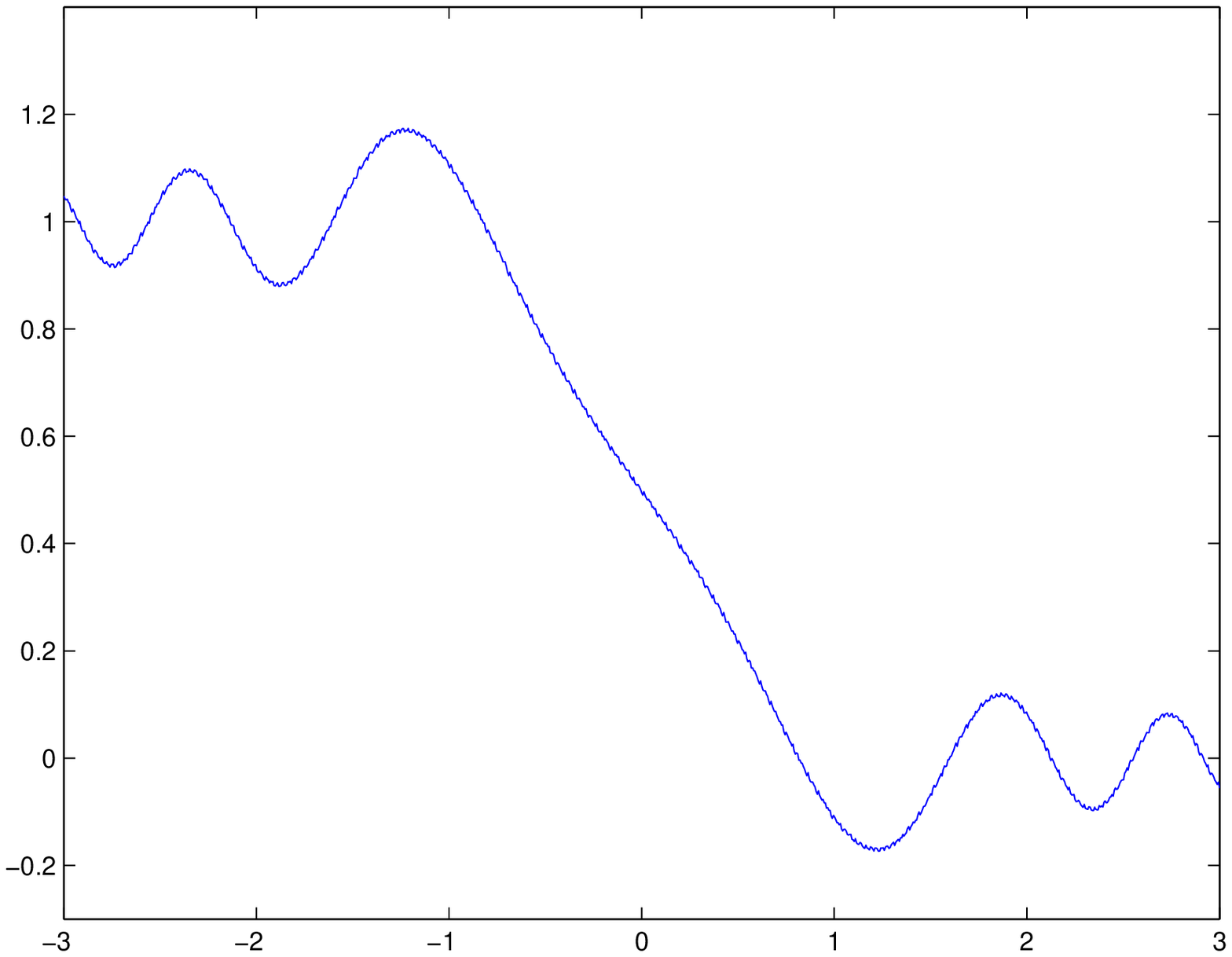}
\includegraphics{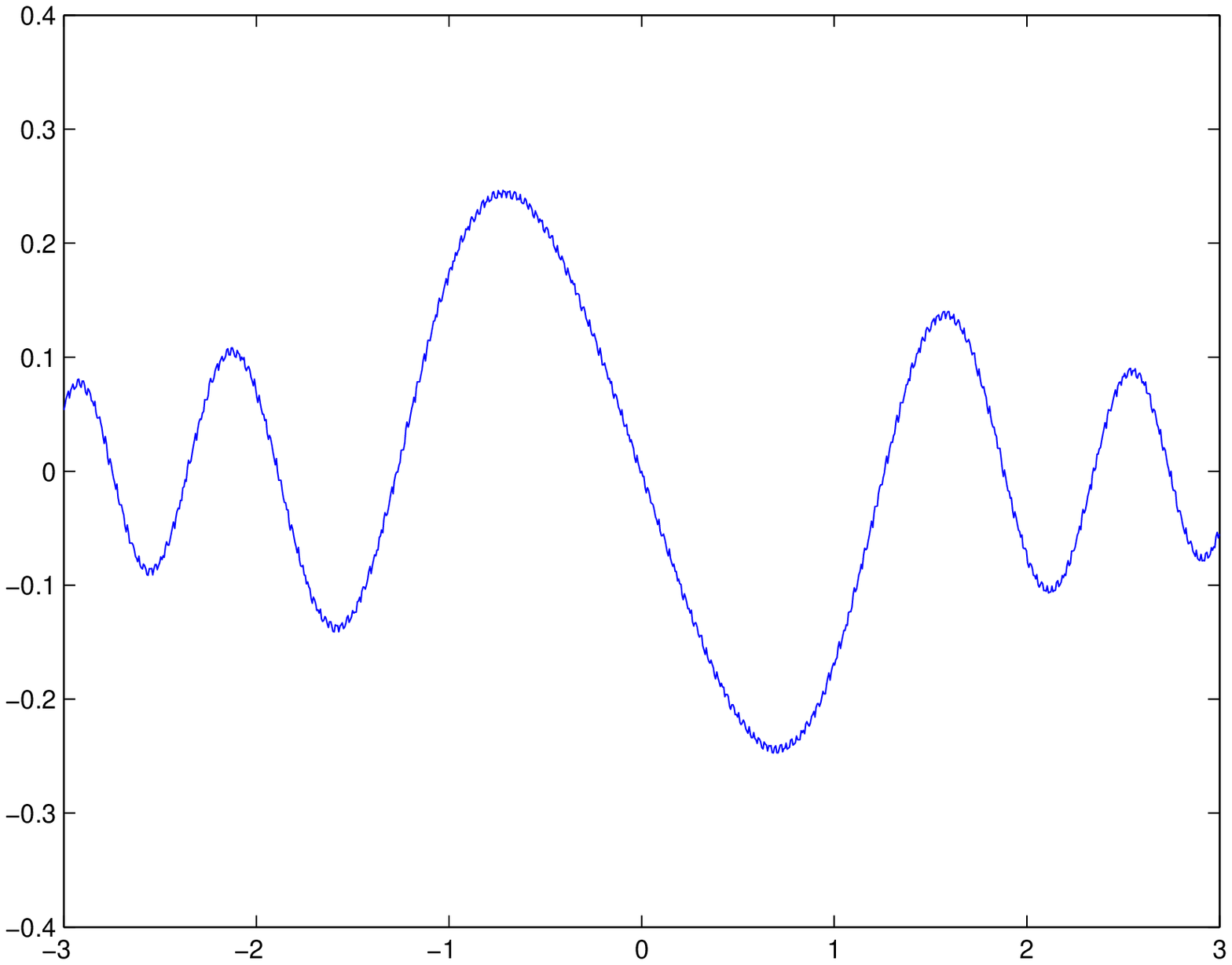}}
\end{center}
\vskip -5mm
\caption{Real and imaginary parts of 
$U(t,\gamma+st^{1/2})$ for $t=1/202,1/1616,1/6464,1/51712$.}
\end{figure}

As mentioned above the proof of Theorem 1.5 uses well-known 
estimates for Weyl sums which are obtained in the literature 
using the Weyl shift method; we state the required results before
proceeding.  If 
$f(k)=\alpha_nk^n+\cdots+ \alpha_0$, $\alpha_n\neq 0$, 
is a polynomial of degree $n\ge 2$ 
with real coefficients, and $k$ varies over an interval of at most $\mu$ 
consecutive integers, then for any $\varepsilon > 0$
\begin{align}
\left|\sum e(f(k))\right|^N
&\ll \mu^{N-1} +\mu^{N-n}\sum_{1\le r_1,\ldots,r_{n-1}\le
\mu-1}\min \left\lbrace \mu, \frac{1}{2\lbrace n!\alpha_n r_1r_2\ldots 
r_{n-1}\rbrace} \right\rbrace \nonumber \\ 
&\ll \mu^{N-1} +\mu^{N-n+\varepsilon}\sum_{1\le m\le (\mu-1)^{n-1}}
\min \left\lbrace \mu, \frac{1}{2\lbrace n!\alpha_n m\rbrace}
\right\rbrace 
\end{align}
where $\lbrace x\rbrace$ denotes the distance from $x$ to the nearest
integer, $N=2^{n-1}$.  (Originally due to Weyl.  See Titchmarsh 1986).  Note also that 
with $\alpha_n=u/q$ we have
\begin{align}
\sum_{1\le m\le(\mu-1)^{n-1}}
\min \left\lbrace \mu, \frac{1}{2\lbrace n!\alpha_n m\rbrace}
\right\rbrace 
=&\sum_{v({\hbox{\smallrm mod}}\,q)}\min \left\lbrace \mu, 
\frac{q}{2v}\right\rbrace
\sum_{\substack{{\scriptstyle{1\le m\le (\mu-1)^{n-1}}} \\
{\scriptstyle{n!um\equiv v({\hbox{\smallrm mod}}\,q)}}}}
1 \nonumber \\
\ll& (\mu^{n-1}q^{-1}+1)(\mu+q\log q) 
\end{align}
and hence 
\begin{equation}
\left|\sum e(f(k))\right|^N
\ll \mu^{N-1} +\mu^{N-n+\varepsilon}(\mu^{n-1}q^{-1}+1)(\mu+q)\log q\ .
\end{equation}
Furthermore, if
$$\left|\alpha_n-\frac{h}{q}\right|\le \frac{1}{q^2}\ ,\ (h,q)=1\ ,$$

\begin{figure}[ht]
\begin{center}
\scalebox{0.7}{\includegraphics{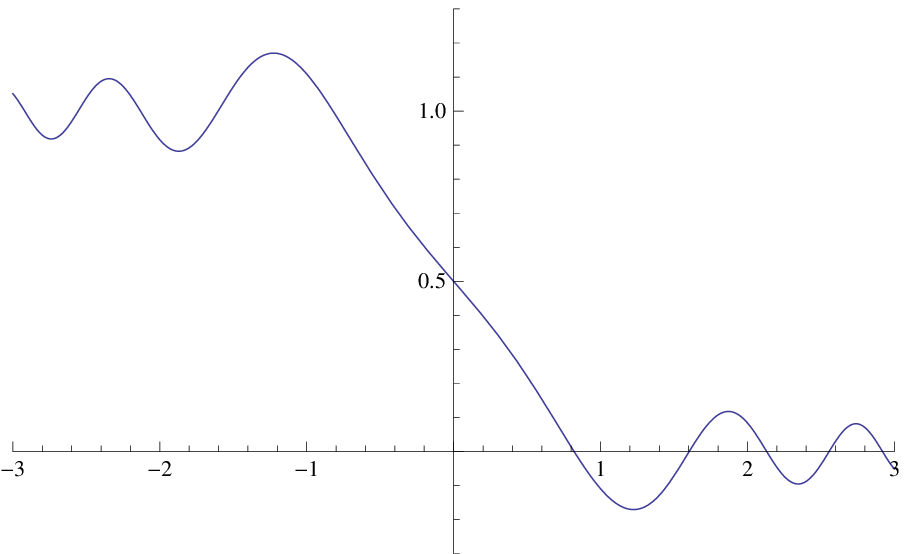}
\includegraphics{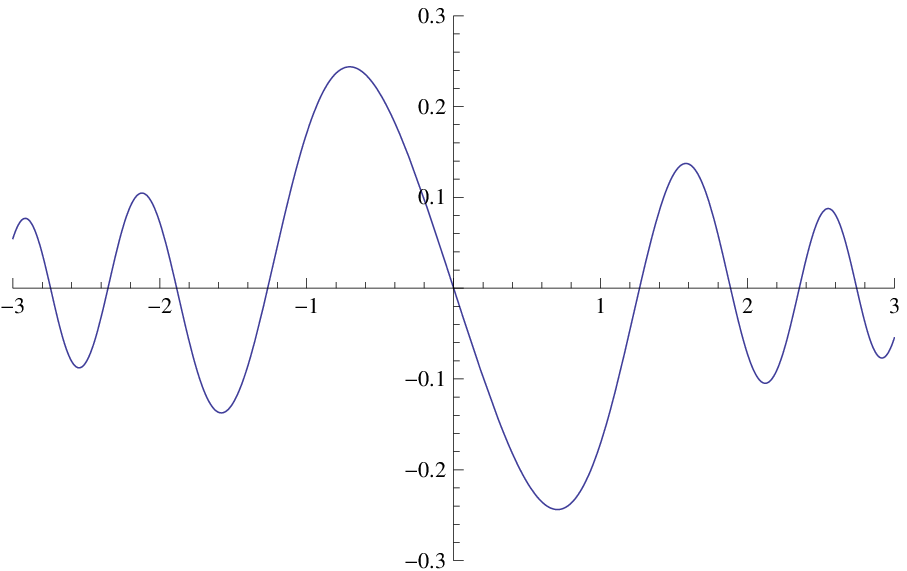}}
\end{center}
\vskip -6mm
\caption{The real and imaginary parts of the integral for $n=2$.}
\end{figure}
\begin{figure}[ht]
\begin{center}
\scalebox{0.7}{\includegraphics{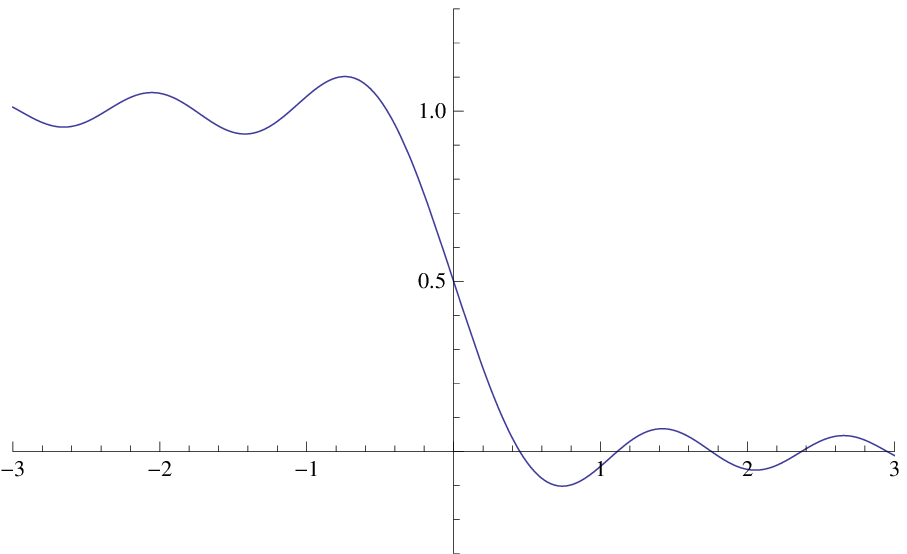}
\includegraphics{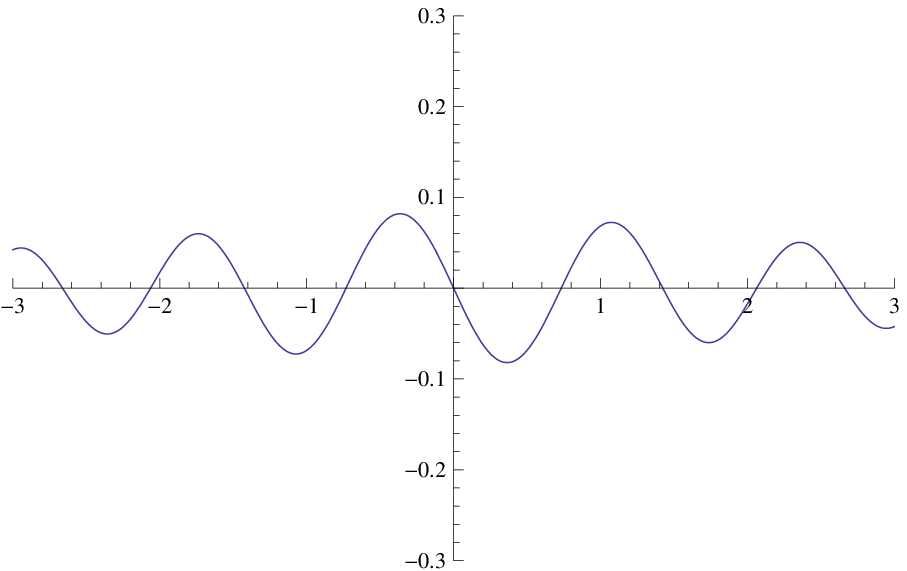}}
\end{center}
\vskip -6mm
\caption{The real and imaginary parts of the integral for $n=6$.}
\end{figure}

\begin{center}
\begin{figure}[ht]
\begin{center}
\scalebox{0.7}{\includegraphics{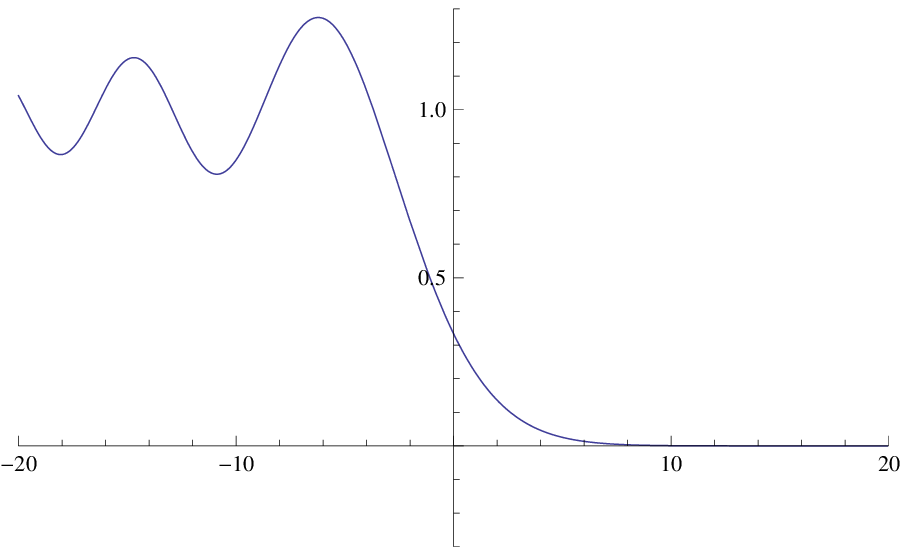}
\includegraphics{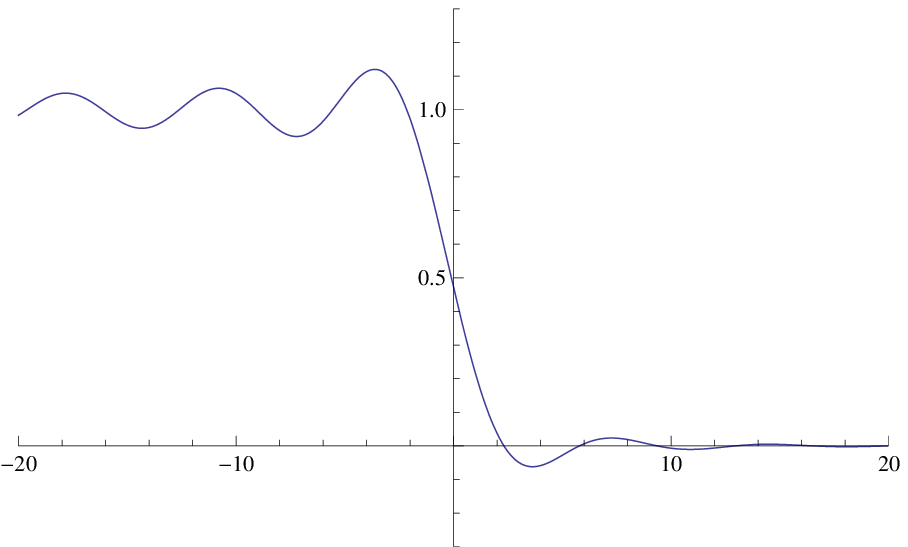}}
\end{center}
\vskip -2mm
\caption{The integral for $n=3$ and $n=17$.}
\end{figure}
\end{center}

\noindent then for any $\varepsilon > 0$,
\begin{equation}
\sum_{m=1}^M e(f(m))\ll \left\lbrace M^{-1}+q^{-1}+qM^{-n}\right\rbrace^{2^{1-n
}}M^{1+\varepsilon}q^{\varepsilon}
\end{equation}
where the implied constant depends on $\varepsilon$, and is uniform in 
the coefficients $\alpha_{n-1},\ldots,\alpha_0$. (See Hua 1965).  

\begin{center}
\begin{figure}[ht]
\begin{center}
\scalebox{0.5}{\includegraphics{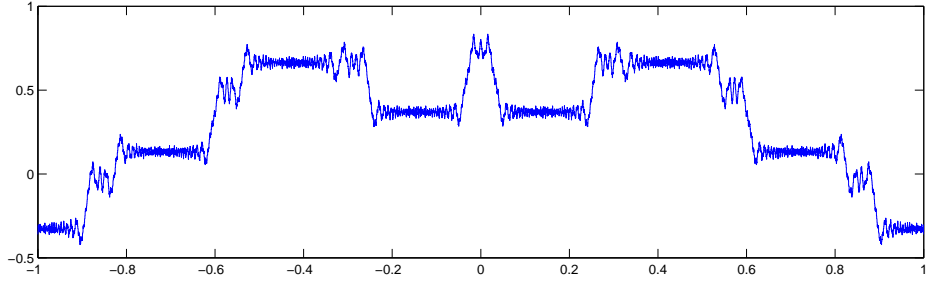}}
\end{center}
\vskip -5mm
\caption{The graph of $U(468/3277,x)$}
\end{figure}
\end{center}

The structure of the paper is as follows.  In section 2 we establish the 
properties of the sets ${\hbox{\script A}}$, ${\hbox{\script A}}_{m}$ 
and ${\hbox{\script B}}_{m,\alpha}$ stated in Lemma 1.2.  Some of these 
are not strictly necessary for the remainder of the paper, but item 6 
in particular is important, since it establishes that there are plenty 
of points in ${\hbox{\script B}}_{m,\alpha}$ through which $t$ may 
tend to zero.  In section 3 we prove Theorem 1.4, which handles the 
convergence of the Fourier series at rational times, and in section 4 
we use $(1.20)$ to consider convergence at irrational times
in ${\hbox{\script A}}$ and prove Theorem 1.3.  The proof of 
Theorem 1.5 follows in sections 5 to 7.  As is common with 
periodic problems, one would like to use harmonic analysis such as the 
Poisson summation formula to replace the series by an integral, which is 
more easily analysed.  In this case, however, it is necessary to first 
separate into large and small frequencies $k$, and treat the large 
frequencies separately, using quite different methods for rational times
and for irrational times in ${\hbox{\script B}}_{m,\alpha}$.
These sections also use $(1.19)$ and $(1.20)$, although the arguments 
are more involved than in the proofs of convergence, and require that 
irrational times be in the set ${\hbox{\script B}}_{m,\alpha}$ rather than
${\hbox{\script A}}$.  In section 7 we apply the Poisson summation formula to 
the small frequencies, and analyse the resulting integral, in essentially 
the same way for rational and irrational times, which completes the proof.

\vskip 2mm

\setcounter{section}{2}
\centerline{\bf 2. Proof of Lemma 1.2.}

\vskip 1mm

Items 1 and 2 are clear from the definition of the sets.
To prove item 3, consider $t\notin {\hbox{\script A}}_{m}$; thus  
there exists $M\ge m$ for which there is no approximant $u/q$ to $t$
with $M^{\Delta}< q\le M^{n-\Delta}$, and hence there must be
consecutive convergents $u_j/q_j$ and $u_{j+1}/q_{j+1}$ for which
$q_j\le M^{\Delta}$ and $q_{j+1}>M^{n-\Delta}$ for some $j$.  
Defining 
$${\hbox{\script S}}_{q,M}
=\left\lbrack 0, \frac{1}{qM^{n-\Delta}}\right\rbrack 
\cup \bigcup_{u=1}^{q-1} \left\lbrack 
\frac{u}{q}-\frac{1}{qM^{n-\Delta}},
\frac{u}{q}+\frac{1}{qM^{n-\Delta}}\right\rbrack
\cup \left\lbrack 1-\frac{1}{qM^{n-\Delta}},1\right\rbrack $$
for any $q$ and $M$, by $(1.9)$ we have 
$$t\in \left\lbrack 0, \frac{1}{q_j q_{j+1}}\right\rbrack 
\cup \bigcup_{u=1}^{q_j-1} \left\lbrack 
\frac{u}{q_j}-\frac{1}{q_jq_{j+1}},
\frac{u}{q_j}+\frac{1}{q_jq_{j+1}}\right\rbrack
\cup \left\lbrack 1-\frac{1}{q_jq_{j+1}},1\right\rbrack $$
and hence $t\in {\hbox{\script S}}_{q,M}$ for some $q,M$ with 
$q\le M^{\Delta}$ and $M \ge m$.  Thus
$$t\in \bigcup_{M\ge m}\bigcup_{q\le M^{\Delta}} {\hbox{\script S}}_{q,M}$$
and since $\mu\left({\hbox{\script S}}_{q,M}\right)=M^{-n+\Delta}$, we
can bound the measure of this union by
$$\mu\left(\bigcup_{M\ge m}\bigcup_{q\le M^{\Delta}} {\hbox{\script S}}_{q,M}
\right)
\le \sum_{M\ge m} M^{2\Delta-n} 
\ll m^{1+2\Delta-n}.$$
This proves item 3.

Suppose now that $t$ is an algebraic irrational, and recall $(1.8)$ 
and $(1.10)$.  Since all but finitely many convergents $u_j/q_j$
differ from $t$ by at least $q_j^{-2-\eta}$, but by at most 
$q_jq_{j+1}^{-1}$, it follows that for all but finitely many convergents
we have $q_{j+1} < q_j^{1+\eta}$.  Now any positive number greater than 
1 must fall between some $q_j$ and $q_{j+1}$, hence for any $M\ge 1$ 
there will exist $j$ for which $q_j < M^{n-\Delta} \le q_{j+1}< q_j^{1+\eta}$
and hence $M^{(n-\Delta)/(1+\eta)} \le q_j < M^{n-\Delta}$.
It follows that for all sufficiently large $M$ we can find a fraction 
$u_j/q_j$ satisfying the conditions of Theorem 1 
as long as we choose $\eta < n/\Delta -2$.  This proves item 4.

Item 5 is again clear from the definitions; item 6 can be proved by 
adapting the argument for item 3; if $0<t<t_0$ and 
$t\notin {\hbox{\script B}}_{m,\alpha}$ 
there must exist $M>m t^{-\alpha}$ such that there is no approximant 
$u/q$ to $t$ with $M^{\Delta}< q\le M^{n-\Delta}$.  Since there
is certainly an approximant with 
$$\left|t-\frac{u}{q}\right| < \frac{1}{qQ}$$
where $q\le Q$ and $Q$ is the integer part of $M^{n-\Delta}$, 
it must be that $q\le M^{\Delta}$.  Were $u=0$ we would have
$t<(qQ)^{-1}$, which implies $t<2/(m t^{-\alpha})^{n-\Delta}$ since 
$Q\ge M^{n-\Delta}/2$.  Thus $t<2m^{-n+\Delta} t^{\alpha(n-\Delta)}$ 
and hence $t^{\alpha(n-\Delta)-1} >m^{n-\Delta}/2$, which is clearly false, 
so $u\neq 0$.  Thus $t$ is in the set 
$${\hbox{\script T}}_{q,M}
\subseteq
\bigcup_{1\le u\le u_0} \left\lbrack 
\frac{u}{q}-\frac{2}{qM^{n-\Delta}},
\frac{u}{q}+\frac{2}{qM^{n-\Delta}}\right\rbrack$$
for some $M\ge m t^{-\alpha}> m t_0^{-\alpha}$ and 
$q\le M^{\Delta}$, where $u_0$ is such that $t_0< u_0/q+2/qM^{n-\Delta}$.
The smallest element of ${\hbox{\script T}}_{q,M}$ is at least $1/2q$, 
so $t>1/2q$ and hence $2qt_0>1$.  Thus we may choose $u_0$ to be
the integer part of $2qt_0$ plus one, which is no more than $4qt_0$.
The measure of ${\hbox{\script T}}_{q,M}$ is thus no more than
$16t_0/M^{n-\Delta}$, and the measure of the union is bounded by 
$$\mu
\left(\bigcup_{M\ge m t_0^{-\alpha}} \bigcup_{q\le M^{\Delta}}
{\hbox{\script T}}_{q,M}\right)
\le \sum_{M\ge m t_0^{-\alpha}} \sum_{q\le M^{\Delta}}
\frac{8t_0}{M^{n-\Delta}}
\ll t_0\left(m t_0^{-\alpha}\right)^{-n+2\Delta +1}.$$
This proves item 6, and completes the proof of Lemma 1.

\vskip 2mm

\setcounter{section}{3}
\setcounter{equation}{0}
\centerline{\bf 3. Convergence at rational times and proof of Theorem 1.4}

\vskip 1mm

This follows from Theorem 2, which we prove.  Fourier modes of the form 
$e(tP(k)+kx)$ solve the equation $LU=U_t$, so given 
$f\in { \hbox{\script D}}$, if $c_k$ are the Fourier coefficients
of the initial condition $U(0,x)=f(x)$ then the Fourier series 
$$\sum_k c_k e(tP(k)+kx)$$
is a solution to the initial value problem in the $L^2$ sense (as 
discussed in the introduction), and we wish to show it is convergent at 
rational times.  From Katznelson 2004, for instance, we know the Fourier 
series for $f(x)$ is convergent, so
$$\lim_{K\rightarrow\infty} 
\frac{1}{q}\sum_{v({\hbox{\smallrm mod}}\,q)}
G(u,v;q)
\sum_{|k|\le K}c_k e\left(\left(x+\frac{v}{q}\right)k\right)$$
exists and is equal to the right-hand side of $(1.15)$.
On the other hand, opening $G$ and changing orders of the two 
sums this is equal to
$$\lim_{K\rightarrow\infty} 
\frac{1}{q}
\sum_{|k|\le K}c_k 
\sum_{v,w({\hbox{\smallrm mod}}\,q)}
e\left(\frac{u}{q}P(w)-vw+\left(x+\frac{v}{q}\right)k\right)
=\lim_{K\rightarrow\infty} 
\sum_{|k|\le K}c_k 
e\left(\frac{u}{q}P(k)+xk\right)
$$
which hence exists, and proves both the convergence and Theorem 1.4.

\vskip 2mm

\setcounter{section}{4}
\setcounter{equation}{0}
\centerline{\bf 4. Convergence at irrational times and completion of the
proof of Theorem 1.3}

\vskip 1mm

Suppose now that $P(k)$ is a monic polynomial of degree $n\ge 2$ 
with real coefficients.  If $t\in {\hbox{\script A}}$ then there 
exists $m$ such that $t\in {\hbox{\script A}}_m$, and we consider $K>m$ 
and $L>K+m$, and sum by parts to obtain
\begin{equation}
\sum_{K<k\le L } \frac{e(tP(k)+xk)}{k}
=\frac{1}{L}\sum_{K<k\le L } e(tP(k)+xk)
+\int\limits_K^L \sum_{K<k\le \xi } e(tP(k)+xk)\,\frac{d\xi}{\xi^2}.
\end{equation}
The interval $K\le \xi\le K+m$ contributes at most $O(K^{-2})$ to the integral
(where the implied constant depends on $m$, and hence on $t$, although this is
unimportant in this section).  In the remaining range $\xi>K+m$ we 
rewrite the sum as
so that 
\begin{equation}
S:=\sum_{K<k\le \xi}e\left(tP(k)+xk\right)
=\sum_{k\le \xi-K}e\left(P^*(k)\right)
\end{equation}
where shifting $k$ has not affected the leading coefficient of $P$, 
so $P^*(k)$ is a polynomial of the form required for $(1.19)$ and 
$(1.20)$, with $\alpha_n=t$.  Choosing an approximant $u_j/q_j$ to $t$
with $(\xi-K)^{\Delta}<q_j \le (\xi-K)^{n-\Delta}$ we apply $(1.20)$ 
and find 
\begin{equation}
S\ll \left\lbrace (\xi-K)^{-1}+q_j^{-1}+q_j(\xi-K)^{-n}\right\rbrace^{2^{1-n
}}(\xi-K)^{1+\varepsilon}q_j^{\varepsilon}
\ll (\xi-K)^{1-2^{1-n}\Delta+(n+1)\varepsilon}.
\end{equation}
Specifying $\varepsilon =\Delta/2^n(n+1)$ this gives
$S\ll (\xi-K)^{1-\delta}$, where $\delta=\Delta/2^n$ and hence
$$\sum_{K<k\le L } \frac{e(tk^n+xk)}{k}
\ll L^{-\delta}
+\int\limits_K^L \xi^{-1-\delta}d\xi
\ll K^{-\delta} +L^{-\delta}\ ,$$
where we note that the constants are independent of $x$.  
It follows that the sequence of partial sums $(1.12)$ is a Cauchy 
sequence, and the series is thus convergent.  The uniformity follows 
since none of the constants appearing depend on $x$, and hence the 
tail of the series for $k>K$ is uniformly bounded by $O(K^{-\delta})$.

\vskip 2mm

\setcounter{section}{5}
\setcounter{equation}{0}
\centerline{\bf 5. Bounding the contribution from large frequencies at irrational times}

\vskip 2mm

Let $t\in {\hbox{\script B}}_{m,\alpha}$ as in Lemma 1.2,
and suppose $(n-\Delta)^{-1}<\alpha<(n-1)^{-1}$.  We define a 
${\hbox{\script C}}^{\,\,\infty}$ smoothing function $\phi(x)$
such that $\phi(x)=1$ for $x\in [-1/2,1/2]$, $\phi(x)=0$ for $x\notin [-2,2]$
and $\phi(x)+\phi(x^{-1})=1$ for all $x\in {\mathbb R}$, so  
\begin{equation}
U(t,x)=\frac{1}{\pi}\lim_{K\rightarrow\infty}\sum_{|k|\le K}
\!\left\lbrace \phi\left(kt^{\alpha}\right)
+\phi\left(\frac{1}{kt^{\alpha}}\right)\right\rbrace\!
\frac{\sin 2\pi\gamma k}{k} e(tk^n+xk) 
=U^*(t,x)+\lim_{K\rightarrow\infty}\! U_K(t,x)
\end{equation}
where
\begin{align}
U^*(t,x)&=\frac{1}{\pi}\sum_{k=-\infty}^{\infty}
\phi\left(kt^{\alpha}\right)
\frac{\sin \gamma k}{k} e(tk^n+xk)\ ,\\
U_K(t,x)&=\frac{1}{\pi}\sum_{|k|\le K}
\phi\left(\frac{1}{kt^{\alpha}}\right)
\frac{\sin 2\pi\gamma k}{k} e(tk^n+xk).
\end{align}
We can write $U_K(t,x)$ as the difference of the two
expressions
$$\frac{1}{2\pi i}
\sum_{|k|\le K}\phi\left(\frac{1}{kt^{\alpha}}\right)
\frac{1}{k} e(tk^n+(x\pm \gamma)k)$$
and will consider only positive values of $k$; the contribution
from negative values can be estimated by the same argument.  By summation 
by parts we obtain
\begin{align*}
&\frac{1}{2\pi i}\int\limits_{1/2t^{\alpha}}^{K}
\phi\left(\frac{1}{\xi t^{\alpha}}\right)
\frac{1}{\xi} d\left( \sum_{1/2t^{\alpha} < k \le \xi}e(tk^n+(x\pm \gamma)k)
\right) \\
=&\frac{1}{2\pi i} \phi\left(\frac{1}{Kt^{\alpha}}\right)
\frac{1}{K}\sum_{1/2t^{\alpha} < k \le K}e(tk^n+(x\pm \gamma)k) \\
&-\frac{1}{2\pi i}\int\limits_{1/2t^{\alpha}}^{K}
\frac{d}{d\xi}\left\lbrace
\frac{1}{\xi} \phi\left(\frac{1}{\xi t^{\alpha}}\right)\right\rbrace
\sum_{1/2t^{\alpha} < k \le \xi}e(tk^n+(x\pm \gamma)k) d\xi.\\
\end{align*}
We can bound the remaining sums over $k$ as
$$\sum_{1/2t^{\alpha} < k \le \xi}\!e(tk^n+(x\pm \gamma)k) d\xi
\ll \left\lbrace \xi^{-1}+q^{-1}+q\xi^{-n}
\right\rbrace^{1/N}\xi^{1+\varepsilon}q^{\varepsilon}$$
by $(1.20)$, where $N=2^{n-1}$.  
If $\xi > m t^{-\alpha}$ then we choose $q=q_j$ 
as in the definition of ${\hbox{\script B}}_{m,\alpha}$, 
and the sum is $O(\xi^{1-\Delta/N+\varepsilon})$; similarly 
for the first term with $K$ in place of $\xi$.  
On the other hand, if $\xi$ is smaller than this
we take the value of $q$ which the lemma gives for $M= (m +1)t^{-\alpha}$,
and the sum is $O\left(t^{-\alpha+\alpha \Delta/N-\varepsilon}\right)$.
On taking the limit as $K\rightarrow\infty$ we now have
\begin{equation}
U(t,x)= U^*(t,x)+O(t^{\alpha\Delta/N-\varepsilon}).
\end{equation}

\vskip 2mm

\setcounter{section}{6}
\setcounter{equation}{0}
\centerline{\bf 6. Bounding the contribution from large frequencies at rational times}

\vskip 1mm

We need to provide an argument for rational times similar to that proved 
in the previous paragraph for irrational times; specifically we will
consider $t=u/q\rightarrow 0^+$.
The analysis is a little 
trickier, however, so we begin by considering the related sum 
given by
$$V^*(t,x)=\frac{1}{\pi}\sum_{k=-\infty}^{\infty}
\phi\left(\frac{k}{q^4}\right)
\frac{\sin 2\pi \gamma k}{k} e(tk^n+xk)$$
where $\phi$ is as used in $(5.1)$.  We begin by showing 
\begin{equation}
U(t,x)=V^*(t,x)+O\left(q^{-1/2}\right).
\end{equation}
which allows us to consider the finite sum $V^*$ in place of $U$. 

By Fourier inversion 
\begin{equation}
V^*(t,x)=\frac{1}{\pi}\lim_{K\rightarrow\infty}
\int\limits_{-\infty}^{\infty}{\hat\phi}(\eta)
\sum_{|k|\le K}
\frac{\sin 2\pi \gamma k}{k} e\left(\frac{uk^n}{q}
+\left(x+\frac{\eta}{q^4}\right)k\right)d\eta
\end{equation}
where ${\hat\phi}$ denotes the Fourier transform $\hat\phi$,
which can be bounded using repeated integration by parts by
\begin{equation}
{\hat\phi}(\eta) = \int\limits_{-\infty}^{\infty}\phi(\xi)e(-\xi\eta)d\xi 
\ll \eta^{-\nu}\ ,\ \nu=0,1,2,\ldots\ .
\end{equation}
In order to change the orders of limit and integral in $(6.2)$ we note that
for any $y$
\begin{equation}
\frac{1}{\pi}\sum_{|k|\le K}
\frac{\sin 2\pi \gamma k}{k} e\left(tk^n+yk\right)
=\int\limits_{0}^1U(t,y+u)\frac{\sin \pi(2K+1)u}{\sin \pi u}du
\end{equation}
where $U(t,y+u)$ is piecewise constant, as discussed above, 
and consider the integral over $[0,1.2]$; the other half can be treated 
similarly.  The function $g(u)=\sin (\pi u) /\pi u$ is analytic
and bounded on $[0,\infty)$, and its antiderivative $G(u)$ 
is also bounded on $[0,\infty)$ as can be seen by observing that it tends
to a finite limit as $u\rightarrow \infty$.  Breaking up $(6.4)$ into
the integrals over subintervals where $U(t,y+u)$ is constant, if 
$[a,b]$ is any such subinterval then it contributes
$$(2K+1)\int\limits_{a}^bU(t,y+u)\frac{\pi u}{\sin \pi u}
g((2K+1)u)du.$$
Integrating by parts the integral in $(6.4)$ is $O(q)$, uniformly
in $K$ and $y$.  The Lebesgue Dominated Convergence Theorem
can then be applied in $(6.2)$ to change the orders of limit and integral, 
and obtain
\begin{equation}
V^*(t,x)
=\frac{1}{\pi}
\int\limits_{-\infty}^{\infty}{\hat\phi}(\eta)
\sum_{k=-\infty}^{\infty}
\frac{\sin 2\pi \gamma k}{k} e\left(\frac{uk^n}{q}
+\left(x+\frac{\eta}{q^4}\right)k\right)d\eta.
\end{equation}
Applying Theorem 1.4 to the sum over $k$ 
it can be evaluated as a value of $U(t,x)$, and 
\begin{equation}
V^*(t,x)
=\frac{1}{q}\sum_{w({\hbox{\smallrm mod}}\,q)}
\sum_{v({\hbox{\smallrm mod}}\,q)}e_q(uv^n-vw)
\int\limits_{-\infty}^{\infty}
{\hat\phi}(\eta)\,U\!\left(0,x+\frac{w}{q}+\frac{\eta}{q^4}\right)d\eta.
\end{equation}
Using the bound $(6.3)$ the integral over $\eta$ is
$$U\!\left(0,x+\frac{w}{q}\right)\int\limits_{-\infty}^{\infty}
\!\!{\hat\phi}(\eta)d\eta
+\int\limits_{|\eta|\le q^2}
\!\!{\hat\phi}(\eta)\left\lbrace U\!\left(0,x+\frac{w}{q}
+\frac{\eta}{q^4}\right)
-U\!\left(0,x+\frac{w}{q}\right)\right\rbrace d\eta
+O\!\left(q^{-A}\right)$$
for any $A> 0$, and applying this in $(6.6)$ the first term
gives $U(t,x)$ precisely, by Theorem 1.4.  The integrand in the 
second term vanishes unless one of the conditions 
$$-\gamma \le x+w/q \le \gamma\ ,\ -\gamma \le x+w/q+\eta/q^4 \le \gamma$$
is true and the other false; this can happen in four ways, each of which
implies that $w/q$ is in an interval of length $O(q^{-2})$, since 
the range of integration in $\eta$ is of length $q^2$.
For sufficiently large $q$ there can be at most one $w/q$ in each of these
intervals, hence these terms contribute $O\left(q^{-1/2}\right)$
to $(6.6)$.  (Here we have used the bound $(1.16)$ for the complete sum
over $v$ modulo $q$.)  This proves $(6.1)$.

Recalling that $\phi(kt^{\alpha})+\phi(1/kt^{\alpha})=1$, for sufficiently 
large $q$ we have 
$\phi(kt^{\alpha})\phi(k/q^4)=\phi(kt^{\alpha})$
and the definition of $V^*$ and $(6.1)$ imply that
$$U(t,x)=U^*(t,x)
+\frac{1}{\pi}
\sum_{k=-\infty}^{\infty} 
\phi\left(\frac{1}{kt^{\alpha}}\right)
\phi\left(\frac{k}{q^4}\right)
\frac{\sin 2\pi \gamma k}{k} e\left(\frac{u}{q}k^n+xk\right)
 +O(q^{-1/2})\ .$$
To estimate the second term it is sufficient to bound the sum
$$W(t,x)=\sum_{k=1}^{\infty} 
\phi\left(\frac{1}{kt^{\alpha}}\right)
\phi\left(\frac{k}{q^4}\right)k^{-1}
e\left(\frac{u}{q}k^n+yk\right)$$
where $y=x\pm \gamma$.  Summing by parts in a similar fashion to section 4,
\begin{align}
W(t,x)&=\int\limits_{1/2t^{\alpha}}^{2q^4}
\phi\left(\frac{1}{\xi t^{\alpha}}\right)
\phi\left(\frac{\xi}{q^4}\right)\xi^{-1}
d\left(\sum_{1/2t^{\alpha} < k \le \xi}
e\left(\frac{u}{q}k^n+yk\right)\right) \nonumber \\
&=-\int\limits_{1/2t^{\alpha}}^{2q^4}
\frac{d}{d\xi}\left\lbrace \phi\left(\frac{1}{\xi t^{\alpha}}\right)
\phi\left(\frac{\xi}{q^4}\right)\xi^{-1}\right\rbrace 
\sum_{1/2t^{\alpha} < k \le \xi}
e\left(\frac{u}{q}k^n+yk\right)\, d\xi
\end{align}
but rather than apply $(1.20)$ as we did for irrational $t$, 
we use $(1.19)$ and obtain
\begin{align*}
W(t,x)\ll &\int\limits_{1/2t^{\alpha}}^{2q^4}
\left\lbrace
\xi^{1-1/N} +\xi^{1-n/N+\varepsilon}
(\xi^{n/N-1/N}q^{-1/N}+1)(\xi^{1/N}+q^{1/N})\log q\right\rbrace
\xi^{-2}d\xi \\
\ll&\left\lbrace t^{\alpha/N}+q^{-1/N}
+t^{\alpha n/N}q^{1/N}\right\rbrace q^{\varepsilon} .
\end{align*}
Since $\alpha < 1/(n-1)\le 1$, we have $t^{\alpha}>t>q^{-1}$ 
so the second term can be dropped, as can the error term from $(6.1)$.
Furthermore the second term dominates the third only if 
$t<q^{-1/\alpha(n-1)}<q^{-1}$, which is false for sufficiently large $q$,
so $t^{\alpha}>t^{1/(n-1)}>q^{-1/(n-1)}$, and 
\begin{equation}
U(t,x)=U^*(t,x)+
O\left(
t^{\alpha n/N}q^{1/N+\varepsilon}
\right).
\end{equation}
Supposing now that $t\le mq^{-1/\alpha (n-\Delta)}$ this error term is 
$O(t^{\alpha\Delta/N-\varepsilon})$ as in $(5.4)$ and Theorem 1.5.

\vskip 2mm

\setcounter{section}{7}
\setcounter{equation}{0}
\centerline{\bf 7. An asymptotic estimate for $U^*(t,x)$}

\vskip 1mm

To estimate $U^*$ we write it as
$$U^*(t,x)=\int\limits_{x-\gamma}^{x+\gamma}
\sum_{k=-\infty}^{\infty}
\phi\left(kt^{\alpha}\right) e(tk^n+\eta k) d\eta$$
and consider the inner sum.  Note that for this part of the analysis 
it makes no difference whether $t$ is rational or irrational.  Applying 
the Poisson summation formula 
$$\sum_{k}\phi\left(kt^{\alpha}\right) e(tk^n+\eta k) 
=\sum_{m} \int\limits_{-\infty}^{\infty}
\phi\left(yt^{\alpha}\right) e(ty^n+(\eta-m)y)dy 
=t^{-\alpha}
\sum_{m} 
\int\limits_{-\infty}^{\infty}
\phi(y)e(F_m(y))dy$$
where $F_m(y)=t^{1-n\alpha}y^n+(\eta-m)t^{-\alpha}y$.
For $m\neq 0$ 
$$\left|\eta -m\right| \ge |m| -(|x|+\gamma) \ge
1 -(|x|+\gamma) > 0$$
so for sufficiently small $t$
$$\left|F_m'(y)\right| 
\ge |m-\eta|t^{-\alpha} 
-n t^{1-n\alpha}|y|^{n-1}
\ge \frac{1}{2}|m-\eta|t^{-\alpha}
\ge \frac{1}{2} \left(|m| -(|x|+\gamma)\right)t^{-\alpha}$$
since $\alpha < 1/(n-1)$.
Thus for $m\neq 0$, $F_m'(y)$ does 
not vanish in the region of integration and we can integrate by parts 
twice and find
$$\int\limits_{-\infty}^{\infty}
\!\phi(y)e(F_m(y))dy 
=-\frac{1}{4\pi^2}\!\int\limits_{-\infty}^{\infty}\!
\frac{d}{dy}\!\left\lbrace \frac{1}{F_m'(y)}\frac{d}{dy}
\!\left\lbrace \frac{\phi(y)}{F_m'(y)}\right\rbrace \right\rbrace 
e(F_m(y))dy \ll \frac{t^{2\alpha}}{(|m|-(|x|+\gamma))^{2}}$$
so the contribution from all $m\neq 0$ is thus
\begin{equation}
t^{-\alpha}\int\limits_{x-\gamma}^{x+\gamma}
\sum_{m\neq 0} \int\limits_{-\infty}^{\infty}
\phi(y)e(F_m(y))dy
\ll t^{\alpha}
\sum_{m=1}^{\infty}  
\left(|m|-(|x|+\gamma)\right)^{-2} \ll t^{\alpha},
\end{equation}
where the implied constant depends on $\gamma$ but not on $x$, and hence
\begin{align}
U^*(t,x)&=t^{-\alpha}\int\limits_{x-\gamma}^{x+\gamma}
\int\limits_{-\infty}^{\infty}
\phi(y)e(t^{1-n\alpha}y^n+\eta t^{-\alpha}y)dy\, d\eta
+O\left(t^{\alpha}\right) \nonumber \\
&=\int\limits_{(x-\gamma)t^{-1/n}}^{(x+\gamma)t^{-1/n}}
\int\limits_{-\infty}^{\infty}
\phi(yt^{\alpha-1/n})e(y^n+\eta y)dy\, d\eta
+O\left(t^{\alpha}\right). 
\end{align}

If $n$ is even, the equation $ny^{n-1}+\eta=0$ 
has a unique solution $y_0=(-\eta/n)^{1/(n-1)}$; if $n$ is odd
there at most two solutions, given by $\pm y_0$ if present.  
In either case the solutions are in the range of support of 
$\phi(yt^{\alpha-1/n})$ if and only if $|\eta|\le n2^{n-1}t^{1-1/n-\alpha(n-1)}$.
If $x > \gamma$ then for sufficiently small $t$ we have 
$(x-\gamma)t^{-1/n}>n2^{n}t^{1-1/n-\alpha(n-1)}$, since 
$\alpha < 1/(n-1)$, so there are no solutions $y_0$ 
in the support of $\phi(yt^{\alpha-1/n})$, and
$$|ny^{n-1}+\eta|\ge \eta -ny^{n-1} > 
\eta -n2^{n-1}t^{1-1/n-\alpha (n-1)}
> \eta -\frac{x-\gamma}{2}t^{-1/n}
>\eta/2.$$ 
Integrating by parts twice in $y$ we can now bound the integral in $(7.2)$ by
$$-\frac{1}{4\pi^2}
\int\limits_{(x-\gamma)t^{-1/n}}^{(x+\gamma)t^{-1/n}}
\int\limits_{-\infty}^{\infty}
\frac{d}{dy}\left\lbrace \frac{1}{ny^{n-1}+\eta}
\frac{d}{dy}\left\lbrace \frac{\phi(yt^{\alpha -1/n})}{ny^{n-1}+\eta}
\right\rbrace\right\rbrace
e(y^n+\eta y)dy\, d\eta
\ll t^{\alpha}(x-\gamma)^{-1}\ .$$
A similar calculation shows that if $x<-\gamma$ then the
integral is $O(t^{\alpha} |\gamma +x|^{-1})$; we can combine these
two cases with $(7.1)$ to give $U^*(t,x)\ll t^{\alpha }|x^2-\gamma^2|^{-1}$
where the implied constant depends on $\gamma$. 

On the other hand, if $-\gamma < x < \gamma$ then this same 
reasoning applies to the two tails $\eta > (x+\gamma)t^{-1/n}$
and $\eta < (x-\gamma)t^{-1/n}$ so that in this case
$$U^*(t,x)=
\int\limits_{-\infty}^{\infty}\int\limits_{-\infty}^{\infty}
\phi(yt^{\alpha-1/n})e(y^n+\eta y)dy\, d\eta
+O\left(t^{\alpha}|x^2-\gamma^2|^{-1}\right)\ .$$
The inner integral in the remaining expression is the Fourier
transform of $\phi(yt^{\alpha -1/n})e(y^n)$, hence the outer integral 
is inverting the transform and gives the value of the original function
at $0$, so combining these estimates with $(7.1)$ we have
$$U^*(t,x)=
\begin{cases}
1 +O\left(t^{\alpha}|x^2-\gamma^2|^{-1}\right) &\text{if $|x|<\gamma$} \\
O\left(t^{\alpha}|x^2-\gamma^2|^{-1}\right) &\text{if $|x|>\gamma$} 
\end{cases}$$
and combining this last with $(5.4)$ and $(6.8)$ we obtain the first 
claim in Theorem 1.5.

The transition between these two cases is the main concern in Theorem 1.5.
We will consider $x$ near $\gamma$, omitting the details for $x$ near
$-\gamma$, which are very similar.
If we renormalise values of $x$ as $x=\gamma+st^{1/n}$
for $s$ in some fixed interval $[-S,S]$, then $(7.2)$ becomes
$$U^*\left(t,\gamma+st^{1/n}\right)
=\int\limits_{s}^{s+2\gamma t^{-1/n}}
\int\limits_{-\infty}^{\infty}
\phi(yt^{\alpha-1/n})e(y^n+\eta y)dy\, d\eta
+O\left(t^{\alpha}\right).$$
Since $1/n < \alpha <1/(n-1)$ by hypothesis, $-1 < (1-\alpha n)(n-1) < 0$, 
and hence for $\eta> s+2\gamma t^{-1/n}$ 
$$\left| ny^{n-1}+\eta\right| \gg t^{-1/n}-t^{(1-\alpha n)(n-1)/n} \gg t^{-1/n}.$$
Integrating by parts twice as above, 
$$\int\limits_{s+2\gamma t^{-1/n}}^{\infty}
\int\limits_{-\infty}^{\infty}
\phi(yt^{\alpha-1/n})e(y^n+\eta y)dy\, d\eta \ll t^{\alpha}$$
and hence
\begin{equation*}
U^*\left(t,\gamma+st^{1/n}\right)
=\int\limits_{s}^{\infty}
\int\limits_{-\infty}^{\infty}
\phi(yt^{\alpha-1/n})e(y^n+\eta y)dy\, d\eta
+O\left(t^{\alpha}\right)\ .
\end{equation*}

Suppose $n$ is even.  Were the integral in $\eta$ over $(-\infty,\infty)$ 
then the double integral would be the Fourier inverse at zero of the 
Fourier transform of $\phi(yt^{\alpha-1/n})e(y^n)$, so the double 
integral would be equal to $1$.  Since the integral over 
$\eta\in (0,\infty)$ is half this, 
\begin{align*}
U^*\left(t,\gamma+st^{1/n}\right)
&=\frac{1}{2}-\int\limits_{0}^{s}
\int\limits_{-\infty}^{\infty}
\phi(yt^{\alpha-1/n})e(y^n+\eta y)dy\, d\eta
+O\left(t^{\alpha}\right) \\
&=\frac{1}{2}-\frac{1}{2\pi}\int\limits_{-\infty}^{\infty}
\phi(yt^{\alpha-1/n})e(y^n)\sin 2\pi sy\, \frac{dy}{y}
+O\left(t^{\alpha}\right).
\end{align*}
The calculation is completed by noting that by integration by parts 
$$\int\limits_{0}^{\infty}
\phi(1/yt^{\alpha-1/n})e(y^n+sy)\, \frac{dy}{y}
=-\int\limits_{0}^{\infty}
\frac{d}{dy} \left\lbrace \frac{\phi(1/yt^{\alpha-1/n})}{y \left( ny^{n-1}+s\right)}\right\rbrace e(y^n+sy)\, dy
\ll t^{\alpha n-1}$$
and similarly for $y^n-sy$, and hence 
$$U^*\left(t,\gamma+st^{1/n}\right)
=\frac{1}{2}-\frac{1}{\pi}\int\limits_{0}^{\infty}
e(y^n)\sin 2\pi sy\, \frac{dy}{y}
+O\left(t^{\alpha n-1}\right).$$

On the other hand, if $n$ is odd then $(7.2)$ becomes
\begin{align*}
U^*\left(t,\gamma+st^{1/n}\right)
&=\int\limits_{s}^{\infty}
\int\limits_{-\infty}^{\infty}
\phi(yt^{\alpha-1/n})\cos 2\pi(y^n+\eta y)dy\, d\eta
+O\left(t^{\alpha}\right) \\
&=\frac{1}{2}-\int\limits_0^{s}
\int\limits_{-\infty}^{\infty}
\phi(yt^{\alpha-1/n})\cos 2\pi(y^n+\eta y)dy\, d\eta
+O\left(t^{\alpha}\right) \\
&=\frac{1}{2}-\frac{1}{2\pi}
\int\limits_{-\infty}^{\infty}
\phi(yt^{\alpha-1/n})\sin 2\pi(y^n+\eta y)\frac{dy}{y}
+O\left(t^{\alpha}\right) 
\end{align*}
which gives the form in Theorem 1.5.

\vskip 3mm

\vbox{\baselineskip=9pt
\noindent {\smallrm K.D.T.-R.M. was supported in part by US National Science 
Foundation grant number DMS-0800979, and in part by CNPq grant 312363/2009-5.  
He gratefully acknowledges the assistance and support of the faculty and staff
of the Universidade de Bras\'{i}lia, and also of the Mathematical Sciences 
Research Institute, where parts of this work were carried out.  

\noindent N.J.E.P. 
was supported in part by FAPDF - the Research Foundation of  the Federal 
District, Brazil.  This paper was completed during the Spring and Summer of 
2011, when he was a member at the Institute for Advanced Study; he thanks 
the Institute and its faculty and staff for their assistance and support. }}

\vskip 7mm

\centerline{\bf References}

\vskip 2mm

\noindent Abramovitz, M \& Stegun, I. 1972, Handbook of Mathematical Functions, US Govt. Printing Office.

\vskip 2mm
\noindent Arkhipov, G. I. \& Oskolkov, K. I. 1989, On a special 
trigonometric series and its applications, {\it Math. USSR Sbornik}, {\bf 62} (1), 145-155

\vskip 2mm
\noindent Berry, M. V. \& Klein, S. 1996,  Integer, fractional and fractal Talbot effects, {\it Journal of Modern Optics} {\bf 43}, no. 10, 2139--2164.

\vskip 2mm
\noindent DiFranco, J. C. \& McLaughlin, K. T.-R., 2005, A nonlinear Gibbs-type phenomenon for the defocusing nonlinear Schr\"odinger equation, {\it Int. Math. Res. Pap.} {\bf 8}, 403-459

\vskip 2mm
\noindent Hardy, G. H. \& E. M. Wright, E. M. 1979, {\it An introduction to the theory of numbers}, 5th edition, OUP.

\vskip 2mm
\noindent Hua, L.-K. 1965, {\it Additive theory of prime numbers}, AMS.

\vskip 2mm
\noindent Hua, L.-K. 1982, {\it Introduction to number theory}, Springer. 

\vskip 2mm
\noindent Kapitanski, L \& Rodnianski, I 1999, Does a quantum particle know the time?, {\it Proceedings of the workshop on Emerging Applications of 
Number Theory, IMA Volumes in Mathematics and its Applications} 109, 355-371.

\vskip 2mm
\noindent Katznelson, Y 2004, {\it An introduction to harmonic analysis}, 3rd edition, CUP

\vskip 2mm
\noindent Olver, P 2010, Dispersive quantization, {\it Amer. Math. Monthly} {\bf 117}, 599-610

\vskip 2mm
\noindent Oskolkov, K. I. 1992, A class of I. M. Vinogradov’s 
series and its applications in harmonic analysis, Springer Series in 
Computational Mathematics, 19, {\it Progress in Approxi mation Theory, an 
International perspective} (A. A. Gonchar and E. B. Sta, eds.), 
Springer-Verlag, New York.

\vskip 2mm
\noindent Rodnianski, I. 2000, Fractal solutions of the Schr\"{o}dinger equation, in: Nonlinear PDE's, dynamics and continuum physics (South Hadley, MA, 1998), 181Ð187, {\it Contemp. Math}. {\bf 255}, Amer. Math. Soc., Providence, RI, 2000. 

\vskip 2mm
\noindent Rodnianski, I. 1999, Continued fractions and Schrödinger 
evolution, Proceedings of the conference on Continued Fractions: from 
Analytic Number Theory to Constructive Approximations (eds. B. Berndt, 
F. Gesztesy), {\it Contemp. Math}. 236 (1999), 311-323.

\vskip 2mm
\noindent Roth, K. F. 1955, Rational approximations to algebraic numbers, {\it Mathematika} {\bf 2} 1-20; corrigendum, 168

\vskip 2mm
\noindent Schmidt, W. 2004, {\it Equations over finite fields}, 2nd edition, 
Kendrick Press

\vskip 2mm
\noindent Stein, E.M. and Wainger, S. 1990, Discrete analogues of 
singular Radon transforms, {\it Bull. Amer. Math. Soc.} 23 , 537-544.

\vskip 2mm
\noindent Talbot, H. F. 1836, Facts relating to Optical Science. 
No. IV, {\it Phil. Mag.} {\bf 9}, 401-407.

\vskip 2mm
\noindent Titchmarsh, E. C. 1986, {\it The theory of the Riemann zeta-function}, 2nd edition, OUP.

\end{document}